\documentclass[11pt]{amsart}
\usepackage[margin=1.15in]{geometry}
\usepackage{comment}

\usepackage[T1]{fontenc}
\usepackage{lmodern,amsfonts,amsmath,amstext,amsbsy,amssymb,amsopn,amsthm,upref,eucal,mathptmx,mathtools,url,thmtools
}
\usepackage[colorlinks,linkcolor=blue,citecolor=blue]{hyperref}

\newtheorem{thm}{Theorem}[section]
\newtheorem{lem}[thm]{Lemma}
\newtheorem{prop}[thm]{Proposition}
\newtheorem{cor}[thm]{Corollary}

\theoremstyle{definition}
\newtheorem{defn}[thm]{Definition}
\newtheorem{oss}[thm]{Remark}
\numberwithin{equation}{section}

\newcommand{\bR}{\mathbb R}
\newcommand{\bT}{\mathbb T}
\newcommand{\cP}{\mathcal P}
\newcommand{\cR}{\mathcal R}
\newcommand{\cC}{\mathcal C}
\newcommand{\cA}{\mathcal A}
\newcommand{\cS}{\mathcal S}
\newcommand{\Reg}{\mathrm{Reg}}
\newcommand{\Leb}{\mathrm{Leb}}
\newcommand{\bZ}{\mathbb Z}

\newcommand{\cU}{\mathcal U}
\newcommand{\cE}{\mathcal E}
\newcommand{\sL}{\mathsf L}
\newcommand{\Dd}{\mathfrak D}

\begin{document}

\title[Mean-field limit for a thermostatted Lorentz gas]{Mean-field limit and Vlasov equation for a Gaussian-thermostatted Lorentz gas}
\author{Roberto Castorrini}
\address{Dipartimento di Economia, Ingegneria, Società e Impresa (DEIM), Università della Tuscia, 01100, Viterbo, Italy}
\email{roberto.castorrini@unitus.it}
 

\begin{abstract}
We study $N$ particles moving in a dispersing billiard under an electric field and coupled
through a Gaussian thermostat. The corresponding mean-field equation is a Vlasov-type
equation with specular reflection. We reformulate the equation as a fixed
point problem for the time-integrated thermostat coefficient. This yields existence of weak solutions for
every absolutely continuous initial datum of finite energy and, for bounded initial densities supported on
velocities that are large compared with the field, uniqueness of the weak solution; under the finite
horizon assumption, the upper bound on the velocities can be replaced by an exponential tail of the initial density. In these regimes, we prove convergence of the empirical measures of the particle system to the
solution of the kinetic equation. In general, convergence holds along subsequences. 
\end{abstract}

\maketitle
 \tableofcontents

\section{Introduction}

The derivation of kinetic equations from microscopic particle systems is a central problem in
mathematical physics. In the mean-field regime, a large class of interacting systems is expected to be
described, as the number of particles $N\to\infty$, by nonlinear transport equations of Vlasov type.
This scenario has been rigorously established in various settings, starting
from the classical works of Braun--Hepp \cite{BraHe} and Dobrushin
\cite{Dob} (see also \cite{Spo,Gol1} for the general theory), and
\cite{GoMoRi} for the relation between empirical measures and Vlasov
hierarchies. Important extensions of the classical theory to rough or
singular interaction forces were obtained in \cite{HaJa,JaWa,Ser}. A distinct
line of work derives Vlasov equations from fermionic quantum many-body dynamics in
a joint mean-field and semiclassical regime, see, e.g.,
\cite{BePoSaSc} and references therein. Taken together, the above works support the general principle that a microscopic
system with mean-field interactions should, in the many-particle limit, be
governed by a macroscopic equation of Vlasov type.

In the present paper we
investigate this principle for a classical mechanical system whose
one-particle dynamics is a billiard flow. The resulting specular collisions
and grazing trajectories introduce singularities that constitute the main novelty of the problem.\\
We study a thermostatted system of particles evolving in a dispersing billiard (Sinai
billiard) under the action of a constant\footnote{The constant field assumption can be dropped and, in some cases, with no major complications. See Section \ref{sec:further}.} external electric field. The dynamics is that of a Lorentz gas with
specular reflections at the obstacles, coupled through a Gaussian thermostat which enforces conservation
of the total kinetic energy.  The case of a single moving
particle is a canonical model of deterministic nonequilibrium statistical
mechanics. For sufficiently weak fields, it admits a physical
nonequilibrium steady state (NESS) with well-understood statistical
properties; in particular, its transport, response and entropy-production
mechanisms have been studied extensively
\cite{BoCh1,BuSi,Chernov2001,Chernov2008,ChernovZhang2015,DeReZh}.
In the many particles case, although the particles don't interact directly, the friction
coefficient depends on their empirical current and therefore couples all
their motions through a mean field while preserving the total kinetic
energy. Formal analysis and numerical evidence point to a self-consistent
kinetic evolution and to a stable small-field distribution of particle
speeds \cite{BoDaLeRi}. However, a rigorous derivation of this macroscopic evolution
from the deterministic many-particle billiard dynamics has remained unavailable.


From the mathematical point of view, the thermostat leads to a nonlinear Vlasov-type equation for the
one-particle distribution, while the obstacles introduce specular reflections and grazing collisions.
The billiard flow is only piecewise smooth, and its derivative is unbounded near grazing trajectories
\cite{CheMa}. As a consequence the Cauchy--Lipschitz theory does not apply to the characteristics, and the
classical Dobrushin stability argument \cite{Dob} is not available. There are results on Vlasov equations with
specular reflection in bounded domains, mainly for Vlasov--Poisson-type systems, such as the boundary
regularity results in \cite{GuoHalf}, the global existence results 
\cite{HwVe}, and the more recent global well-posedness result \cite{CeIa}. In these works the behaviour near grazing velocities plays a
crucial role, and conditions such as flatness of the datum near the grazing set are imposed, which do not applied in the present setting. The closest
rigorous analogue of our setting is the thermostatted kinetic model of \cite{BoCAEsLe}, where the Gaussian
thermostat also generates a mean-field interaction but collisions are modelled by a stochastic mechanism
instead of deterministic reflections. To the best of our knowledge, a rigorous mean-field limit for a
deterministic thermostatted Lorentz gas has not been established before.\\
Our approach is mainly dynamical and it is  based on the observation that the nonlinearity is of rank one, i.e. all particles feel the same
friction coefficient $a(t)$, and when the field $E_0$ is constant the dynamics depends on $a$ only through
its primitive $A(t)=\int_0^ta$. For a given {history} $A$ the one-particle dynamics is a
non-autonomous billiard flow $\Phi_A^t$, with explicit affine free flights, and it turns out taht the kinetic equation is
equivalent to the scalar fixed point problem
\[
A(t)=\int\langle E_0,\widetilde X^A_t(z)-x\rangle\,d\mu_0(z),
\]
where $\widetilde X^A_t(z)-x$ is the displacement of the particle started at $z=(x,v)$, measured in the
universal cover.\footnote{This formulation also has a  physical interpretation: denoting $V_s^A(z)$ the velocity at time
$s$ of the trajectory generated by $A$, since
\[
\langle E_0,\widetilde X_t^A(z)-x\rangle
=
\int_0^t\langle E_0,V_s^A(z)\rangle\,ds,
\]
the right-hand side is the average work performed by the external field up
to time $t$; under the normalizing the kinetic energy, self-consistency
identifies it with the thermostat coefficient $A(t)$.}

The right-hand side involves positions only, which are continuous in time, so that no trace
theory at the boundary is needed. The $N$-particle system satisfies exactly the same identity with $\mu_0$
replaced by the initial empirical measure.\\
The results follow the classical scheme of a mean-field limit.
\begin{itemize}
\item[--] \emph{Particle system} (Theorem~\ref{thm:micro}). For almost every initial configuration the
$N$-particle dynamics is well defined and unique, and its empirical measure $\delta_{N,t}$ is a weak solution
of the kinetic equation.
\item[--] \emph{Kinetic equation} (Theorem~\ref{thm:macro}). For every absolutely continuous initial datum
$\mu_0$ with unit kinetic energy there exists a weak solution $\mu_t$. If $\mu_0$ has a bounded density supported in
$\{c_0\le|v|\le L\}$ and $T$ is such that all speeds remain large compared with
$\sqrt{|E_0|/\kappa_{\min}}$, where $\kappa_{\min}$ is the minimal curvature of the obstacles, the weak solution
is unique among absolutely continuous weak solutions. Under finite horizon (see Sec.\ref{sec:nobound}), the upper bound $L$ can be replaced by an exponential tail of the initial density (see Remark \ref{rem:tail}) compatible with the predictions of \cite{BoCh2}.
\item[--] \emph{Mean-field limit} (Theorem~\ref{thm:mfl}). In the regime of uniqueness, if the initial empirical
measures converge to $\mu_0$ then $\delta_{N,t}$ converges to $\mu_t$ for every $t$. Without the
uniqueness assumption, the empirical measures are relatively compact and all their limit points are weak
solutions with initial datum $\mu_0$.
\end{itemize}
The idea of the proof in the regime $\{c_0\le|v|\le L\}$ is actually naive, its implementation is nevertheless quite technical. It suffices to compare two histories on a short time
window after the last time at which they agree. In the fast regime each particle has at most one
collision in such a window, and for a single collision the loss of regularity is a factor $1/\nu$, where
$\nu$ is the normal velocity at the collision. 
It turns out that this singularity is balanced exactly because the Lebesgue measure of the initial conditions is transported
to the flux measure,\footnote{It measures the rate at which phase points cross the boundary.} which has the factor $\nu$. Hence $1/\nu$ is integrable, with an explicit bound, and
the nearly grazing collisions gives a higher order term. Thanks to this observation, we can prove the result in this regime 
without requiring hyperbolicity,  transversality of
singularity manifolds or flatness of the datum near the grazing set.\\
Without the upper bound on the velocities, a very fast particle has many collisions in a window, and the
singular factors of the collisions multiply and the above argument fails. However, if the table has finite horizon, we show that the Lebesgue average of
their product grows at most exponentially in the number of collisions, uniformly in the history, by following
one-parameter families of initial data and counting the pieces into which the collisions cut them. An exponential tail of the initial density, compatible with previous numerical results, compensates this growth (Theorem~\ref{thm:C-tail}).

The main question left open is global uniqueness, that is, uniqueness without the lower cut-off on the
speeds.\footnote{We stress that the existence and mean-field result are independent from this cut-off.} Slow particles may return several times to the same obstacle in a short time, hopping along its
boundary, and for such trajectories a different argument is needed (we discuss this in
Remark~\ref{rem:global-uniqueness}, where we explain why we expect uniqueness to hold in general). Other natural questions
are the extension to non-constant fields (we discuss this possibility in Section \ref{sec:further}), quantitative rates and propagation of chaos, and the long-time
behaviour and stationary states of the limit equation, in connection with the nonequilibrium steady states
studied in \cite{BoCAEsLe,BoCh1,BoCh2,BoCh3,BoDaLeRi}. 

Finally, the thermostat couples the particles through a global observable, which
places the model in the broader circle of ideas surrounding globally
coupled dynamical systems and self-consistent transfer operators, initiated by Keller \cite{Kel1,Kel2}. For
general results on self-consistent transfer operators see
\cite{GalSTO,TanziReview}; for hyperbolic globally coupled systems and
their statistical properties see \cite{BaLiSa}, and for bifurcations and
phase transitions in mean-field coupled chaotic systems see \cite{BaLi}.
Recent developments concerning the differential and stability of
self-consistent transfer operators are in
\cite{CasGaTa,CasGaTa1}. We develop this point of view in
Appendix~\ref{app:STO}. All of the aforementioned works treat the uniform hyperbolic case, while the non-uniformly hyperbolic case has proven to be essential in the understanding of physically-relevant models. In particular, an interesting connection to the derivation of macroscopic transport and stochastic behaviour from deterministic chaotic dynamics is the investigation of
fast--slow partially hyperbolic systems (see, for example,
\cite{CaLi,DeSiLi,LiTransport} and the references therein); a first evidence of that is the recent remarkable result \cite{CLO} in which the heat equation has been derived for a harmonic chain
under fast deterministic chaotic forcing, with no introduction of randomness. Those
works provide a complementary dynamical-systems perspective on the
emergence of effective transport laws, but they still concern smooth dynamics and we hope that the present paper can be a useful step forward to the study of more general discontinuous cases.

\medskip\noindent\textbf{Structure of the paper.}
The length of the paper is largely due to the fact that it is entirely self-contained, although the proofs of some lemmas follow a standard approach, particularly those relating to Sinai billiards.\\
Section~\ref{sec:model} introduces the particle system and the kinetic equation, and
Section~\ref{sec:results} collects the definitions and a summary of the results. In
Section~\ref{sec:flow} we construct the billiard flow for a given history and prove that almost every
point is regular. Section~\ref{sec:lagr} contains the fixed point formulation, the existence theorem and the equivalence
between weak and Lagrangian solutions,
Section~\ref{sec:mf} the mean-field limit, and Section~\ref{sec:fast} the uniqueness theorem in the fast
regime. In Sections~\ref{sec:nobound}--\ref{sec:fastcomp} we remove the upper bound on the velocities for tables with
finite horizon, and Section~\ref{sec:further} collects further questions. The paper concludes with two appendices: Appendix~\ref{app:transported-tests} is the technical proof of Lemma \ref{lem:invariant-test}; Appendix~\ref{app:STO} explores the
connection with self-consistent transfer operators.

\section{The model}\label{sec:model}

\subsection{Microscopic dynamics}
Let $\bT^2=\bR^2/\mathbb Z^2$ and ${\bf D}:=\bT^2\setminus\cup_{i=1}^rQ_i$, where $Q_1,\dots,Q_r$ are open,
strictly convex subsets of $\bT^2$ with pairwise disjoint closures and $\cC^3$ boundaries of positive
curvature. The one-particle phase space is $\Omega:={\bf D}\times\bR^2$. We consider $N$ particles of unit
mass with positions $x_i\in{\bf D}$ and velocities $v_i\in\bR^2$, subject to a constant field
$E_0\in\bR^2\setminus\{0\}$ and to the Gaussian thermostat. Setting
\[
J_N:=\sum_{j=1}^N\langle E_0,v_j\rangle,\qquad U_N:=\sum_{j=1}^N|v_j|^2,\qquad\alpha_N:=\frac{J_N}{U_N},
\]
the equations of motion in the interior of ${\bf D}$ are
\begin{equation}\label{eq:equation of motion}
\dot x_i=v_i,\qquad\dot v_i=E_0-\alpha_N\,v_i,\qquad i=1,\dots,N,
\end{equation}
and when a particle reaches $\partial{\bf D}$ at a point $x$ its velocity is reflected specularly,
\begin{equation}\label{eq:elastic-coll}
v^+=\cR_x(v^-),\qquad\cR_x(v):=v-2\langle v,n(x)\rangle\,n(x),
\end{equation}
where $n(x)$ is the unit normal at $x\in\partial{\bf D}$ pointing into ${\bf D}$. A direct computation gives
$\tfrac{d}{dt}U_N=2J_N-2\alpha_NU_N=0$, and reflections preserve $|v_i|$, so the total kinetic energy is
conserved. The flow is not Hamiltonian: the divergence of the vector field in
\eqref{eq:equation of motion} is $-(2N-1)\alpha_N$ (it is however conformally symplectic, see \cite{LiMi}).
The collisions with the obstacles do not correlate the particles and the interaction is due to the
thermostat alone.

We fix the energy per particle to be one, i.e.\ we work on the invariant shell
\begin{equation}\label{eq:shell-intro}
\frac1N\sum_{i=1}^N|v_i|^2=1 .
\end{equation}
This is the normalisation compatible with a nontrivial limit $N\to\infty$ in which the empirical measure
$\frac1N\sum_i\delta_{(x_i,v_i)}$ converges to a probability measure with finite second moment in velocity.
On the shell \eqref{eq:shell-intro} the thermostat coefficient is an empirical average,
$\alpha_N=\frac1N\sum_j\langle E_0,v_j\rangle$, and \eqref{eq:equation of motion} reads
\begin{equation}\label{eq:equation of motion2}
\dot x_i=v_i,\qquad\dot v_i=E_0-\Bigl[\frac1N\sum_{j=1}^N\langle E_0,v_j\rangle\Bigr]v_i,\qquad i=1,\dots,N.
\end{equation}

\subsection{The kinetic equation}
Equation \eqref{eq:equation of motion2} has the mean-field form
$$\dot z_i=\frac1N\sum_jV(z_i,z_j)$$ with $$V((x,v),(x',v'))=(v,E_0-\langle E_0,v'\rangle v).$$ Replacing the
empirical average by an average with respect to the law $f_t$ of one particle leads formally to the
Vlasov-type equation
\begin{equation}\label{eq:Vlasov}
\partial_tf+\operatorname{div}_{(x,v)}\bigl(V_f\,f\bigr)=0,\qquad
V_f(t,x,v)=\Bigl(v,\ E_0-\frac{j_f(t)}{u_f(t)}\,v\Bigr),
\end{equation}
\[
j_f(t):=\int_\Omega\langle E_0,v\rangle f(t,x,v)\,dx\,dv,\qquad u_f(t):=\int_\Omega|v|^2f(t,x,v)\,dx\,dv,
\]
with initial datum $f(0,\cdot)=f_0$ and the specular boundary condition
\begin{equation}\label{eq:condizione al bordo}
f(t,x,v)=f\bigl(t,x,\cR_x(v)\bigr),\qquad x\in\partial{\bf D} .
\end{equation}
Mass and kinetic energy are conserved. For a smooth solution decaying fast enough in $v$, multiply
\eqref{eq:Vlasov} by $\psi(v)=1$ or $\psi(v)=|v|^2$ and integrate over $\Omega$. The boundary term
$\int_{\partial{\bf D}\times\bR^2}\langle v,n(x)\rangle\,\psi(v)\,f\,d\sigma(x)\,dv$ vanishes, because the change of
variables $v\mapsto\cR_xv$ preserves $\psi$, $dv$ and, by \eqref{eq:condizione al bordo}, $f$, while it changes the
sign of $\langle v,n(x)\rangle$. The remaining term vanishes for $\psi=1$, and for $\psi=|v|^2$ it gives
\[
\frac{d}{dt}u_f=2j_f-2\,\frac{j_f}{u_f}\,u_f=0 .
\]
Hence, for $u_f(0)=1$, the friction coefficient is $j_f(t)$, the macroscopic counterpart of $\alpha_N$. For weak
solutions this is Proposition~\ref{prop:conservation}. Solutions will in general not be smooth, even for smooth data, because of the
grazing trajectories, we therefore work with weak solutions and with their Lagrangian description.

\section{Definitions and results}\label{sec:results}

We denote by $\cP(\Omega)$ the set of Borel probability measures on $\Omega$. A sequence $\mu_n\in\cP(\Omega)$
\emph{converges narrowly} to $\mu$, written $\mu_n\rightharpoonup\mu$, if $\int\varphi\,d\mu_n\to\int\varphi\,d\mu$
for every $\varphi\in\cC_b(\Omega)$, the bounded continuous functions. A function $\varphi$ on $\Omega$ (or on
$[0,T)\times\Omega$) is \emph{symmetric} if $\varphi(x,v)=\varphi(x,\cR_xv)$ for $x\in\partial{\bf D}$.

\begin{defn}[Weak solution]\label{def:w-sol}
Let $\mathcal D_{\mathrm{sym}}$ be the set of symmetric functions in $\cC_c^1([0,T)\times\Omega)$.
A family $\{\mu_t\}_{t\in[0,T]}\subset\cP(\Omega)$ is a \emph{weak solution} of
\eqref{eq:Vlasov}--\eqref{eq:condizione al bordo} with initial datum $\mu_0$ if:
\begin{enumerate}
\item[(i)] for every $\varphi\in\cC_b(\Omega)$ the function $t\mapsto\int\varphi\,d\mu_t$ is Borel measurable, and it
is continuous whenever $\varphi$ is symmetric;
\item[(ii)] $0<\inf_t\int|v|^2d\mu_t\le\sup_t\int|v|^2d\mu_t<\infty$;
\item[(iii)] for every $\varphi\in\mathcal D_{\mathrm{sym}}$,
\begin{equation}\label{eq:weak-form-final}
\begin{gathered}
\int_0^T\!\!\int_\Omega\Bigl(\partial_t\varphi+\langle\nabla_z\varphi,V_{\mu_t}\rangle\Bigr)d\mu_t\,dt
+\int_\Omega\varphi(0,z)\,d\mu_0(z)=0,\\
V_{\mu_t}(x,v):=\Bigl(v,\ E_0-\frac{\int\langle E_0,v'\rangle d\mu_t}{\int|v'|^2d\mu_t}\,v\Bigr).
\end{gathered}
\end{equation}
\end{enumerate}
We denote by $\mathcal W_T(\mu_0)$ the set of weak solutions with initial datum $\mu_0$ such that $\mu_t$ is
absolutely continuous for every $t$.
\end{defn}

\begin{oss}\label{rem:sym-cont}
Continuity is required only against symmetric functions because the velocity of a particle jumps at a
collision: the empirical measure of the particle system (see the scheme \eqref{eq:mfscheme} below) is not narrowly continuous in time, but it is
continuous in the sense of {\rm(i)}. The symmetry of the test functions in {\rm(iii)} encodes the specular
boundary condition \eqref{eq:condizione al bordo}.
\end{oss}

\begin{prop}[Conservation of mass and energy]\label{prop:conservation}
Let $\{\mu_t\}_{t\in[0,T]}$ be a family of finite Borel measures on $\Omega$ satisfying {\rm(i)--(iii)} of
Definition~\ref{def:w-sol}. Then, for every $t\in[0,T]$,
\[
\mu_t(\Omega)=\mu_0(\Omega),\qquad\int_\Omega|v|^2\,d\mu_t=\int_\Omega|v|^2\,d\mu_0 .
\]
In particular, if $\{\mu_t\}$ is a weak solution with $\int|v|^2d\mu_0=1$, then $\int|v|^2d\mu_t=1$ for every $t$,
the function $j_\mu(t):=\int_\Omega\langle E_0,v\rangle\,d\mu_t$ is Borel measurable with $|j_\mu|\le|E_0|$, and
$V_{\mu_t}(x,v)=(v,\,E_0-j_\mu(t)\,v)$.
\end{prop}

\begin{proof}
Let $M(t):=\mu_t(\Omega)$, $u(t):=\int|v|^2d\mu_t$ and $j_\mu$ as above. By {\rm(i)} applied to the symmetric
function $1$, $M$ is continuous, hence bounded; by {\rm(ii)}, $u$ is bounded above and below by positive
constants; by {\rm(i)} and dominated convergence $j_\mu$ is Borel, and $|j_\mu|\le|E_0|M^{1/2}u^{1/2}$ by the
Cauchy--Schwarz inequality. Hence the friction coefficient $a:=j_\mu/u$ of $V_{\mu_t}$ is a bounded Borel function
of $t$.

Let $\chi\in\cC^\infty([0,\infty))$ be nonincreasing, with $\chi=1$ on $[0,1]$ and $\chi=0$ on $[2,\infty)$, and for
$R\ge1$ let $\chi_R(v):=\chi(|v|/R)$. Let $\psi$ be either $\chi_R$ or $|v|^2\chi_R$. Since $\psi$ depends on $|v|$ only
and $|\cR_xv|=|v|$, for every $\zeta\in\cC^1_c([0,T))$ the function $\zeta(t)\psi(v)$ belongs to $\mathcal D_{\mathrm{sym}}$
(recall that ${\bf D}$ is compact). Its gradient has no $x$-component, so \eqref{eq:weak-form-final} gives
\[
\int_0^T\zeta'(t)\,m_\psi(t)\,dt+\zeta(0)\,m_\psi(0)+\int_0^T\zeta(t)\,g_\psi(t)\,dt=0,
\]
where $m_\psi(t):=\int\psi\,d\mu_t$ and $g_\psi(t):=\int\langle\nabla\psi(v),E_0-a(t)v\rangle\,d\mu_t$. The function
$m_\psi$ is continuous by {\rm(i)}, and $g_\psi$ is bounded. With $G_\psi(t):=m_\psi(0)+\int_0^tg_\psi(s)\,ds$, the
identity reads $\int_0^T\zeta'(m_\psi-G_\psi)\,dt=0$ for every $\zeta\in\cC^1_c([0,T))$; hence $m_\psi-G_\psi$ is almost
everywhere equal to a constant, which is $0$ (take $\zeta(0)=1$). Both functions being continuous,
\begin{equation}\label{eq:m-psi}
m_\psi(t)=m_\psi(0)+\int_0^tg_\psi(s)\,ds\qquad\text{for every }t\in[0,T].
\end{equation}
Now $\nabla\chi_R(v)=R^{-1}\chi'(|v|/R)\,v/|v|$ vanishes outside $\{R\le|v|\le2R\}$, so that, for $R\ge1$,
\[
\bigl|\langle\nabla\chi_R(v),E_0-a v\rangle\bigr|\le\|\chi'\|_\infty\bigl(|E_0|+2|a|\bigr)\,\mathbf 1_{\{|v|\ge R\}} .
\]
For $\psi=\chi_R$ this gives $|g_\psi(s)|\le\|\chi'\|_\infty(|E_0|+2|a(s)|)\,\mu_s(\{|v|\ge R\})$, which tends to $0$ for
every $s$ as $R\to\infty$ and is bounded uniformly in $R$ and $s$. Letting $R\to\infty$ in \eqref{eq:m-psi}, by
monotone convergence on the left and dominated convergence on the right, we get $M(t)=M(0)$. For
$\psi=|v|^2\chi_R$,
\[
\langle\nabla\psi,E_0-a v\rangle=2\chi_R(v)\bigl(\langle E_0,v\rangle-a|v|^2\bigr)+|v|^2\langle\nabla\chi_R(v),E_0-a v\rangle .
\]
As $R\to\infty$ the first term converges to $2(\langle E_0,v\rangle-a|v|^2)$ and the second to $0$, and both are
bounded by $C(1+|v|^2)$ with $C$ independent of $R$ and $s$. Hence $g_\psi(s)\to2j_\mu(s)-2a(s)u(s)=0$ for every
$s$, with $|g_\psi(s)|\le C(M(s)+u(s))$, while $m_\psi(t)\uparrow u(t)$ because $\chi_R$ increases with $R$. Letting
$R\to\infty$ in \eqref{eq:m-psi} gives $u(t)=u(0)$.

If $\{\mu_t\}$ is a weak solution with $u(0)=1$, then $M\equiv1$ and $u\equiv1$, so that $a=j_\mu$, and
$|j_\mu|\le|E_0|$.
\end{proof}

The results of the paper are organised as in the following diagram, in which the horizontal arrows are
the two evolutions and the vertical arrows are the limit $N\to\infty$:
\begin{equation}\label{eq:mfscheme}
\begin{array}{ccc}
\delta_{N,0}=\dfrac1N\displaystyle\sum_{i=1}^N\delta_{z^N_i} & \xrightarrow{\ \text{particle dynamics}\ } &
\delta_{N,t}=\dfrac1N\displaystyle\sum_{i=1}^N\delta_{z^N_i(t)}\\[4mm]
\Big\downarrow & & \Big\downarrow\\[2mm]
\mu_0 & \xrightarrow{\ \text{kinetic equation}\ } & \mu_t
\end{array}
\end{equation}

\begin{thm}[Particle system]\label{thm:micro}
Let $N\ge1$ and $T>0$. For almost every initial configuration on the shell \eqref{eq:shell-intro}, with
respect to its surface measure, the equations \eqref{eq:equation of motion2} with the reflection law
\eqref{eq:elastic-coll} have a unique solution on $[0,T]$, with finitely many collisions, each transversal
and involving a single particle. The empirical measures $\{\delta_{N,t}\}_{t\in[0,T]}$ of this solution form
a weak solution of the kinetic equation with initial datum $\delta_{N,0}$.
\end{thm}

\begin{thm}[Kinetic equation]\label{thm:macro}
Let $\mu_0=f_0\,dz\in\cP(\Omega)$ with $\int_\Omega|v|^2d\mu_0=1$, and let $T>0$.
\begin{enumerate}
\item[(1)] There exists a weak solution $\{\mu_t\}_{t\in[0,T]}$ with initial datum $\mu_0$. It is narrowly
continuous, absolutely continuous, and $\int|v|^2d\mu_t=1$ for every $t$.
\item[(2)] Assume in addition that, for some $0< c_0<L<\infty$, $f_0\in L^\infty$ is supported in $\{c_0\le|v|\le L\}$. Denoting by $\kappa_{\min}$ the minimal curvature of $\partial{\bf D}$,  if
\begin{equation}\label{eq:fast-intro}
(c_0+1)e^{-|E_0|T}>1\qquad\text{and}\qquad\kappa_{\min}\bigl((c_0+1)e^{-|E_0|T}-1\bigr)^2>2|E_0|,
\end{equation}
then $\mathcal W_T(\mu_0)$ consists of a
single element, that is,  the weak solution is unique among absolutely continuous weak solutions.
\item[(3)] Assume instead that the table has finite horizon {\rm(FH)} (Section~\ref{sec:nobound}), that
$f_0(x,v)\le C_0e^{-\alpha|v|}$ for some $C_0,\alpha>0$, that $f_0$ is supported in $\{|v|\ge c_0\}$, and that
\eqref{eq:fast-intro} holds. Then $\mathcal W_T(\mu_0)$ consists of a single element.
\end{enumerate}
\end{thm}

\begin{thm}[Mean-field limit]\label{thm:mfl}
Let $\mu_0$ and $T$ be as in Theorem~\ref{thm:macro}{\rm(2)} or {\rm(3)}, and let $\{\mu_t\}$ be the unique element of
$\mathcal W_T(\mu_0)$. For each $N$ let $(z^N_1,\dots,z^N_N)$ be an initial configuration as in
Theorem~\ref{thm:micro}, and assume that $\delta_{N,0}\rightharpoonup\mu_0$ as $N\to\infty$. Then
\[
\delta_{N,t}\rightharpoonup\mu_t\qquad\text{for every }t\in[0,T].
\]
If $\mu_0$ only satisfies the assumptions of Theorem~\ref{thm:macro}{\rm(1)}, then from every subsequence
one can extract a further subsequence along which $\delta_{N,t}$ converges narrowly, for every
$t\in[0,T]$, to a weak solution with initial datum $\mu_0$.
\end{thm}

\begin{oss}
For a given $c_0>\sqrt{2|E_0|/\kappa_{\min}}$, condition \eqref{eq:fast-intro} holds for all
$$T<T_*:=|E_0|^{-1}\log\left(\dfrac{1+c_0}{1+\sqrt{2|E_0|/\kappa_{\min}}}\right).$$ Since $\int|v|^2d\mu_0=1$ forces
$c_0\le1$, this is a condition of small field. For fixed $c_0$ and small field,
$T_*\sim|E_0|^{-1}\log(1+c_0)$, which is comparable to $c_0/|E_0|$, the time scale on which the thermostat
changes the velocities by an amount of order one; we have $T_*\to0$ as $c_0$ approaches the threshold
$\sqrt{2|E_0|/\kappa_{\min}}$. We expect uniqueness, and hence the
convergence of the whole sequence, to hold without condition \eqref{eq:fast-intro} (see Remark~\ref{rem:global-uniqueness}).
\end{oss}

The proofs go through a Lagrangian description. For a given \emph{history} $A$ (the time integral of
the friction coefficient) we construct in Section~\ref{sec:flow} the billiard flow $\Phi^t_A$, defined at
almost every point. A \emph{Lagrangian solution} (Definition~\ref{def:lagrangian}) is a curve
$\mu_t=(\Phi^t_A)_\#\mu_0$ whose history is the one it generates. Lagrangian solutions are weak solutions
(Proposition~\ref{prop:lagr-weak}) and, conversely, every element of $\mathcal W_T(\mu_0)$ is a Lagrangian
solution (Proposition~\ref{prop:weak-lagr}). Theorem~\ref{thm:micro} is proved in
Proposition~\ref{prop:admissible} and Lemma \ref{lem:empirical}; Theorem~\ref{thm:macro} is proved in
Theorems~\ref{thm:A}, \ref{thm:C-fast} and~\ref{thm:C-tail} and Corollaries~\ref{cor:weak-unique} and~\ref{cor:tail};
Theorem~\ref{thm:mfl} is proved in Theorem~\ref{thm:B} and Corollaries~\ref{cor:full-limit-fast} and~\ref{cor:tail}.

\section{The billiard flow for a given history}\label{sec:flow}
This section exploit many tools which are more or less standard in dispersing billiards. Nevertheless, they must be adapted to the present setting and, moreover, we prove every result to keep the section completely self-contained.
\subsection{Standing assumptions and notation}\label{subsec:standing}
With the notations of Section~\ref{sec:model}
we denote by $n(y)$ the unit normal pointing {into} ${\bf D}$,
and by $d(x)=\operatorname{dist}(x,\partial{\bf D})$ the distance to the boundary, which is
$\cC^2$ in a neighbourhood of $\partial{\bf D}$ in ${\bf D}$, with $\nabla d=n$ on $\partial{\bf D}$.\footnote{Note that
no finite horizon assumption is needed here (the collision time of a particle can be $\infty$).}

Fix $T>0$. The thermostat enters only through a scalar function of time. We call
\emph{history} any element of
\begin{equation}\label{eq:def-AT}
\cA_T:=\Bigl\{A\in \cC^0([0,T]):\ A(0)=0,\ |A(t)-A(s)|\le |E_0|\,|t-s|\ \ \forall s,t\in[0,T]\Bigr\}.
\end{equation}
$\cA_T$ is a convex subset of $\cC^0([0,T])$, compact by the Ascoli--Arzel\`a theorem, and
$\|A\|_\infty\le |E_0|T$ for $A\in\cA_T$. The function $A$ plays the role of $\int_0^t a(s)\,ds$,
where $a$ is the thermostat coefficient.

\subsection{Free flights, reflections and regular points}\label{subsec:linear-flow}

For reasons that will be clear in a moment, it is convenient to define the dynamics in the following coordinates.\\
 Let $A\in\cA_T$. Writing $v=e^{-A(t)}w$, the equations $\dot x=v$, $\dot v=E_0-\dot A(t)v$
become
\begin{equation}\label{eq:xw-system}
\dot x=e^{-A(t)}\,w,\qquad \dot w=e^{A(t)}\,E_0 ,
\end{equation}
which only involve $A$ (not $\dot A$) and have continuous coefficients. In particular, uniform convergence of histories implies uniform
convergence of the corresponding free-flight coefficients. Notice that
this change of variables does not remove the velocity jump at a collision because,
since specular reflection commutes with scalar multiplication,
\(w^+=\mathcal R_xw^-\). Its purpose is instead to obtain continuous
dependence on the history \(A\) and volume-preserving free-flight maps in
the variables \((x,w)\). We take \eqref{eq:xw-system} as the definition of the dynamics.\\
For $0\le s\le t\le T$ set
\[
\alpha_A(s,t):=\int_s^t e^{-A(r)}dr,\qquad
\gamma_A(s,t):=\int_s^t e^{A(r)}dr,\qquad
\beta_A(s,t):=\int_s^t e^{-A(r)}\gamma_A(s,r)\,dr .
\]
The {free flight} issued from $(x,w)$ at time $s$ is the explicit solution of \eqref{eq:xw-system}
\begin{equation}\label{eq:free-flight}
w(t)=w+\gamma_A(s,t)\,E_0,\qquad
x(t)=x+\alpha_A(s,t)\,w+\beta_A(s,t)\,E_0 ,
\end{equation}
understood in the universal cover $\bR^2$ of $\bT^2$. It is easy to check that, for fixed $s,t$, this is an affine map of
$(x,w)$ with unit Jacobian.

Let $z=(x,v)\in\mathring\Omega$ (so that $w=v$ at time $0$). The trajectory of $z$ is defined by
following \eqref{eq:free-flight} from time $0$ until the first {contact time} $\tau_1$, i.e. the
first time at which $x(\tau_1)\in\partial{\bf D}$. The contact is \emph{transversal} if
$\langle w(\tau_1^-),n(x(\tau_1))\rangle<0$ and \emph{tangential} if this scalar product vanishes.
At a transversal contact we set $w(\tau_1^+)=\cR_{x(\tau_1)}w(\tau_1^-)$, restart the free flight
from $(x(\tau_1),w(\tau_1^+))$ at time $\tau_1$, and so on.

\begin{oss} We stress that formula~\eqref{eq:free-flight} is defined in the whole cover
$\bR^2$, independently of the obstacles. We shall occasionally continue
a free flight beyond its first intersection with $\partial{\bf D}$,
ignoring reflections. This ``{auxiliary free continuation}" may pass
through the interior of the obstacles and should not be confused with the
physical billiard trajectory, which is reflected when it reaches the
boundary.
\end{oss}

\begin{defn}[Regular points]\label{def:regular}
A point $z\in\mathring\Omega$ is \emph{regular for $A$ on $[0,T]$} if the above procedure
produces a trajectory on the whole of $[0,T]$ with finitely many contact times
$0<\tau_1<\dots<\tau_m<T$, all transversal, and $x(T)\in\mathring{\bf D}$. We denote by
$\Reg_A\subset\mathring\Omega$ the set of regular points.
\end{defn}

For $z\in\Reg_A$ we write $(\widetilde X^A_t(z),W^A_t(z))$ for the trajectory in the cover, and
\begin{equation}\label{eq:def-Delta}
\begin{gathered}
V^A_t(z):=e^{-A(t)}W^A_t(z),\qquad
\Delta^A_t(z):=\widetilde X^A_t(z)-x=\int_0^tV^A_s(z)\,ds,\\
\Phi_A^{t}(z):=\bigl(\widetilde X^A_t(z)\ \mathrm{mod}\ \mathbb Z^2,\ V^A_t(z)\bigr).
\end{gathered}
\end{equation}
The {lifted displacement} $\Delta^A_t$ is continuous in $t$ (the velocity jumps, the position
does not). From \eqref{eq:free-flight} and $|\cR_xw|=|w|$ we get, for $z=(x,v)\in\Reg_A$ and $t\in[0,T]$,
\begin{equation}\label{eq:speed-bound}
|W^A_t(z)|\le |v|+c_T,\qquad |V^A_t(z)|\le C_T(1+|v|),\qquad |\Delta^A_t(z)|\le T\,C_T(1+|v|),
\end{equation}
with $c_T:=|E_0|Te^{|E_0|T}$ and $C_T:=e^{|E_0|T}(1+c_T)$, uniformly in $A\in\cA_T$.\\

Although the billiard flow is discontinuous near tangential contacts, its
collision pattern is stable along a regular trajectory.\footnote{Indeed, transversal
contacts persist under small perturbations of the history and of the initial
state, and their times and states vary continuously.} The following lemma is
the basic continuity statement used throughout the paper, both in the local
Jacobian analysis and in the fixed-point and mean-field arguments. The
displacement converges uniformly, whereas velocity convergence is asserted
only away from contact times because the velocity jumps at reflection.
\begin{lem}[Stability at regular points]\label{lem:stability}
Let $A\in\cA_T$ and $z\in\Reg_A$, with contact times $\tau_1<\dots<\tau_m$. Let $A_k\in\cA_T$ with
$A_k\to A$ uniformly, and $z_k\to z$ in $\mathring\Omega$. Then, for $k$ large enough,
$z_k\in\Reg_{A_k}$, the trajectory of $z_k$ has exactly $m$ contact times
$\tau^k_j\to\tau_j$, and
\[
\sup_{t\in[0,T]}\bigl|\Delta^{A_k}_t(z_k)-\Delta^A_t(z)\bigr|\to0,
\qquad
V^{A_k}_t(z_k)\to V^A_t(z)\quad\text{for every }t\notin\{\tau_1,\dots,\tau_m\}.
\]
In particular $\Reg_A$ is open.
\end{lem}

\begin{proof}
If $A_k\to A$ uniformly, then $\alpha_{A_k},\beta_{A_k},\gamma_{A_k}$ converge uniformly on
$\{0\le s\le t\le T\}$ together with their derivatives in $t$. Hence, by \eqref{eq:free-flight}, if
$s_k\to s$ and $(x_k,w_k)\to(x,w)$, the free flights $\gamma_k$ of $A_k$ issued from $(x_k,w_k)$ at
time $s_k$ converge to the free flight $\gamma$ of $A$ issued from $(x,w)$ at time $s$, in $\cC^1$ on
every compact subinterval of $(s,T]$, and with $\dot\gamma_k(s_k)\to\dot\gamma(s)$.

We argue by induction on the contacts. Assume that the data at the beginning of a flight converge,
where the limit datum is either the initial point (at $s=0$, with $x\in\mathring{\bf D}$) or a
post-collision state ($x\in\partial{\bf D}$, $\langle w,n(x)\rangle>0$). Let $\tau$ be the first contact
time of $\gamma$ after $s$ (transversal by assumption), or $\tau=T$ if there is none, in which case
$\gamma(T)\in\mathring{\bf D}$. Choose $\delta>0$ small. We thus have the following three situations:\\
\indent (i) In the post-collision case, $\tfrac{d}{dr}d(\gamma(r))\ge c>0$ on $[s,s+\delta]$, hence the same
holds for $\gamma_k$ on $[s_k,s+\delta]$ with $c/2$, and $\gamma_k$ has no contact on $(s_k,s+\delta]$.\\
\indent (ii) On $[s+\delta,\tau-\delta]$ (or $[s+\delta,T]$ if there is no contact) $d\circ\gamma$ is bounded
below by a positive constant, hence so is $d\circ\gamma_k$ for $k$ large.\\
\indent (iii) Let $Q$ be the obstacle met at time $\tau$ and let $d_Q^s$ be the signed distance to $\partial Q$,
positive in ${\bf D}$ and negative inside $Q$, which is $\cC^2$ in a neighbrhood of $\partial Q$. The flight
$\gamma$, continued beyond $\tau$, enters $Q$, and by transversality, after shrinking $\delta$,
$$\tfrac{d}{dr}d_Q^s(\gamma(r))\le -c<0\qquad \text{on} \qquad  [\tau-\delta,\tau+\delta].$$ Hence
$\tfrac{d}{dr}d_Q^s(\gamma_k(r))\le-c/2$ there for $k$ large, and $d_Q^s\circ\gamma_k$ has a unique zero
$\tau^k\to\tau$ in this interval; by (ii) it is the first contact of $\gamma_k$, it is transversal, and the
pre-collision states converge.
Since $(x,w)\mapsto(x,\cR_xw)$ is continuous, the post-collision states converge as well, and the
induction proceeds. After $m$ steps we obtain the first claim; the convergence of $V$ and of
$\Delta$ follows from the convergence of the flights and from \eqref{eq:speed-bound}.
\end{proof}

The previous construction defines the flow only for regular initial
conditions. To use it for absolutely continuous measures, we must therefore
show that the exceptional trajectories (infinitely many
collisions, tangential contact, or a contact at the final time) form a
null set. This is provided in the next proposition which also
shows that collisions at any given time are negligible and gives the
Jacobian formula needed to transport densities. In particular, the flow is
defined almost everywhere and is almost everywhere continuous at each fixed
time.\footnote{Another possibility would have been to work in
the quotient space obtained by identifying incoming and outgoing boundary states, but at the price of changing the topology. This setting allows us to work in the ``true" phase space.}

\begin{prop}[Almost every point is regular]\label{lem:ae-regular}
Let $A\in\cA_T$. Then:
\begin{enumerate}
\item[(a)] $\Leb(\mathring\Omega\setminus\Reg_A)=0$;
\item[(b)] for every $t\in[0,T]$ the set of $z\in\Reg_A$ having a contact at time $t$ is Lebesgue null;
\item[(c)] for every $t\in[0,T]$ the map $\Phi_A^t$ is injective on $\Reg_A$ and
$(\Phi_A^t)_\#(f_0\,dz)=f_t\,dz$ with
$f_t=e^{2A(t)}f_0\circ(\Phi^t_A)^{-1}$. In particular
$\|f_t\|_{L^\infty}\le e^{2|E_0|T}\|f_0\|_{L^\infty}$.
\end{enumerate}
\end{prop}

The proof uses the trichotomy below and two lemmas on the measure-theoretic structure of the flow,
which we prove first. For $0\le s,t\le T$ we denote by $F^{t,s}_A$ the free-flight map
\eqref{eq:free-flight} from time $s$ to time $t$ (also for $t<s$, by inverting it), acting on
$\bR^2\times\bR^2$; it is affine with unit Jacobian, $\cC^1$ in $(s,t)$, and
$F^{t,r}_A\circ F^{r,s}_A=F^{t,s}_A$. For the vector field of \eqref{eq:xw-system}, we write
$$ X_t(x,w):=(e^{-A(t)}w,e^{A(t)}E_0),$$ so that
$$\partial_tF^{t,s}_A=X_t\circ F^{t,s}_A  \qquad \text{and} \qquad  \partial_sF^{t,s}_A=-DF^{t,s}_A\,X_s.$$
A flight arriving at $\partial{\bf D}$ from ${\bf D}$ at time $\tau$ satisfies
$$\tfrac{d}{dt}d(x(t))|_{\tau^-}=e^{-A(\tau)}\langle w(\tau^-),n\rangle\le0,$$ so every contact is either
transversal or tangential; after a transversal contact $\langle w(\tau^+),n\rangle>0$ and the next flight
enters $\mathring{\bf D}$. Hence every $z\in\mathring\Omega\setminus\Reg_A$ belongs to at least one of the
following classes:
\begin{enumerate}
\item[(N1)] the trajectory has infinitely many transversal contacts in $(0,T)$;
\item[(N2)] after finitely many transversal contacts, the trajectory has a tangential contact at some
time $r\in(0,T]$;
\item[(N3)] the trajectory has finitely many contacts in $(0,T)$, all transversal, and
$x(T)\in\partial{\bf D}$.
\end{enumerate}

For $q\in[0,T]$ let $O_q\subset\mathring\Omega$ be the set of points whose trajectory is defined on
$[0,q]$ with finitely many contacts in $(0,q)$, all transversal, and position in $\mathring{\bf D}$ at
time $q$; recalling that $\widetilde X$ denote a trajectory in the cover, let $\Psi^q:O_q\to\mathring{\bf D}\times\bR^2$ be defined via 
$$\Psi^q(z):=(\widetilde X_q(z)\bmod\mathbb Z^2,W_q(z))$$
By Lemma~\ref{lem:stability} applied on $[0,q]$, $O_q$ is open.

\begin{lem}[Volume preservation]\label{lem:F1}
For every $q\in[0,T]$, $\Psi^q$ is an injective local $\cC^1$ diffeomorphism of $O_q$ with
$|\det D\Psi^q|=1$.
\end{lem}

\begin{proof}
\emph{Injectivity.} Let $z,z'\in O_q$ with $\Psi^q(z)=\Psi^q(z')$, and let $\tau$ be the largest of
the last contact times of the two trajectories before $q$ (if neither has a contact, both are the same
free flight on $[0,q]$ and $z=z'$). On $[\tau,q]$ the two trajectories coincide with the same free
flight run backwards from the common state at time $q$. At time $\tau$ this flight is on
$\partial{\bf D}$, so $\tau$ is a contact time for both, with the same outgoing state and therefore, the
reflection being an involution, the same incoming state. Repeating the argument on $[0,\tau]$ gives
$z=z'$ after finitely many steps.

\emph{Jacobian.} Let $z_0\in O_q$ have contact times $\tau_1<\dots<\tau_m$ in $(0,q)$, and fix times
$0=q_0<\tau_1<q_1<\dots<\tau_m<q_m=q$. If \(m=0\), then \(\Psi^q=F_A^{q,0}\), so the conclusion follows
immediately. Assume therefore that \(m\ge1\). By Lemma~\ref{lem:stability}, for $z$ near $z_0$ the
trajectory has exactly one contact in each $(q_{j-1},q_j)$, so that locally $\Psi^q$ is the composition
of $m$ maps, each describing a flight from time $q_{j-1}$ to time $q_j$ with exactly one transversal
contact. It suffices to treat one such map; to keep the notation light we take $q_{j-1}=0$, $q_j=q$.

Let $\pi$ be the nearest-point projection onto $\partial{\bf D}$, of class $\cC^2$ near
$\partial{\bf D}$, and set $\widetilde R(x,w):=(x,\cR_{\pi(x)}w)$, a $\cC^1$ map near
$\partial{\bf D}\times\bR^2$ extending the reflection. Its differential is block triangular,
\[
D\widetilde R=\begin{pmatrix}I&0\\ \ast&\cR_{\pi(x)}\end{pmatrix},
\qquad|\det D\widetilde R|=1 ,
\]
and $D\widetilde R^{-1}$ has the same form; in particular both preserve the $x$-component of a vector.
Let $d^s$ be the signed distance to the boundary of the obstacle met at the contact, positive in ${\bf D}$ and
negative inside the obstacle; it is of class $\cC^2$ near that boundary, with $\nabla d^s=n$ on it. The contact
time $\tau(z)$ is defined by $d^s\bigl(\pi_xF^{\tau,0}_A(z)\bigr)=0$. Since
$$\partial_\tau\,d^s(\pi_xF^{\tau,0}_A(z))=e^{-A(\tau)}\langle n,w^-\rangle\ne0,$$ by transversality, the
implicit function theorem gives $\tau\in\cC^1$ with
\begin{equation}\label{eq:grad-tau}
D\tau(z)[\delta z]=-\frac{\bigl\langle n,\pi_x\,DF^{\tau,0}_A(z)[\delta z]\bigr\rangle}{e^{-A(\tau)}\langle n,w^-\rangle},
\end{equation}
where $n=n(y)$, and $(y,w^-)=F^{\tau,0}_A(z)$ is the incoming state. Define, for $\sigma$ near $\tau(z_0)$,
\[
G(\sigma,z):=F^{q,\sigma}_A\circ\widetilde R\circ F^{\sigma,0}_A(z),\qquad\text{so that}\qquad
\Psi^q(z)=G(\tau(z),z).
\]
Then $D\Psi^q=M+a\otimes D\tau$,\footnote{Here \(a\otimes D\tau\) denotes the rank-one linear map
$
\delta z\longmapsto a\,D\tau[\delta z].
$} with $M:=D_zG$ and $a:=\partial_\sigma G$ evaluated at $\sigma=\tau(z)$. 
First, $M=DF^{q,\tau}_A\,D\widetilde R\,DF^{\tau,0}_A$ has $|\det M|=1$. Second, writing
$p:=(y,w^-)$ and using the formulas for $\partial_tF$ and $\partial_sF$,
\[
a=DF^{q,\tau}_A\bigl[D\widetilde R(p)\,X_\tau(p)-X_\tau(\widetilde Rp)\bigr],
\qquad
M^{-1}a=\bigl(DF^{\tau,0}_A\bigr)^{-1}\bigl[X_\tau(p)-D\widetilde R(p)^{-1}X_\tau(\widetilde Rp)\bigr].
\]
By \eqref{eq:grad-tau}, and since $D\widetilde R^{-1}$ preserves $x$-components,
\[
D\tau\bigl[M^{-1}a\bigr]
=-\frac{\bigl\langle n,\ e^{-A(\tau)}w^--e^{-A(\tau)}\cR_yw^-\bigr\rangle}{e^{-A(\tau)}\langle n,w^-\rangle}
=-\frac{2\langle n,w^-\rangle}{\langle n,w^-\rangle}=-2 .
\]
Recall the matrix determinant lemma: if \(M\) is invertible,
\(u\) is a vector and \(\ell\) is a linear functional, then
\[
\det(M+u\otimes\ell)
=
\det M\bigl(1+\ell(M^{-1}u)\bigr).
\]
Applying this with \(u=a\) and \(\ell=D\tau\), we obtain
\[
\det D\Psi^q
=
\det M\bigl(1+D\tau[M^{-1}a]\bigr)
=
-\det M.
\]
Since the free-flight maps have determinant one and
\(\mathcal R_{\pi(x)}\) is an orthogonal reflection, \(M\) is
invertible and
\[
|\det M|
=
|\det DF_A^{q,\tau}|\,
|\det D\widetilde R|\,
|\det DF_A^{\tau,0}|
=1.
\]

Therefore \(|\det D\Psi^q|=1\), and the inverse function theorem
shows that \(\Psi^q\) is a local \(\mathcal C^1\) diffeomorphism.$\qedhere$
\end{proof}

The next lemma is the last  ingredient needed to prove Proposition~\ref{lem:ae-regular}. It gives a flow-box parametrisation of initial
conditions in terms of the time, position and incoming velocity at a
collision. Under this change of variables, Lebesgue measure becomes the
natural flux measure \eqref{eq:flux-meas} through the collision space. The factor given by the
normal velocity will be essential both for proving that non-regular
trajectories form a null set and, later, for controlling nearly grazing
collisions in the uniqueness argument. We first need some more notations and definitions.\\
\indent Let 
\begin{equation}\label{eq:SigmaT}
\Sigma_T:=\{(t,y,w): t\in(0,T),\ y\in\partial{\bf D},\ \langle w,n(y)\rangle<0\}
\end{equation}
 be the set of
incoming boundary states, endowed with the flux measure
\begin{equation}\label{eq:flux-meas}
d\lambda:=e^{-A(t)}\,|\langle w,n(y)\rangle|\,dt\,d\sigma(y)\,dw .
\end{equation}
Above, the factor \(e^{-A(t)}|\langle w,n(y)\rangle|
=|\langle v,n(y)\rangle|\) is the normal physical speed at which
trajectories cross the boundary.\\
For $k\ge1$ let $R_k\subset\mathring\Omega$ be the set of points whose trajectory has at least $k$ contacts
in $(0,T)$, the first $k$ of which are transversal, and let $C_k:R_k\to\Sigma_T$ be the map giving time,
position and incoming velocity at the $k$-th contact. $R_k$ is open by Lemma~\ref{lem:stability}, applied on
an interval ending shortly after the $k$-th contact.

\begin{lem}[Flow-box formula]\label{lem:F2}
For every $k\ge1$, $C_k$ is injective and $$(C_k)_\#\bigl(\Leb|_{R_k}\bigr)=\lambda|_{C_k(R_k)}.$$
Moreover the sets $C_k(R_k)$, $k\ge1$, are pairwise disjoint.
\end{lem}

\begin{proof}
\emph{The case of a free flight.}
Fix \(q\in[0,T)\), and let \(R^{(q)}_1\) be the set of
\(\zeta\in\mathring{\bf D}\times\mathbb R^2\) whose free flight issued at
time \(q\) has a first contact in \((q,T)\), which is transversal. Let
\(C^{(q)}_1\) be the corresponding contact map and set
\[
\Gamma_q:=C^{(q)}_1(R^{(q)}_1)\subset\Sigma_T,
\]
recalling $\Sigma_T$ from \eqref{eq:SigmaT}.
The inverse of \(C^{(q)}_1\) on \(\Gamma_q\) is
\[
H(t,y,w):=F_A^{q,t}(y,w).
\]
Indeed, running the free flight backwards from the incoming contact state
\((y,w)\) at time \(t\) uniquely recovers its state at time \(q\).

Parametrise the relevant component of \(\partial{\bf D}\) by arclength,
and denote its unit tangent by \(\mathbf t(y)\). Using
\[
\partial_tF_A^{q,t}
=
-DF_A^{q,t}X_t,
\]
the columns of \(DH\), corresponding respectively to the variables
\(t\), the arclength coordinate of \(y\), and the two components of \(w\),
are
\[
-DF_A^{q,t}X_t(y,w),\qquad
DF_A^{q,t}(\mathbf t(y),0),\qquad
DF_A^{q,t}(0,e_1),\qquad
DF_A^{q,t}(0,e_2).
\]
Factoring out \(DF_A^{q,t}\) and using
\(|\det DF_A^{q,t}|=1\), we obtain
\[
\begin{aligned}
|\det DH|
&=
\left|
\det\bigl(
-X_t(y,w),(\mathbf t(y),0),(0,e_1),(0,e_2)
\bigr)
\right|\\
&=
\left|\det\bigl(-e^{-A(t)}w,\mathbf t(y)\bigr)\right|\\
&=
e^{-A(t)}|\langle w,n(y)\rangle|.
\end{aligned}
\]
In the second equality, the velocity component \(e^{A(t)}E_0\) of
\(X_t\) does not contribute because the last two columns form a basis in
the velocity variables. In the last equality we used that
\((\mathbf t(y),n(y))\) is an orthonormal basis.

Transversality implies that this determinant is nonzero. Hence \(H\), and
therefore \(C^{(q)}_1\), is a local \(\mathcal C^1\) diffeomorphism.
Moreover, for every Borel set \(B\subset\Gamma_q\), the change-of-variables
formula gives
\[
\begin{aligned}
\Leb\bigl((C^{(q)}_1)^{-1}(B)\bigr)
&=\Leb(H(B))\\
&=\int_B
e^{-A(t)}|\langle w,n(y)\rangle|\,
dt\,d\sigma(y)\,dw\\
&=\lambda(B).
\end{aligned}
\]
Thus
\[
(C^{(q)}_1)_\#
\bigl(\Leb|_{R^{(q)}_1}\bigr)
=
\lambda|_{\Gamma_q}.
\]

\emph{The general case.}
Let \(z_0\in R_k\), and choose \(q\) strictly between the
\((k-1)\)-st and the \(k\)-th contact times of the trajectory of \(z_0\),
with \(q=0\) when \(k=1\). By Lemma~\ref{lem:stability}, there is a
neighbourhood \(U\) of \(z_0\) such that, for every \(z\in U\),
\[
z\in O_q,\qquad
\Psi^q(z)\in R^{(q)}_1,
\]
and the \(k\)-th contact of \(z\) is the first contact after time \(q\).
Consequently,
\[
C_k=C^{(q)}_1\circ\Psi^q
\qquad\text{on }U.
\]
By Lemma~\ref{lem:F1}, \(\Psi^q\) is a local
volume-preserving \(\mathcal C^1\) diffeomorphism. Combining this with the
free-flight case shows that \(C_k\) is locally a \(\mathcal C^1\)
diffeomorphism and locally transports Lebesgue measure to \(\lambda\).

\emph{Injectivity and disjointness.}
Suppose that
$
C_k(z)=C_j(z')=(t,y,w).
$
Starting from the common incoming contact state \((y,w)\) at time \(t\),
run the trajectory backwards. The free dynamics is uniquely invertible,
and at every preceding contact the incoming state is uniquely recovered
from the outgoing one because the reflection is an involution. Thus the
entire trajectory on \([0,t]\), including its sequence of contact times,
is uniquely determined by \((t,y,w)\). It follows that
\[
z=z'
\qquad\text{and}\qquad
k=j.
\]
Therefore each \(C_k\) is injective, and the sets \(C_k(R_k)\) are pairwise
disjoint.

Since \(C_k\) is an injective local diffeomorphism, \(C_k(R_k)\) is open
in \(\Sigma_T\). The local change-of-variables formula, together with
injectivity, now gives, for every Borel set
\(B\subset C_k(R_k)\),
\[
\Leb(C_k^{-1}(B))=\lambda(B).
\]
Equivalently,
\[
(C_k)_\#\bigl(\Leb|_{R_k}\bigr)
=
\lambda|_{C_k(R_k)}. \qedhere
\]
\end{proof}

We are now ready to prove Proposition~\ref{lem:ae-regular}.

\begin{proof}[Proof of Proposition~\ref{lem:ae-regular}]
We divide the proof in three steps. In steps 1 and 2 we prove, respectively, that the classes (N1) and (N2) are null. In Step 3 we prove that the class (N3) is null and items \textup{(a)} and \textup{(b)}. In Step 4 we prove
\textup{(c)}. \\

\emph{Step 1.}
Points of class (N1) belong to \(R_k\) for every \(k\). Fix \(R>0\).
By \eqref{eq:speed-bound}, if \(z=(x,v)\in R_k\) with \(|v|\le R\),
then the incoming velocity at its \(k\)-th contact satisfies
$
|w|\le R+c_T.
$
By Lemma~\ref{lem:F2},
\[
\Leb\bigl(R_k\cap\{|v|\le R\}\bigr)
=
\lambda\Bigl(C_k\bigl(R_k\cap\{|v|\le R\}\bigr)\Bigr).
\]
Moreover, the sets
\[
C_k\bigl(R_k\cap\{|v|\le R\}\bigr),
\qquad k\ge1,
\]
are pairwise disjoint and are all contained in
$
\Sigma_T\cap\{|w|\le R+c_T\}.
$
Therefore
\[
\begin{aligned}
\sum_{k\ge1}
\Leb\bigl(R_k\cap\{|v|\le R\}\bigr)
&\le
\lambda\bigl(\Sigma_T\cap\{|w|\le R+c_T\}\bigr)\\
&\le
T e^{|E_0|T}|\partial{\bf D}|
\int_{|w|\le R+c_T}|w|\,dw
<\infty.
\end{aligned}
\]
The Borel--Cantelli lemma implies that almost every \(z\) with
\(|v|\le R\) belongs to only finitely many of the sets \(R_k\).
Since every point of class (N1) belongs to all the \(R_k\), the set of
points of class (N1) with \(|v|\le R\) is null. Letting \(R\to\infty\)
shows that (N1) is null.\\

\emph{Step 2}.
Fix \(q\in[0,T)\), and let \(E_q\) be the set of
\((x,\omega)\in\mathring{\bf D}\times\mathbb R^2\) whose free flight
issued at time \(q\) has a tangential contact with \(\partial{\bf D}\)
at some time \(r\in(q,T]\). Suppose that the contact occurs at
\(y\in\partial{\bf D}\) and that
$
w(r)=c\,\mathbf t(y)
$
for some \(c\in\mathbb R\). By \eqref{eq:free-flight}, the state
\((x,\omega)\) at time \(q\) is then given by
\[
\omega
=
c\,\mathbf t(y)-\gamma_A(q,r)E_0,
\qquad
x
=
y-\alpha_A(q,r)\omega-\beta_A(q,r)E_0.
\]
Thus, after parametrising each lifted boundary component by arclength,
the set \(E_q\) is contained in the image of the map
\[
(r,s,c)\longmapsto
\left(
y(s)-\alpha_A(q,r)\omega(r,s,c)-\beta_A(q,r)E_0,\,
\omega(r,s,c)
\right),
\]
where
\[
\omega(r,s,c)
=
c\,\mathbf t(y(s))-\gamma_A(q,r)E_0.
\]
This map is of class \(\mathcal C^1\) from a three-dimensional
parameter space into \(\mathbb R^4\). Covering the parameter space by
countably many bounded coordinate sets, its image has four-dimensional
Lebesgue measure zero. Since the boundary in the cover is a countable
union of lifted \(\mathcal C^3\) curves, it follows that
\[
\Leb(E_q)=0.
\]

Now let \(z\) be of class (N2), with transversal contacts
$
\tau_1<\dots<\tau_m
$
followed by a tangential contact at time $r$. Choose
$
q\in\mathbb Q\cap(\tau_m,r)
$
if \(m\ge1\), and \(q\in\mathbb Q\cap[0,r)\) if \(m=0\). Then
\(z\in O_q\) and
$
\Psi^q(z)\in E_q.
$
Consequently,
\[
(\mathrm{N2})
\subset
\bigcup_{q\in\mathbb Q\cap[0,T)}
(\Psi^q)^{-1}(E_q).
\]
Each preimage on the right-hand side is null. Indeed, \(O_q\) can be
covered by countably many open sets on which \(\Psi^q\) is a
\(\mathcal C^1\) diffeomorphism, by Lemma~\ref{lem:F1}, and a
\(\mathcal C^1\) diffeomorphism preserves null sets in both directions.
Hence (N2) is null.\\

\emph{Step 3}
Fix \(t\in(0,T]\) and \(q\in[0,t)\). Let \(E_{q,t}\) be the set of
states \((x,w)\) whose free flight issued at time \(q\) lies on
\(\partial{\bf D}\) at time \(t\). Then
\[
E_{q,t}
=
(F_A^{t,q})^{-1}
\bigl(\partial{\bf D}\times\mathbb R^2\bigr).
\]
Since \(F_A^{t,q}\) is an affine bijection and
\(\partial{\bf D}\times\mathbb R^2\) is Lebesgue null in the
four-dimensional phase space,
\[
\Leb(E_{q,t})=0.
\]

Suppose that \(z\) has finitely many transversal contacts in \((0,t)\)
and that its position lies on \(\partial{\bf D}\) at time \(t\).
Choose a rational \(q\) strictly between the last contact before \(t\)
and \(t\), or any rational \(q\in[0,t)\) if there is no previous
contact. Then
$
z\in O_q$ and $
\Psi^q(z)\in E_{q,t}.$ 
It follows that the set of such points is contained in
\[
\bigcup_{q\in\mathbb Q\cap[0,t)}
(\Psi^q)^{-1}(E_{q,t}),
\]
and is therefore null by Lemma~\ref{lem:F1}. For \(t=T\), this proves
that (N3) is null. For an arbitrary \(t\in(0,T]\), it proves item
\textup{(b)}. The case \(t=0\) is immediate, since the points of
\(\Reg_A\) have initial position in \(\mathring{\bf D}\). Steps 1--3, together with the trichotomy, prove item \textup{(a)}.\\

\emph{Step 4}.
Let \(\Reg_A^t\) be the set of points of \(\Reg_A\) having no contact
at time \(t\). By items \textup{(a)} and \textup{(b)},
\[
\Leb(\mathring\Omega\setminus\Reg_A^t)=0.
\]
Moreover, \(\Reg_A^t\subset O_t\). On \(O_t\), the map \(\Psi^t\) is
injective and satisfies
$
|\det D\Psi^t|=1
$
by Lemma~\ref{lem:F1}. Injectivity of \(\Phi_A^t\) on all of
\(\Reg_A\), including points having a contact exactly at time \(t\),
follows from the same backward-determinism argument: from the outgoing
state at the contact one first applies the inverse reflection and then
runs the preceding free flight backwards.

At time \(t\), the change from \(W_t\) to \(V_t\) is
$
V_t=e^{-A(t)}W_t.
$
Consequently, on \(\Reg_A^t\),
\[
|\det D\Phi_A^t|
=
e^{-2A(t)}|\det D\Psi^t|
=
e^{-2A(t)}.
\]
The change-of-variables formula therefore gives
\[
(\Phi_A^t)_\#(f_0\,dz)=f_t\,dz,
\]
where
\[
f_t(z)
=
e^{2A(t)}
f_0\bigl((\Phi_A^t)^{-1}(z)\bigr)
\]
on \(\Phi_A^t(\Reg_A^t)\), and \(f_t=0\) outside this set. Finally,
since
$
A(t)\le |A(t)|\le |E_0|T,
$
we obtain
\[
\|f_t\|_{L^\infty}
\le
e^{2|E_0|T}\|f_0\|_{L^\infty}.\qedhere
\]

\end{proof}
\section{Lagrangian solutions and existence}\label{sec:lagr}

\subsection{Current, energy and the fixed point}\label{subsec:lagrangian}

Let $\mu_0\in\cP(\Omega)$ be absolutely continuous with respect to Lebesgue, with
\begin{equation}\label{eq:mu0-normalised}
\int_\Omega|v|^2\,d\mu_0=1 .
\end{equation}
By Proposition~\ref{lem:ae-regular}(a), for every $A\in\cA_T$ the maps $V^A_t,\Delta^A_t,\Phi^t_A$ (see \eqref{eq:def-Delta}) are defined
$\mu_0$-almost everywhere, and we set
\begin{equation}\label{eq:def-j-u-A}
j_A(t):=\int_\Omega\langle E_0,V^A_t\rangle\,d\mu_0,\qquad
u_A(t):=\int_\Omega|V^A_t|^2\,d\mu_0,\qquad
\mu^A_t:=(\Phi^t_A)_\#\mu_0 .
\end{equation}

\begin{lem}\label{lem:j-u}
Let $A\in\cA_T$. Then $j_A$ and $u_A$ are continuous on $[0,T]$, $u_A>0$,
$|j_A|\le|E_0|\,u_A^{1/2}$, and for every $t\in[0,T]$
\begin{align}
\int_0^tj_A(s)\,ds&=\int_\Omega\langle E_0,\Delta^A_t\rangle\,d\mu_0, \label{eq:j-displacement}\\
u_A(t)&=1+\int_0^t\bigl(2j_A(s)-2\dot A(s)\,u_A(s)\bigr)\,ds . \label{eq:u-balance}
\end{align}
Moreover $t\mapsto\mu^A_t$ is narrowly continuous in $\cP(\Omega)$.
\end{lem}

\begin{proof}
For $z\in\Reg_A$ the function $t\mapsto V^A_t(z)$ is continuous except at the contact times, and
$t\mapsto|V^A_t(z)|$ is continuous. For fixed $t$, almost every $z$ has no contact at time $t$
(Proposition~\ref{lem:ae-regular}(b)); continuity of $j_A$, $u_A$ and of $t\mapsto\int\varphi\,d\mu^A_t$,
$\varphi\in\cC_b(\Omega)$, then follows from dominated convergence, using \eqref{eq:speed-bound} and
\eqref{eq:mu0-normalised}. The bound on $j_A$ is the Cauchy--Schwarz inequality.
Identity \eqref{eq:j-displacement} is Fubini's theorem applied to \eqref{eq:def-Delta}.
Between contacts $V=e^{-A}W$ is Lipschitz in $t$ with $\dot V=E_0-\dot A\,V$ for a.e.\ $t$, so that
$\tfrac{d}{dt}|V|^2=2\langle E_0,V\rangle-2\dot A|V|^2$, while $|V|$ is continuous across contacts;
integrating in time and then in $\mu_0$ gives \eqref{eq:u-balance}. Finally $u_A(t)>0$ because
$\mu^A_t$ is absolutely continuous (Proposition~\ref{lem:ae-regular}(c)) and cannot be concentrated on $\{v=0\}$.
\end{proof}
We can now provide the Lagrangian formulation of a solution.
\begin{defn}[Lagrangian solution]\label{def:lagrangian}
A \emph{Lagrangian solution} on $[0,T]$ with initial datum $\mu_0$ is a pair
$(A,\{\mu_t\}_{t\in[0,T]})$ with $A\in\cA_T$,
\begin{equation}\label{eq:fixed-point-original}
A(t)=\int_0^t\frac{j_A(s)}{u_A(s)}\,ds\qquad\forall t\in[0,T],
\end{equation}
and $\mu_t=\mu^A_t=(\Phi^t_A)_\#\mu_0$.
\end{defn}
One advantage of the above definition is the following equivalent forms of the fixed point, direct consequence of Lemma \ref{lem:j-u}.
\begin{lem}\label{lem:equivalent-fp}
For $A\in\cA_T$ the following are equivalent:
\begin{enumerate}
\item[(i)] $A(t)=\int_0^tj_A/u_A\,ds$ for all $t$;
\item[(ii)] $A(t)=\int_\Omega\langle E_0,\Delta^A_t\rangle\,d\mu_0$ for all $t$;
\item[(iii)] $A(t)=\int_0^tj_A\,u_A^{-1/2}\,ds$ for all $t$.
\end{enumerate}
In each case $u_A\equiv1$ on $[0,T]$ and $\dot A=j_A$.
\end{lem}

\begin{proof}
Set $y:=u_A-1$, so $y(0)=0$. In case (i), \eqref{eq:u-balance} gives $\dot y=0$. In case (ii),
by \eqref{eq:j-displacement} $\dot A=j_A$ and \eqref{eq:u-balance} gives $\dot y=-2j_A\,y$. In case (iii),
$\dot y=2j_A(1-u_A^{1/2})=-2j_A\,y/(1+u_A^{1/2})$. In the last two cases $|\dot y|\le C|y|$ with
$C=2|E_0|\sup_{[0,T]}u_A^{1/2}$, hence $y\equiv0$ by Gr\"onwall's lemma. Once $u_A\equiv1$ the three
identities coincide.
\end{proof}
Lemmas~\ref{lem:j-u} and~\ref{lem:equivalent-fp} allow us to easily prove that Lagrangian solutions are weak solutions, as shown in the next proposition.
\begin{prop}\label{prop:lagr-weak}
If $(A,\{\mu_t\})$ is a Lagrangian solution, then $\{\mu_t\}$ is a weak solution in the sense of
Definition~\ref{def:w-sol}, with $\int|v|^2d\mu_t=1$ for all $t$.
\end{prop}
 
\begin{proof}
Narrow continuity, which implies Definition~\ref{def:w-sol}{\rm(i)}, and the second moment come from
Lemmas~\ref{lem:j-u} and~\ref{lem:equivalent-fp}.
Let $\varphi\in\mathcal D_{\mathrm{sym}}$ and $z\in\Reg_A$. The function
$g(t):=\varphi(t,\Phi^t_A(z))$ is continuous on $[0,T]$, because $\varphi$ is symmetric at the contact
times, and Lipschitz between contacts; hence it is absolutely continuous, with
$$\dot g=\partial_t\varphi+\langle\nabla_x\varphi,V\rangle+\langle\nabla_v\varphi,E_0-\dot A\,V\rangle$$
a.e., and $g(T)=0$. Thus $-\varphi(0,z)=\int_0^T\dot g\,dt$. Integrating with respect to $\mu_0$
(the integrand is bounded since $\varphi$ has compact support) and using
$\dot A=j_A/u_A=j_{\mu_t}/u_{\mu_t}$ gives \eqref{eq:weak-form-final}.
\end{proof}

\begin{oss}\label{rem:converse}
The converse holds for absolutely continuous weak solutions, see
Proposition~\ref{prop:weak-lagr}.
\end{oss}

\subsection{Existence}\label{subsec:existence}  We now prove existence of weak solutions of the Vlasov equation, concluding
the existence part (1) of Theorem~\ref{thm:macro}. Our
strategy is to construct a self-consistent Lagrangian solution and then use
Proposition~\ref{prop:lagr-weak} to conclude that it is a weak solution.
The construction works on every finite time interval and requires no
restriction on the size of the field or on the velocity support of the
initial datum.
\begin{thm}[Existence]\label{thm:A}
Let $\mu_0=f_0\,dz\in\cP(\Omega)$ satisfy \eqref{eq:mu0-normalised}. Then for every $T>0$ there exists a
Lagrangian solution $(A,\{\mu_t\})$ on $[0,T]$ with initial datum $\mu_0$. It satisfies
$\int|v|^2d\mu_t=1$, $\mu_t=f_t\,dz$, and if $f_0\in L^\infty$ then
$\|f_t\|_{L^\infty}\le e^{2|E_0|T}\|f_0\|_{L^\infty}$. In particular, \(\{\mu_t\}\) is a weak solution of
the Vlasov equation.
\end{thm}

\begin{proof}
Define $\cS:\cA_T\to\cC^0([0,T])$ by
\[
\cS(A)(t):=\int_0^t\frac{j_A(s)}{u_A(s)^{1/2}}\,ds .
\]
By Lemma~\ref{lem:j-u} the integrand is continuous and bounded by $|E_0|$, so $\cS(\cA_T)\subset\cA_T$.\\
 Let $A_k\to A$ uniformly in $\cA_T$ and fix $s\in[0,T]$. For $\mu_0$-a.e.\ $z$ we
have $z\in\Reg_A$ with no contact at time $s$ (Proposition~\ref{lem:ae-regular}), hence
$V^{A_k}_s(z)\to V^A_s(z)$ by Lemma~\ref{lem:stability} (with $z_k=z$). By \eqref{eq:speed-bound} the
functions $|V^{A_k}_s|^2$ are dominated by $C_T^2(1+|v|)^2\in L^1(\mu_0)$, so
$j_{A_k}(s)\to j_A(s)$ and $u_{A_k}(s)\to u_A(s)>0$. By dominated convergence theorem in $s$ we have
$\cS(A_k)(t)\to\cS(A)(t)$ for every $t$, and the convergence is uniform because the functions
$\cS(A_k)$ are equi-Lipschitz.

Since $\cA_T$ is convex and compact, Schauder's theorem yields $A\in\cA_T$ with $\cS(A)=A$. By
Lemma~\ref{lem:equivalent-fp}, (iii)$\Rightarrow$(i), the pair $(A,\{\mu^A_t\})$ is a Lagrangian
solution with $u_A\equiv1$. The statements on the density are Proposition~\ref{lem:ae-regular}(c). 
By Proposition \ref{prop:lagr-weak}, this solution is a weak solution, whereby concluding part (1) of Theorem~\ref{thm:macro}. 
\end{proof}

\subsection{Absolutely continuous weak solutions}\label{subsec:weak-lagr}

We conclude this section by proving a converse of Proposition~\ref{prop:lagr-weak} in the class $\mathcal W_T(\mu_0)$ of
Definition~\ref{def:w-sol}. 
Note that this implication is not needed for the
existence result above, which holds with no assumption on the velocities.
 Its purpose is to identify weak and Lagrangian
solutions in the class \(\mathcal W_T(\mu_0)\), so that uniqueness of the
Lagrangian fixed point will later imply uniqueness of weak solutions in
that class. The target of this subsection is thus the following Proposition, which proof is given at the end of the subsection. .

\begin{prop}[Weak solutions are Lagrangian]\label{prop:weak-lagr}
Let $\mu_0=f_0\,dz$ with $\int|v|^2d\mu_0=1$ and let $\{\mu_t\}\in\mathcal W_T(\mu_0)$. Let
$A(t):=\int_0^tj_\mu(s)\,ds$. Then $\mu_t=(\Phi^t_A)_\#\mu_0$ for every $t\in[0,T]$, and $(A,\{\mu_t\})$ is a
Lagrangian solution.
\end{prop}
To prove it, we use Proposition~\ref{prop:conservation} and one more lemma.
The next lemma provides test functions that are constant along the
characteristics of the prescribed billiard flow. These functions will be
used to recover the initial measure from a weak solution and hence to prove
that weak solutions are transported by the flow. The only delicate point is
their regularity near a collision, where the incoming and outgoing states
must be treated by separate local inverse branches. Since this construction
is technical, its proof is provided in
Appendix~\ref{app:transported-tests} for readability.
\begin{lem}[Test functions transported by the flow]\label{lem:invariant-test}
Let $A\in\cA_T$ and $\Theta\in\cC^1_c(\Reg_A)$. For $t\in[0,T]$ and $z\in\Omega$ set
\[
\varphi_\Theta(t,z):=\Theta(z_0)\quad\text{if }z=\Phi^t_A(z_0)\text{ for some }z_0\in\operatorname{supp}\Theta,
\qquad\varphi_\Theta(t,z):=0\quad\text{otherwise},
\]
where at a contact time of $z_0$ both the incoming and the outgoing state are regarded as $\Phi^t_A(z_0)$.
Then $\varphi_\Theta$ is well defined, continuous, symmetric and compactly supported in $[0,T]\times\Omega$; it
is Lipschitz in $t$; for every $t$ it is of class $\cC^1$ in $z$ on $\Omega$, regarded as a manifold with
boundary, and $\nabla_z\varphi_\Theta$ is continuous on $[0,T]\times\Omega$. Finally,
\begin{equation}\label{eq:transported}
\partial_t\varphi_\Theta+\bigl\langle\nabla_z\varphi_\Theta,(v,\,E_0-\dot A(t)\,v)\bigr\rangle=0
\qquad\text{for almost every }(t,z)\in(0,T)\times\mathring\Omega .
\end{equation}
\end{lem}

We are now ready to prove Proposition \ref{prop:weak-lagr}.
\begin{proof}[Proof of Proposition \ref{prop:weak-lagr}.]
By Proposition~\ref{prop:conservation}, $j_\mu$ is Borel with $|j_\mu|\le|E_0|$; hence $A\in\cA_T$, $\dot A=j_\mu$ a.e.\ and $V_{\mu_t}(x,v)=(v,E_0-\dot A(t)v)$ for a.e.\ $t$.

The identity \eqref{eq:weak-form-final} holds for every function
$\varphi$ on $[0,T)\times\Omega$ which is continuous, symmetric, compactly supported, Lipschitz in $t$, and
such that $\nabla_z\varphi$ is continuous. Indeed, extend $\varphi$ by zero for $t\ge T$, let $\rho_\delta$ be a mollifier supported in $[0,\delta]$
and $\varphi_\delta(t,z):=\int\rho_\delta(s)\varphi(t+s,z)\,ds$, which belongs to $\mathcal D_{\mathrm{sym}}$ for
$\delta$ small. As $\delta\to0$, $\varphi_\delta\to\varphi$ and $\nabla_z\varphi_\delta\to\nabla_z\varphi$
uniformly, while $\partial_t\varphi_\delta$ is uniformly bounded and converges to $\partial_t\varphi$ at
almost every $(t,z)$. By Definition~\ref{def:w-sol}{\rm(i)}, $t\mapsto\mu_t$ is a Borel probability
kernel, so $M(dt\,dz):=\mu_t(dz)\,dt$ is a well-defined measure on $[0,T]\times\Omega$; since $\mu_t\ll dz$ for
every $t$, $M\ll dt\,dz$. Hence $\partial_t\varphi_\delta\to\partial_t\varphi$ $M$-almost everywhere, and
dominated convergence allows us to pass to the limit in \eqref{eq:weak-form-final}.

Next, let $\Theta\in\cC_c^1(\Reg_A)$ and $\zeta\in\cC_c^1([0,T))$. By
Lemma~\ref{lem:invariant-test} and the previous paragraph we may use $\zeta(t)\varphi_\Theta(t,z)$ as a test function. By
\eqref{eq:transported}, which holds $M$-almost everywhere because $M\ll dt\,dz$ and $M$ does not charge
$[0,T]\times\partial\Omega$, the transport terms cancel, and we are left with
\[
\int_0^T\zeta'(t)\,m(t)\,dt+\zeta(0)\int_\Omega\Theta\,d\mu_0=0,\qquad m(t):=\int_\Omega\varphi_\Theta(t,z)\,d\mu_t(z).
\]
Hence $m(t)=\int\Theta\,d\mu_0$ for almost every $t$. The function $m$ is continuous: $\varphi_\Theta(t,\cdot)$ is
symmetric, bounded and continuous, it depends continuously on $t$ in the uniform norm, and $\{\mu_t\}$
satisfies Definition~\ref{def:w-sol}{\rm(i)}. Therefore
\begin{equation}\label{eq:pullback}
\int_\Omega\varphi_\Theta(t,z)\,d\mu_t(z)=\int_\Omega\Theta\,d\mu_0\qquad\text{for every }t\in[0,T]\text{ and }\Theta\in\cC^1_c(\Reg_A).
\end{equation}

We can now conclude. Fix $t$ and let $\Reg^t_A$ be the set of points of $\Reg_A$ with no contact at time
$t$. By Proposition~\ref{lem:ae-regular}, $\Reg^t_A$ is open with null complement  and $\Phi^t_A$ is a homeomorphism of
$\Reg^t_A$ onto an open set $G_t\subset\mathring\Omega$. For $\Theta\in\cC^1_c(\Reg^t_A)$ we have
$\varphi_\Theta(t,\cdot)=\Theta\circ(\Phi^t_A)^{-1}$ on $G_t$ and $0$ elsewhere, so \eqref{eq:pullback} says
that the measures $((\Phi^t_A)^{-1})_\#(\mu_t|_{G_t})$ and $\mu_0|_{\Reg^t_A}$ coincide as Radon measures on
$\Reg^t_A$. Since $\mu_0$ is absolutely continuous, $\mu_0(\Reg^t_A)=1$; exhausting $\Reg^t_A$ by compact sets we get
$\mu_t(G_t)=\mu_0(\Reg^t_A)=1$, and therefore $\mu_t=(\Phi^t_A)_\#\mu_0$. In particular
$j_A=j_\mu$ and $u_A\equiv1$, so $A=\int_0^{\cdot}j_A/u_A$, which is \eqref{eq:fixed-point-original}.
\end{proof}

\section{The particle system and the mean-field limit}\label{sec:mf}
This section is devoted to the full proof of Theorem~\ref{thm:micro} and the proof of Theorem \ref{thm:mfl} under the cut-off on the velocity.\\
\subsection{Particle system} We consider $N$ particles $\mathbf z^N=(z^N_1,\dots,z^N_N)\in\mathring\Omega^N$,
$z^N_i=(x^N_i,v^N_i)$, on the energy shell
\begin{equation}\label{eq:shell}
\frac1N\sum_{i=1}^N|v^N_i|^2=1,
\end{equation}
evolving according to \eqref{eq:equation of motion2} with specular reflection. On this shell the
thermostat coefficient is $$a_N(t)=\frac1N\sum_i\langle E_0,v_i(t)\rangle,$$ and $|a_N|\le|E_0|$ by the
Cauchy--Schwarz inequality and the conservation of \eqref{eq:shell}.

\begin{defn}[Admissible configurations]\label{def:admissible}
A configuration $\mathbf z^N$ satisfying \eqref{eq:shell} is \emph{admissible on $[0,T]$} if the
$N$-particle dynamics issued from it is defined on $[0,T]$ and each particle has finitely many
contacts with $\partial{\bf D}$ in $[0,T]$, all transversal. \\

Note that the above definition is equivalent to the existence of $A_N\in\cA_T$ such that the
trajectory of each $z^N_i$ under the flow of $A_N$ is defined on $[0,T]$ with finitely many contacts, all transversal, and
\begin{equation}\label{eq:AN-fp}
A_N(t)=\frac1N\sum_{i=1}^N\bigl\langle E_0,\Delta^{A_N}_t(z^N_i)\bigr\rangle\qquad\forall t\in[0,T].
\end{equation}
\end{defn}

Indeed, given the $N$-particle trajectory, $A_N(t):=\int_0^ta_N(s)\,ds$ belongs to $\cA_T$, each particle
follows by definition the one-particle flow of $A_N$, and \eqref{eq:AN-fp} is the time integral of
$a_N=\frac1N\sum_i\langle E_0,v_i\rangle$; conversely, if $A_N\in\cA_T$ satisfies \eqref{eq:AN-fp} then
$\dot A_N=\frac1N\sum_i\langle E_0,v_i\rangle$, and the energy identity gives $\frac1N\sum_i|v_i|^2= 1$ as in
Lemma~\ref{lem:equivalent-fp}, so that the particles solve \eqref{eq:equation of motion2}. Also, note that particles interact only through $A_N$, i.e. simultaneous
contacts of different particles play no role. We denote
\begin{equation}\label{eq:empirical}
\delta_{N,0}:=\frac1N\sum_{i=1}^N\delta_{z^N_i},\qquad
\delta_{N,t}:=\frac1N\sum_{i=1}^N\delta_{\Phi^t_{A_N}(z^N_i)} 
\end{equation}
the initial empirical measure and its push-forward by the flow, respectively. \\

The following proposition state that almost every configuration is admissible, and it is the main technical and key part of the section. 
The proof is the many-particle analogue of the argument used for
Proposition~\ref{lem:ae-regular} and it might look technical, but the main idea is quite naive, which is the follow. 
Ordinary transversal collisions can be parametrised by their time and by the
state of the system at the boundary. The total flux of trajectories through
the boundary during a finite time interval is finite, which excludes an
infinite number of such collisions for almost every initial configuration.
The remaining problematic configurations are tangential contacts and simultaneous contacts
of two particles. Nevertheless, it turns out that they form a
lower-dimensional set and hence, as the dynamics is smooth and locally invertible
between transversal collisions, their preimages are null as well.
\begin{prop}\label{prop:admissible}
Let $S_N$ be the shell \eqref{eq:shell}, endowed with its $(4N-1)$-dimensional surface measure
$\mathfrak m$. For every $T>0$, $\mathfrak m$-almost every configuration in $S_N$ is admissible on $[0,T]$, and along its
trajectory no two particles are on $\partial{\bf D}$ at the same time.
\end{prop}

\begin{proof}
We follow the proof of Proposition~\ref{lem:ae-regular}. Configurations with a particle on $\partial{\bf D}$ at time
$0$ form an $\mathfrak m$-null set and are discarded. In the other cases,  there are two possible obstructions to admissibility: a trajectory may have
infinitely many ordinary transversal contacts in finite time, or it may reach
a singular contact, namely a tangency or a simultaneous contact of two
particles. The first possibility will be excluded by estimating the total
flux through the contact hypersurfaces, and the second by observing that the
set of singular contact states has lower dimension.\\
\indent It is convenient to distinguish the physical shell $S_N$ from the
unrestricted energy shell
\[
\widehat S_N
:=
(\mathbb T^2)^N
\times
\left\{
\mathbf v\in\bR^{2N}:|\mathbf v|^2=N
\right\},
\]
or its lift to the cover in the position variables. The free flow is defined
on $\widehat S_N$, without regard to the obstacles, and is used between
successive contacts. Let $\mathcal E:=(E_0,\dots,E_0)\in\bR^{2N}$, so that
$
\alpha_N(\mathbf v)
=
\frac{\langle\mathcal E,\mathbf v\rangle}{|\mathbf v|^2},
$
and let
$
\mathbf X(\mathbf x,\mathbf v)
:=
\bigl(\mathbf v,\mathcal E-\alpha_N(\mathbf v)\mathbf v\bigr).
$
On $\widehat S_N$ one has $|\mathbf v|^2=N$ and
$
|\alpha_N(\mathbf v)|
\le
\frac{|\mathcal E|}{|\mathbf v|}
=
|E_0|.
$
Thus $\mathbf X$ is smooth and bounded. Moreover,
$
\left\langle
\mathbf v,\mathcal E-\alpha_N(\mathbf v)\mathbf v
\right\rangle
=0,
$
so $\mathbf X$ is tangent to every level set of
$
U_N(\mathbf x,\mathbf v):=|\mathbf v|^2.
$
It therefore generates a complete flow $\varphi^t$ of smooth
diffeomorphisms of $\widehat S_N$. Since $\alpha_N$ is homogeneous of degree $-1$,
$
\langle\mathbf v,\nabla_{\mathbf v}\alpha_N\rangle=-\alpha_N,
$
and hence
\[
\operatorname{div}\mathbf X
=
-\langle\mathbf v,\nabla_{\mathbf v}\alpha_N\rangle
-2N\alpha_N
=
-(2N-1)\alpha_N.
\]
Note that, if a $\cC^1$ diffeomorphism $F$ preserves $U_N$, then its Jacobian
on a regular level set of $U_N$ is
\[
J_{U_N}F(\mathbf z)
=
|\det DF(\mathbf z)|
\frac{|\nabla U_N(F(\mathbf z))|}
     {|\nabla U_N(\mathbf z)|}.
\]
Here $U_N\circ\varphi^t=U_N$ and
$
|\nabla U_N|=2|\mathbf v|=2\sqrt N
$
on $\widehat S_N$. Therefore the tangential Jacobian of $\varphi^t$
coincides with its Jacobian. Liouville's formula gives
\begin{equation}\label{eq:surface-jac}
J_t(\mathbf z)
=
\exp\left(
-(2N-1)\int_0^t
\alpha_N(\varphi^s\mathbf z)\,ds
\right)
\in[\Lambda^{-1},\Lambda],
\qquad
\Lambda:=e^{(2N-1)|E_0|T},
\quad |t|\le T .
\end{equation}
Next, recall that a {contact} is a time at which some particle is on $\partial{\bf D}$ and it is
{transversal} if exactly one particle, say the $k$-th, is on the boundary, with
$\langle v_k,n(x_k)\rangle<0$; then $v_k$ is replaced by $\cR_{x_k}v_k$, which preserves the shell, and the free
flow is restarted. As in the trichotomy of Section~\ref{sec:flow}, a configuration for which the conclusion of
the proposition fails belongs to one of the following classes: (N1) its trajectory has infinitely many
transversal contacts in $(0,T)$; (N2) after finitely many transversal contacts, at some time $r\in(0,T]$ a
particle touches $\partial{\bf D}$ tangentially or two particles are on $\partial{\bf D}$ simultaneously. A
transversal contact at the final time $T$ is allowed by Definition~\ref{def:admissible}.\\
Let $O_q$ be the set of configurations whose trajectory on $[0,q]$ has
finitely many contacts, all transversal, and none at time $q$, and let $\Psi^q_N$ be the corresponding time-$q$
map. The proof of Lemma~\ref{lem:stability} applies to the flow $\varphi^t$, which is $\cC^1$ in time and in the
initial point, and shows that $O_q$ is open, $\Psi^q_N$ is continuous and injective.
 For the Jacobian it suffices, as in Lemma~\ref{lem:F1}, to consider one contact, of the $k$-th
particle, at time $\tau(\mathbf z)$. Let $\widetilde R_k(\mathbf x,\mathbf v)$ be the map which replaces $v_k$ by
$\cR_{\pi(x_k)}v_k$; it preserves the shell and $\mathfrak m$, and both $D\widetilde R_k$ and its inverse leave
the position components unchanged. With $G(\sigma,\mathbf z):=\varphi^{q-\sigma}\circ\widetilde R_k\circ\varphi^{\sigma}(\mathbf z)$
we have $\Psi^q_N(\mathbf z)=G(\tau(\mathbf z),\mathbf z)$ and $D\Psi^q_N=M+a\otimes D\tau$, where, writing
$p$ for the incoming state,
\begin{gather*}
M=D\varphi^{q-\tau}\,D\widetilde R_k\,D\varphi^{\tau},\qquad
a=D\varphi^{q-\tau}\bigl[D\widetilde R_k(p)\mathbf X(p)-\mathbf X(\widetilde R_kp)\bigr],\\
D\tau[\delta\mathbf z]=-\frac{\langle n,\pi_{x_k}D\varphi^{\tau}[\delta\mathbf z]\rangle}{\langle n,v_k^-\rangle}.
\end{gather*}
The $x_k$-components of $\mathbf X(p)$ and of $\mathbf X(\widetilde R_kp)$ are $v_k^-$ and $\cR v_k^-$, so
\[
D\tau[M^{-1}a]=-\frac{\langle n,v_k^--\cR v_k^-\rangle}{\langle n,v_k^-\rangle}=-2 ,
\]
and the matrix determinant lemma, applied on the tangent space of the shell, shows that the Jacobian of
$\Psi^q_N$ with respect to $\mathfrak m$ equals, in absolute value, that of $M$, i.e.\ the product of the
Jacobians \eqref{eq:surface-jac} of the free flows before and after the contact. By induction on the contacts,
$\Psi^q_N$ is a local $\cC^1$ diffeomorphism of $O_q$ whose Jacobian with respect to $\mathfrak m$ lies in
$[\Lambda^{-1},\Lambda]$.

Now, let $\partial_kS$ be the set of points of the shell with
$x_k\in\partial{\bf D}$, a hypersurface of the shell with unit normal $(n(x_k),0)$ in the $x_k$-slot, and let
$\mathfrak m_\partial$ be its surface measure. For the free flow started at time $q$, the map sending a
configuration to the time $t$ and the state $\mathbf y\in\partial_kS$ of its first contact has inverse
$(t,\mathbf y)\mapsto\varphi^{-(t-q)}(\mathbf y)$, whose Jacobian with respect to $dt\,d\mathfrak m_\partial$ is
$$|\langle\mathbf X(\mathbf y),(n,0)\rangle|\,J_{-(t-q)}(\mathbf y)=|\langle v_k,n(y_k)\rangle|\,J_{-(t-q)}(\mathbf y).$$
Composing with $\Psi^q_N$ as in Lemma~\ref{lem:F2}, the $j$-th contact map $C_j$ is an injective local
diffeomorphism of the open set $R_j$ of configurations having at least $j$ contacts in $(0,T)$, the first $j$
transversal, and
\[
(C_j)_\#\bigl(\mathfrak m|_{R_j}\bigr)\le\Lambda^2\,|\langle v_k,n(y_k)\rangle|\,dt\,d\mathfrak m_\partial
\quad\text{on }C_j(R_j) .
\]
The images $C_j(R_j)$, $j\ge1$, are pairwise disjoint by backward determinism, and
$|\langle v_k,n\rangle|\le\sqrt N$ on the shell, which is compact modulo translations of the cover. Hence
\[
\sum_{j\ge1}\mathfrak m(R_j)\le\Lambda^2\,T\sqrt N\sum_{k=1}^N\mathfrak m_\partial(\partial_kS)<\infty ,
\]
and by the Borel--Cantelli lemma (N1) is $\mathfrak m$-null.
Finally, let $G$ be the set of points of the shell at which either $x_k\in\partial{\bf D}$ and
$\langle v_k,n(x_k)\rangle=0$ for some $k$, or $x_k,x_l\in\partial{\bf D}$ for some $k\ne l$. Each of these
finitely many sets is the image of a $\cC^1$ map of $4N-3$ variables: in the first case $x_k$ is given by
its arclength and $v_k=c\,\mathbf t(x_k)$, in the second $x_k$ and $x_l$ are given by their arclengths, and in
both cases one variable is removed by the constraint $|\mathbf v|^2=N$. Therefore
$E:=\bigcup_{r\in[0,T]}\varphi^{-r}(G)$ is the image of a $\cC^1$ map of $4N-2$ variables into the shell, of
dimension $4N-1$, and $\mathfrak m(E)=0$. If a configuration is of class (N2), with non-transversal contact at
time $r$, and $q$ is a rational number between the last transversal contact (or $0$) and $r$, then the
configuration belongs to $O_q$ and its state at time $q$ belongs to $E$. Thus (N2) is contained in
$\bigcup_{q\in\mathbb Q\cap[0,T)}(\Psi^q_N)^{-1}(E)$, which is $\mathfrak m$-null by Step~2.
\end{proof}
The previous Lemma, together with the following, conclude the proof of Theorem~\ref{thm:micro}.
\begin{lem}[Empirical measures are weak solutions]\label{lem:empirical}
Let $\mathbf z^N$ be admissible on $[0,T]$. Then the $N$-particle dynamics issued from $\mathbf z^N$ is unique
on $[0,T]$, and $\{\delta_{N,t}\}_{t\in[0,T]}$ is a weak solution with initial datum $\delta_{N,0}$, in the sense
of Definition~\ref{def:w-sol}.
\end{lem}
\begin{proof}
Between two contacts the particles solve the smooth system \eqref{eq:equation of motion2} on the shell,
and at a contact the reflection is prescribed; uniqueness follows by induction on the finitely many
contacts. For $\varphi\in\cC_b(\Omega)$ the function $t\mapsto\int\varphi\,d\delta_{N,t}=\frac1N\sum_i\varphi(z_i(t))$
is piecewise continuous, and it is continuous if $\varphi$ is symmetric, since at a contact only the velocity
of the colliding particle changes, from $v$ to $\cR_xv$. This is {\rm(i)}, and {\rm(ii)} holds because
$\int|v|^2d\delta_{N,t}=1$. For {\rm(iii)}, let $\varphi\in\mathcal D_{\mathrm{sym}}$. For each $i$ the
function $g_i(t):=\varphi(t,z_i(t))$ is continuous on $[0,T]$, by the symmetry of $\varphi$, and piecewise
$\cC^1$, with $g_i(T)=0$ and
$\dot g_i=\partial_t\varphi+\langle\nabla_x\varphi,v_i\rangle+\langle\nabla_v\varphi,E_0-\alpha_Nv_i\rangle$ between
contacts. Since $\alpha_N=J_N/U_N$ is the friction coefficient of $V_{\delta_{N,t}}$, summing
$-\varphi(0,z_i(0))=\int_0^T\dot g_i\,dt$ over $i$ and dividing by $N$ gives \eqref{eq:weak-form-final}.
\end{proof}
\subsection{Mean-field limit}
We are now ready to prove the last claim of Theorem \ref{thm:mfl}. Thanks to the previous section, the proof will be a consequence of the following result (which is standard, but we include the proof for completeness), and the next Theorem \ref{thm:B}.

\begin{lem}\label{lem:cmt}
Let $\nu_N\rightharpoonup\nu$ narrowly in $\cP(\Omega)$, and let $h_N,h:\Omega\to\bR^d$ be Borel.
Assume that there is a Borel set $\Omega_0$ with $\nu(\Omega_0)=1$ such that $h_N(z_N)\to h(z)$
whenever $z_N\to z\in\Omega_0$, and that $|h_N(z)|\le C(1+|v|)$ with
$\sup_N\int|v|^2d\nu_N<\infty$. Then $\int h_N\,d\nu_N\to\int h\,d\nu$.
\end{lem}

\begin{proof}
$\Omega$ is a Polish space, so by Skorokhod's representation theorem \cite[Section~6]{Bil} there exist random
variables $Z_N,Z$ on a common probability space, with laws $\nu_N,\nu$, such that $Z_N\to Z$ almost surely.
Since $Z\in\Omega_0$ almost surely, $h_N(Z_N)\to h(Z)$ almost surely. The bound on $h_N$ and the uniform
second moments make the family $\{h_N(Z_N)\}$ uniformly integrable, and Vitali's theorem gives
$\int h_N\,d\nu_N=\mathbb E\,h_N(Z_N)\to\mathbb E\,h(Z)=\int h\,d\nu$.
\end{proof}
\begin{thm}[Mean-field limit]\label{thm:B}
Let $\mu_0=f_0\,dz\in\cP(\Omega)$ satisfy \eqref{eq:mu0-normalised}, and let $\mathbf z^N$ be
configurations admissible on $[0,T]$ with $\delta_{N,0}\rightharpoonup\mu_0$ narrowly. Then:
\begin{enumerate}
\item[(i)] the sequence $\{A_N\}$ is relatively compact in $\cC^0([0,T])$, and every limit point $A$ is
the history of a Lagrangian solution $(A,\{\mu_t\})$ with initial datum $\mu_0$;
\item[(ii)] if $A_{N_k}\to A$ uniformly, then $\delta_{N_k,t}\rightharpoonup\mu_t$ narrowly in
$\cP(\Omega)$ for every $t\in[0,T]$;
\item[(iii)] if the Lagrangian solution with initial datum $\mu_0$ is unique, then
$A_N\to A$ uniformly and $\delta_{N,t}\rightharpoonup\mu_t$ for every $t\in[0,T]$, for the whole sequence.
\end{enumerate}
\end{thm}

\begin{proof}
(i) Since $A_N\in\cA_T$, which is compact, we may assume $A_{N_k}\to A\in\cA_T$ uniformly. Fix $t$ and
let $h_k(z):=\langle E_0,\Delta^{A_{N_k}}_t(z)\rangle$ if $z\in\Reg_{A_{N_k}}$ or $z$ is one of the
points $z^{N_k}_i$, and $h_k(z):=0$ otherwise; let $h:=\langle E_0,\Delta^A_t\rangle$ on
$\Omega_0:=\Reg_A$ and $h:=0$ elsewhere. Then $\mu_0(\Omega_0)=1$ by Proposition~\ref{lem:ae-regular}(a), and if
$z_k\to z\in\Reg_A$ then, by Lemma~\ref{lem:stability}, $z_k\in\Reg_{A_{N_k}}$ for $k$ large and
$h_k(z_k)\to h(z)$. The bound $|h_k|\le C(1+|v|)$ is \eqref{eq:speed-bound} and the second moments of
$\delta_{N,0}$ equal $1$ by \eqref{eq:shell}. Lemma~\ref{lem:cmt} and \eqref{eq:AN-fp} give
\[
A(t)=\lim_kA_{N_k}(t)=\lim_k\int h_k\,d\delta_{N_k,0}=\int_\Omega\langle E_0,\Delta^A_t\rangle\,d\mu_0 .
\]
This holds for every $t$, so $A$ satisfies (ii) of Lemma~\ref{lem:equivalent-fp}, and
$(A,\{\mu^A_t\})$ is a Lagrangian solution.

(ii) Let $\varphi\in\cC_b(\Omega)$ and $t\in[0,T]$. Apply Lemma~\ref{lem:cmt} to
$h_k:=\varphi\circ\Phi^t_{A_{N_k}}$ (defined as above) and $h:=\varphi\circ\Phi^t_A$, with $\Omega_0$ the set
of points of $\Reg_A$ having no contact at time $t$, which has full $\mu_0$-measure by
Proposition~\ref{lem:ae-regular}(b); the continuity assumption is again Lemma~\ref{lem:stability}. This gives
$\int\varphi\,d\delta_{N_k,t}\to\int\varphi\,d\mu_t$.

(iii) If the solution is unique, every subsequence of $\{A_N\}$ has a further subsequence converging to
the same $A$, so the whole sequence converges, and (ii) applies.
\end{proof}

\begin{oss}\label{rem:B-comments}
(a) We stress that no stability estimate, truncation in velocity or cut-off near grazing is needed for
Theorem~\ref{thm:B}:, but only compactness of $\cA_T$ and almost-everywhere continuity of the flow.
The convergence in (ii) holds for every bounded continuous test function, not only for symmetric test functions,
although the velocities of the particles jump at the collisions.


(b) By Proposition~\ref{prop:admissible} the admissibility assumption is satisfied by almost every
configuration. For instance, if $z_1,z_2,\dots$ are independent with law $\mu_0$ and the velocities of the first
$N$ points are divided by $(\frac1N\sum_{i\le N}|v_i|^2)^{1/2}$, the resulting configurations lie on the shell,
their empirical measures converge narrowly to $\mu_0$ almost surely, and their law on the shell is absolutely
continuous with respect to $\mathfrak m$.
\end{oss}

\section{Uniqueness in the fast regime}\label{sec:fast}

In this section $\mu_0=f_0\,dz$ with
\begin{equation}\label{eq:fast-data}
f_0\in L^\infty(\Omega),\qquad \operatorname{supp}f_0\subset\{c_0\le|v|\le L\},\qquad\int_\Omega|v|^2d\mu_0=1,
\end{equation}
for some $0<c_0\le L$. We denote by $\kappa_{\min}>0$ and $\kappa_{\max}$ the extreme values of the curvature of
$\partial{\bf D}$, by $d_{\min}>0$ the minimal distance between two distinct obstacles in the cover
$\bR^2$, and we set $\lambda_\pm:=e^{\pm|E_0|T}$, so that $\lambda_-\le e^{\pm A}\le\lambda_+$ for
$A\in\cA_T$. Define
\begin{equation}\label{eq:c1L1}
c_1:=(c_0+1)e^{-|E_0|T}-1,\qquad L_1:=(L+1)e^{|E_0|T}-1 .
\end{equation}
The \emph{fast regime} is the assumption
\begin{equation}\label{eq:fast-regime}
c_1>0\qquad\text{and}\qquad\kappa_{\min}\,c_1^2>2|E_0| .
\end{equation}
It says that on $[0,T]$ every particle is fast enough for its trajectory to be less curved than the obstacles.\\
The uniqueness argument will compare two characteristic flows on short time
intervals. We first establish that particles remain
in a fixed range of speeds and that two consecutive contacts are separated by a
uniform positive time. We may therefore divide $[0,T]$ into windows in which
each trajectory has at most one contact.\\
Recalling the notation introduced in \eqref{eq:def-Delta}, we have he following easy Lemma.

\begin{lem}\label{lem:speed-window}
Let $A\in\cA_T$, $z=(x,v)\in\Reg_A$ with $c_0\le|v|\le L$. Then $c_1\le|V^A_t(z)|\le L_1$ for all $t\in[0,T]$.
\end{lem}

\begin{proof}
The function $t\mapsto|V^A_t(z)|$ is continuous (reflections preserve the speed) and Lipschitz, and between
contacts $\dot V=E_0-\dot AV$ a.e.\ with $|\dot A|\le|E_0|$. Hence
$\bigl|\tfrac{d}{dt}|V|\bigr|\le|E_0|(1+|V|)$ a.e., that is
$\bigl|\tfrac{d}{dt}\log(1+|V|)\bigr|\le|E_0|$, which gives the claim.
\end{proof}
The point of the next Lemma is that, immediately after the first contact, the distance from the
obstacle is increasing. The fast-regime condition prevents its derivative
from returning to zero in a uniformly short time, and hence prevents the
trajectory from returning to the same obstacle.
\begin{lem}\label{lem:one-contact}
Under assumption \eqref{eq:fast-regime} there exists $h_0>0$, depending only on $c_1,L_1,|E_0|$ and the geometry, such
that for every $A\in\cA_T$ and every $z\in\Reg_A$ with $c_0\le|v|\le L$, any two contact times
$\tau<\tau'$ of the trajectory of $z$ satisfy $\tau'-\tau>h_0$.
\end{lem}

\begin{proof}
It suffices to consider consecutive contacts. If they occur on two distinct obstacles of the cover, the
trajectory travels a distance at least $d_{\min}$ at speed at most $L_1$
(Lemma~\ref{lem:speed-window}), so $\tau'-\tau\ge d_{\min}/L_1>d_{\min}/(2L_1)$.

Assume now that both contacts are on the same obstacle $Q$ of the cover, and let $d:=\operatorname{dist}(\cdot,Q)$,
which is of class $\cC^2$ on $\bR^2\setminus Q$. At a point at distance $d$ from $Q$, with nearest point
of curvature $\kappa$, the Hessian of $d$ has eigenvalue $0$ in the direction $\nabla d$ and
$\kappa/(1+\kappa d)\ge\kappa_{\min}/(1+\kappa_{\min}d)$ in the orthogonal direction. Let
$\rho(t):=d(\widetilde X_t)$ for $t\in[\tau,\tau']$. Since there are no contacts in $(\tau,\tau')$, the velocity is Lipschitz on
this interval and  the distance function has bounded Hessian outside the convex
obstacle, it follows that $\rho$ is continuously differentiable and $\rho'$ is Lipschitz on
$[\tau,\tau']$. In particular,
 $\rho'=\langle\nabla d,V\rangle$,
$\rho'(\tau^+)>0$ by transversality, $0\le\rho(t)\le L_1(t-\tau)$, and for a.e.\ $t$
\[
\rho''=D^2d[V,V]+\langle\nabla d,E_0-\dot A\,V\rangle
\ \ge\ \frac{\kappa_{\min}\,(|V|^2-\rho'^2)}{1+\kappa_{\min}\rho}-|E_0|\,(1+|\rho'|) .
\]
Let $h_1>0$ be such that $\kappa_{\min}c_1^2/(1+\kappa_{\min}L_1h_1)\ge\tfrac32|E_0|$, which is possible by
\eqref{eq:fast-regime}. Using $|V|\ge c_1$ and $\rho\le L_1h_1$, for a.e.\ $t\in[\tau,\min(\tau',\tau+h_1)]$ we get
\[
\rho''\ge G(\rho'),\qquad G(p):=\frac{\kappa_{\min}}{1+\kappa_{\min}L_1h_1}\,(c_1^2-p^2)-|E_0|\,(1+|p|) ,
\]
and $G(0)\ge\tfrac12|E_0|$ by the choice of $h_1$. By continuity there is $p_0>0$ with $G\ge\tfrac14|E_0|$ on
$[0,p_0]$. Suppose that $\rho'$ vanishes somewhere in $(\tau,\min(\tau',\tau+h_1)]$, and let $t_1$ be the first
such time. Since $\rho'$ is continuous and positive on $(\tau,t_1)$, there is $t_2<t_1$ with
$0<\rho'\le p_0$ on $[t_2,t_1)$; on this interval $(\rho')'=\rho''\ge\tfrac14|E_0|>0$ a.e., hence
$\rho'(t_1)>\rho'(t_2)>0$, a contradiction. Thus $\rho'>0$ and $\rho$ is increasing on
$[\tau,\min(\tau',\tau+h_1)]$, and since $\rho(\tau')=0$ we conclude $\tau'-\tau>h_1$. The lemma holds with
$h_0:=\min(h_1,d_{\min}/(2L_1))$.
\end{proof}

We now compare, on a short window, the trajectories of the same state under two histories. Fix
$t_0\in[0,T)$, $h\in(0,h_0]$ with $h\le1$, $I:=[t_0,\min(t_0+h,T)]$, and $A_0,A_1\in\cA_T$ with
\begin{equation}\label{eq:agree}
A_0=A_1\ \text{ on }[0,t_0],\qquad\varepsilon:=\sup_{I}|A_1-A_0| .
\end{equation}
Since both histories are $|E_0|$-Lipschitz, $\varepsilon\le2|E_0|h$. Let $z\in\Reg_{A_0}\cap\Reg_{A_1}$
with $c_0\le|v|\le L$ and no contact at time $t_0$. By \eqref{eq:agree} the two trajectories coincide on
$[0,t_0]$; let $\zeta=(x,w)\in\mathring{\bf D}\times\bR^2$ be their common state at time $t_0$, and denote
by $\widetilde X^i_t$, $i=0,1$, the lifted position at time $t\in I$ under $A_i$. By
Lemmas~\ref{lem:speed-window} and~\ref{lem:one-contact}, $|w|\le W:=\lambda_+L_1$ and each of the two
trajectories has at most one contact in $I$; when it exists we denote it by $\tau_i$, and by
$\nu_i:=|\langle n,w_i^-\rangle|$ the normal component of the incoming $w$-velocity.

\begin{lem}[Comparison on a window]\label{lem:comparison}
Assume \eqref{eq:fast-regime}. There exist constants $C,c_*\ge1$, depending only on $c_1,L_1,|E_0|,T$ and the
geometry, such that, with the notation above and $\eta:=c_*\sqrt{\varepsilon h}$, for every $t\in I$:
\begin{enumerate}
\item[(a)] if neither trajectory has a contact in $[t_0,t]$, then
$|\widetilde X^1_t-\widetilde X^0_t|\le C\varepsilon h$;
\item[(b)] if at least one trajectory has a contact in $[t_0,t]$, and $\nu_i\ge\eta$ for every contact
in $I$, then
\[
|\widetilde X^1_t-\widetilde X^0_t|\le C\varepsilon h\Bigl(1+\frac1{\nu_*}\Bigr),
\]
where $\nu_*$ is the normal velocity at the earlier of the contacts;
\item[(c)] in all cases $|\widetilde X^1_t-\widetilde X^0_t|\le2L_1h$.
\end{enumerate}
\end{lem}

\begin{proof}
First of all, (c) follows directly from Lemma~\ref{lem:speed-window}, so that we are left with (a) and (b).\\
The mechanism behind the proof is the following. As long as both
trajectories are in free flight, their positions differ by $O(\varepsilon
h)$. If one trajectory reaches an obstacle transversally with normal velocity
$\nu$, this error produces a delay of order
$O(\varepsilon h/\nu)$ in the contact time of the other trajectory. During
this delay the positions remain controlled, and after the reflections of both trajectories,
 the regularity of the reflection map propagates the same estimate.
The threshold $\nu\ge c_*\sqrt{\varepsilon h}$ ensures that the perturbation
is small enough for the second contact to occur with the same obstacle.\\ 

Let us start with the free flight, i.e. we prove (a). For $s\le t$ in $I$, since $|e^{\pm A_1}-e^{\pm A_0}|\le\lambda_+\varepsilon$ on $I$,
\begin{equation}\label{eq:coeff-diff}
|\alpha_{A_1}-\alpha_{A_0}|(s,t)\le\lambda_+\varepsilon h,\qquad
|\gamma_{A_1}-\gamma_{A_0}|(s,t)\le\lambda_+\varepsilon h,\qquad
|\beta_{A_1}-\beta_{A_0}|(s,t)\le2\lambda_+^2\varepsilon h^2 .
\end{equation}
Let $\gamma_i(r)$, $r\in I$, be the position along the free flight \eqref{eq:free-flight} of $A_i$
issued from $\zeta$ at time $t_0$, continued through the obstacles. By \eqref{eq:coeff-diff} and
$\dot\gamma_i=e^{-A_i}(w+\gamma_{A_i}E_0)$,
\begin{equation}\label{eq:flight-diff}
|\gamma_1(r)-\gamma_0(r)|\le C_1\varepsilon h,\qquad|\dot\gamma_1(r)-\dot\gamma_0(r)|\le C_1\varepsilon,
\qquad|\dot\gamma_i(r)|\le L_2:=\lambda_+\bigl(W+\lambda_+|E_0|\bigr),
\end{equation}
for $r\in I$; note that $L_2\ge L_1$. All the estimates below are on the window $I$, where \eqref{eq:coeff-diff} holds. This proves (a), since without contacts
$\widetilde X^i_t=\gamma_i(t)$.

We are left with (b). Exchanging the roles of $A_0$ and $A_1$ if necessary, we may assume that the
trajectory of $A_0$ has a contact at $\tau_0\in I$, on an obstacle $Q$, and that the trajectory of $A_1$ has
no contact before $\tau_0$; thus $\nu_*=\nu_0\ge\eta$. For $t\le\tau_0$ both trajectories are free flights
and (a) applies. Let $d_s$ be the signed distance to $\partial Q$ (positive outside $Q$), of class $\cC^2$ on
$\{|d_s|<r_0\}$ with $r_0<d_{\min}/2$, and $\rho_i:=d_s\circ\gamma_i$. Then $\rho_i$ is continuously differentiable and $\rho_i'$ is Lipschitz
whenever $\gamma_i$ remains in the tubular neighbourhood. In particular,
$\rho_i''$ exists almost everywhere there. Choosing
$B$ large enough,
\begin{gather*}
|\rho_i''|\le B\ \text{ a.e. where }|\rho_i|<r_0,\qquad
|\rho_1-\rho_0|\le C_1\varepsilon h,\qquad|\rho_1'-\rho_0'|\le C_2\varepsilon,\\
\frac{\lambda_-(W+\lambda_+|E_0|)\,L_2}{2B}\le\frac{r_0}2 .
\end{gather*}
Here and below $\nu_{\max}:=W+\lambda_+|E_0|$ is an upper bound for $\nu_0$, and the last condition reads
$\lambda_-\nu_{\max}L_2/(2B)\le r_0/2$.
Moreover $\rho_0(\tau_0)=0$ and $\rho_0'(\tau_0)=-e^{-A_0(\tau_0)}\nu_0\le-\lambda_-\nu_0$. Set
\[
\theta:=\frac{\lambda_-\nu_0}{2B},\qquad\theta':=\frac{4C_1\varepsilon h}{\lambda_-\nu_0} .
\]
We choose $c_*$ so large that $\nu_0\ge c_*\sqrt{\varepsilon h}$ implies
\[
\theta'\le\theta,\qquad C_2\varepsilon\le\frac{\lambda_-\nu_0}4,\qquad C_1\varepsilon h\le\frac{r_0}2 .
\]
The first holds if $c_*^2\ge8BC_1/\lambda_-^2$; the second if $c_*\ge4C_2\sqrt{2|E_0|}/\lambda_-$, because
$\varepsilon\le\sqrt{2|E_0|}\sqrt{\varepsilon h}$; the third if $c_*^2\ge2C_1\nu_{\max}^2/r_0$, because
$\varepsilon h\le\nu_0^2/c_*^2\le\nu_{\max}^2/c_*^2$. With these choices, on $[\tau_0-\theta,\tau_0+\theta]\cap I$ we have
$|\rho_0|\le L_2\theta\le r_0/2$ and $|\rho_1|\le|\rho_0|+C_1\varepsilon h\le r_0$, so both flights stay in the
tubular neighbourhood of $\partial Q$, where the bounds above apply and which does not meet any other obstacle.

\emph{The second trajectory collides shortly after.} On $[\tau_0,\tau_0+\theta]\cap I$ we have
$\rho_0'\le-\lambda_-\nu_0/2$, hence $\rho_1'\le-\lambda_-\nu_0/4$, while
$0\le\rho_1(\tau_0)\le C_1\varepsilon h$ (the trajectory of $A_1$ has not entered $Q$ before $\tau_0$).
Therefore $\rho_1$ vanishes at a time $\tau_1^*\in[\tau_0,\tau_0+\theta']$, unless $\tau_0+\theta'$ exceeds
the endpoint of $I$. On $[\tau_0,\tau_1^*]$ the flight $\gamma_1$ stays within distance $r_0$ of $Q$, hence
away from the other obstacles, so $\tau_1^*$ is the first contact of the trajectory of $A_1$. In
all cases:
\begin{equation}\label{eq:tau-diff}
\begin{gathered}
\text{for every }t\in I\text{ with }t>\tau_0:\quad\text{either }\ t-\tau_0\le\theta'\ \text{ and }A_1\text{ has no contact in }[t_0,t],\\
\text{or }\ \tau_1\in[\tau_0,\min(t,\tau_0+\theta')] .
\end{gathered}
\end{equation}

\emph{Between the two contacts.} In the first alternative of \eqref{eq:tau-diff},
\[
|\widetilde X^1_t-\widetilde X^0_t|\le|\gamma_1(t)-\gamma_1(\tau_0)|+|\gamma_1(\tau_0)-\gamma_0(\tau_0)|
+|\widetilde X^0_t-\widetilde X^0_{\tau_0}|\le2L_2\theta'+C_1\varepsilon h .
\]

\emph{After both contacts.} In the second alternative, let $y_i=\gamma_i(\tau_i)$ and
$w_i^-=w+\gamma_{A_i}(t_0,\tau_i)E_0$ be the contact points and incoming velocities. By
\eqref{eq:coeff-diff}, \eqref{eq:flight-diff} and $|\tau_1-\tau_0|\le\theta'$,
\[
|y_1-y_0|\le L_2\theta'+C_1\varepsilon h,\qquad|w_1^--w_0^-|\le\lambda_+|E_0|(\theta'+\varepsilon h),
\]
and, since $|n(y_1)-n(y_0)|\le C_{\mathrm{geo}}|y_1-y_0|$ for a constant $C_{\mathrm{geo}}$ depending only on the
geometry, the outgoing velocities
$w_i^+=\cR_{y_i}w_i^-$ satisfy $|w_1^+-w_0^+|\le C_3(\theta'+\varepsilon h)$. By Lemma~\ref{lem:one-contact}
there is no further contact in $I$, so for $t\ge\tau_1$
\[
\widetilde X^i_t=y_i+\alpha_{A_i}(\tau_i,t)\,w_i^++\beta_{A_i}(\tau_i,t)\,E_0 .
\]
By \eqref{eq:coeff-diff}, $|\alpha_{A_1}(\tau_1,t)-\alpha_{A_0}(\tau_0,t)|\le\lambda_+(\theta'+\varepsilon h)$
and similarly for $\beta$, hence $|\widetilde X^1_t-\widetilde X^0_t|\le C_4(\theta'+\varepsilon h)$.

In both alternatives the bound is $C(\theta'+\varepsilon h)\le C\varepsilon h(1+1/\nu_0)$, which is (b).
\end{proof}

We are finally in the position to state and prove the last ingredient needed to the proofs of Theorem \ref{thm:macro}(2) and the first statement of Theorem \ref{thm:mfl}, which is the following Theorem.
\begin{thm}[Uniqueness in the fast regime]\label{thm:C-fast}
Let $\mu_0$ satisfy \eqref{eq:fast-data} and let $T>0$ be such that \eqref{eq:fast-regime} holds. Then the
Lagrangian solution on $[0,T]$ with initial datum $\mu_0$ is unique.
\end{thm}

\begin{proof}
Let $A_0,A_1$ be the histories of two Lagrangian solutions and let
$t_0:=\sup\{t\in[0,T]:A_0=A_1\text{ on }[0,t]\}$. Assume by contradiction that $t_0<T$, and let $h$, $I$,
$\varepsilon$ be as in \eqref{eq:agree}. By Lemma~\ref{lem:equivalent-fp}(ii), for $t\in I$,
\begin{equation}\label{eq:deltaA-window}
|A_1(t)-A_0(t)|\le|E_0|\int_\Omega\bigl|\Delta^{A_1}_t-\Delta^{A_0}_t\bigr|\,d\mu_0
=|E_0|\int\bigl|\widetilde X^1_t-\widetilde X^0_t\bigr|\,\varrho(\zeta)\,d\zeta ,
\end{equation}
where we used that the two trajectories coincide up to time $t_0$, and $\varrho$ is the density of the
common state $\zeta$ at time $t_0$ in the variables $(x,w)$. By Lemma~\ref{lem:F1},
$\|\varrho\|_\infty\le\|f_0\|_\infty$, and $\varrho$ is supported in $\{|w|\le W\}$.

For $i=0,1$ and $\eta>0$ let $\mathcal C_i$ be the set of $\zeta$ with $|w|\le W$ whose trajectory under $A_i$,
started at time $t_0$, has a contact in $I$, and $\mathcal C_i^\eta\subset\mathcal C_i$ the subset where $\nu_i<\eta$. Applying the free-flight case of
Lemma~\ref{lem:F2} with $q=t_0$, and noting that $|w_i^-|\le W':=W+\lambda_+|E_0|$,
\begin{align}
\int_{\mathcal C_i}\frac{d\zeta}{\nu_i}
&\le\int_I\int_{\partial{\bf D}}\int_{|w|\le W'}e^{-A_i(t)}\,dw\,d\sigma\,dt\le C_5\,h, \label{eq:flux-1}\\
\Leb\bigl(\mathcal C_i^\eta\bigr)
&\le\int_I\int_{\partial{\bf D}}\int_{|w|\le W',\,|\langle n,w\rangle|<\eta}e^{-A_i(t)}|\langle n,w\rangle|\,dw\,d\sigma\,dt
\le C_5\,h\,\eta^2 . \label{eq:flux-2}
\end{align}
Take $\eta=c_*\sqrt{\varepsilon h}$ and split the integral in \eqref{eq:deltaA-window} according to
Lemma~\ref{lem:comparison}. On the set where neither trajectory has a contact in $[t_0,t]$ we use (a); on
$\mathcal C_0^\eta\cup\mathcal C_1^\eta$ we use (c) and \eqref{eq:flux-2}; on the rest we use (b), the
bound $1/\nu_*\le1/\nu_0+1/\nu_1$ and \eqref{eq:flux-1}. This gives, for every $t\in I$,
\[
|A_1(t)-A_0(t)|\le|E_0|\,\|f_0\|_\infty\Bigl[C\varepsilon h\,\pi W^2+2C\varepsilon h\,C_5h
+2L_1h\cdot2C_5\,h\,c_*^2\varepsilon h\Bigr]\le C_6\,h\,\varepsilon ,
\]
using $|{\bf D}|\le1$. Taking the supremum over $t\in I$ we obtain $\varepsilon\le C_6h\,\varepsilon$, with
$C_6$ independent of $t_0$ and $h$. Choosing $h<1/C_6$ gives $\varepsilon=0$, i.e.\ $A_0=A_1$ on
$[0,t_0+h]\cap[0,T]$, which contradicts the definition of $t_0$.
\end{proof}

\begin{cor}[Uniqueness of weak solutions]\label{cor:weak-unique}
Under the assumptions of Theorem~\ref{thm:C-fast}, the Lagrangian solution is the unique element of
$\mathcal W_T(\mu_0)$.
\end{cor}

\begin{proof}
The Lagrangian solution belongs to $\mathcal W_T(\mu_0)$ by Proposition~\ref{prop:lagr-weak}
and Proposition~\ref{lem:ae-regular}{\rm(c)}. Conversely, every element of
$\mathcal W_T(\mu_0)$ is a Lagrangian solution by Proposition~\ref{prop:weak-lagr}, and Lagrangian solutions are
unique by Theorem~\ref{thm:C-fast}.
\end{proof}

\begin{cor}[Mean-field limit in the fast regime]\label{cor:full-limit-fast}
Let $\mu_0$ and $T$ be as in Theorem~\ref{thm:C-fast}, and let $\mathbf z^N$ be configurations admissible on
$[0,T]$ with $\delta_{N,0}\rightharpoonup\mu_0$. Then $A_N\to A$ uniformly on $[0,T]$ and
$\delta_{N,t}\rightharpoonup\mu_t$ narrowly in $\cP(\Omega)$ for every $t\in[0,T]$, where $(A,\{\mu_t\})$ is the
unique Lagrangian solution.
\end{cor}
\begin{proof}
By Theorem~\ref{thm:C-fast}, the Lagrangian solution with initial datum
$\mu_0$ is unique. The conclusion therefore follows directly from
Theorem~\ref{thm:B}(iii).
\end{proof}

\begin{oss}\label{rem:fast-comments}
(a) The proof uses neither the differential of the flow, nor hyperbolicity, nor any transversality of
singularity manifolds. It rests on three facts: by causality only short windows matter; on a short window
each particle has at most one contact; and for one contact the singular factor $1/\nu$ is integrable
\emph{exactly}, because the flux measure carries the factor $\nu$.

(b) Slow particles as such are not the difficulty: the comparison lemma only requires at most one contact
per window. Assumption \eqref{eq:fast-regime} is used in Lemma~\ref{lem:one-contact} alone, to exclude a
quick return to the same obstacle, and the bound $|v|\le L$ is used to exclude many contacts of a very fast
particle within one window. This bound is removed in Section~\ref{sec:nobound}, under a finite horizon assumption and an
exponential tail of the initial density.

(c) The particles of the $N$-particle system need not satisfy any speed bound in
Corollary~\ref{cor:full-limit-fast}.
\end{oss}

\begin{oss}[Beyond the fast regime]\label{rem:global-uniqueness}
We expect the fast-regime assumption in
Theorem~\ref{thm:C-fast} to be technical. More precisely, we conjecture that
if $\mu_0=f_0\,dz$, with $f_0\in L^\infty(\Omega)$ compactly supported in
velocity and
$
\int_\Omega |v|^2\,d\mu_0=1,
$
then the Lagrangian solution is unique on every finite time interval.

The main obstruction to proving this with the present argument is the
possible occurrence of \emph{hopping}. A particle may leave an obstacle with
a very small normal velocity on a portion of the boundary where the field
pushes it back towards the obstacle. It can then return after a very short
time and undergo many contacts in a single window. Applying the
one-contact comparison estimate separately at every contact produces a
nonintegrable loss involving the inverse normal velocities.

This loss appears to result from treating the contacts separately. If the
reflections are unfolded in the normal direction, an entire hopping episode
becomes a continuous oscillatory motion across the boundary. This suggests
that the effect of perturbing the history should be controlled by the total
duration of the episode, and not by the number of contacts. A uniform
Lipschitz estimate for the complete hopping scenario would allow the argument of Theorem~\ref{thm:C-fast} to be extended
beyond the fast regime.
\end{oss}


\section{Uniqueness without an upper bound on the velocities}\label{sec:nobound}

In Theorem~\ref{thm:C-fast} the velocity support of the initial datum is contained in $\{c_0\le|v|\le L\}$. In this and the
following sections we remove the upper bound $L$, at the price of two assumptions: the table has finite horizon, and the
initial density has an exponential tail (see Theorem \ref{thm:C-tail} below). We recall the notation of Section~\ref{sec:fast}: the curvatures
$0<\kappa_{\min}\le\kappa\le\kappa_{\max}$, the distance $d_{\min}$, the constants $\lambda_\pm$ and $c_1$, and the fast regime
\eqref{eq:fast-regime}; also, recall that $c_T=|E_0|Te^{|E_0|T}$ is the constant in \eqref{eq:speed-bound}. We denote by $D_{\max}$ the maximal
diameter of the obstacles, and we add the following assumption, standard for periodic Sinai billiards (see e.g \cite{CheMa, BaDeLi}). By periodicity and compactness,
Lemma~\ref{lem:FH} below turns it into a uniform bound on the free-path
length.

\smallskip
\noindent\textbf{(FH) Finite horizon.} Every ray $\{x+se:\ s\ge0\}\subset\bR^2$, $|e|=1$, meets the interior of an
obstacle of the cover.\\

The scope of this section is to prove the following theorem.
\begin{thm}\label{thm:C-tail}
Assume {\rm(FH)}. Let $\mu_0=f_0\,dz$ with $\int_\Omega|v|^2d\mu_0=1$, $\operatorname{supp}f_0\subset\{|v|\ge c_0\}$,
and
\begin{equation}\label{eq:tail}
f_0(x,v)\le C_0\,e^{-\alpha|v|}\qquad\text{for a.e. }(x,v)\in\Omega ,
\end{equation}
for some $C_0,\alpha>0$. Let $T>0$ satisfy \eqref{eq:fast-regime}. Then the Lagrangian solution on $[0,T]$ with initial datum
$\mu_0$ is unique.
\end{thm}
Once Theorem \ref{thm:C-tail} is proved, the following direct Corollary gives Theorem \ref{thm:macro}(3) and the first statement of Theorem \ref{thm:mfl}, thus concluding all the proofs of the paper.
\begin{cor}\label{cor:tail}
Under the assumptions of Theorem~\ref{thm:C-tail}, the Lagrangian solution is the unique element of $\mathcal W_T(\mu_0)$.
If $\mathbf z^N$ are configurations admissible on $[0,T]$ with $\delta_{N,0}\rightharpoonup\mu_0$, then $A_N\to A$ uniformly on
$[0,T]$ and $\delta_{N,t}\rightharpoonup\mu_t$ narrowly for every $t\in[0,T]$, where $(A,\{\mu_t\})$ is the unique Lagrangian
solution.
\end{cor}
\begin{proof}
By Proposition~\ref{prop:lagr-weak} and
Proposition~\ref{lem:ae-regular}{\rm(c)}, every Lagrangian solution belongs
to $\mathcal W_T(\mu_0)$. Conversely, Proposition~\ref{prop:weak-lagr}
shows that every element of $\mathcal W_T(\mu_0)$ is a Lagrangian
solution. The first conclusion therefore follows from
Theorem~\ref{thm:C-tail}. The mean-field conclusions then follow directly
from Theorem~\ref{thm:B}{\rm(iii)}.
\end{proof}
\smallskip
The proof of Theorem \ref{thm:C-tail}  is given in Section~\ref{subsec:proof-tail} after several preliminary steps, which follow the following startegy, outlined below for the convenience of the reader.\\

\emph{Outline of the proof of Theorem \ref{thm:C-tail}}. As in the proof of Theorem~\ref{thm:C-fast}, let $A_0,A_1$ be the histories of two Lagrangian solutions,
$t_0$ the last time up to which they agree, $I=[t_0,t_0+h]$ and $\varepsilon=\sup_I|A_1-A_0|$. For $t\in I$,
\begin{equation}\label{eq:basic}
|A_1(t)-A_0(t)|\le|E_0|\int\bigl|\widetilde X^1_t(\zeta)-\widetilde X^0_t(\zeta)\bigr|\,\varrho(\zeta)\,d\zeta ,
\end{equation}
where $\zeta=(x,w)$ is the common state at time $t_0$, $\varrho$ its density and $\widetilde X^i_t(\zeta)$ the lifted position at
time $t$ under $A_i$. We split the initial states according to
$$
|w|<W_0
\qquad\text{and}\qquad
|w|\ge W_0,
$$
which we call respectively the bounded-velocity and fast-velocity parts. Here $W_0$ is a velocity threshold depending only on $\Dd$, whose precise
choice will be made in \eqref{eq:choice-W0}.
The bounded-velocity part is controlled by
$C(W_0)\,h\,\varepsilon$ exactly as in Section~\ref{sec:fast}. The fast part rests on three propositions:
\begin{enumerate}
\item[(A)] (\emph{Geometry of fast particles}, Proposition~\ref{prop:fast}) a particle with $|w|\ge W_0$ moves along nearly straight
lines, its consecutive contacts are on distinct obstacles, and it has at most $1+C|w|h$ contacts in $I$;
\item[(B)] (\emph{Product lemma}, Proposition~\ref{prop:product}) for every history,
\[
\int_{\{V\le|w|\le V+1\}}\prod_{k=1}^M\frac1{\cos\varphi_k}\,d\zeta\le C(V+1)\Lambda^M ;
\]
\item[(C)] (\emph{Comparison for fast particles}, Proposition~\ref{prop:fastcomp})
$\int_{\{|w|\ge W_0\}}|\widetilde X^1_t-\widetilde X^0_t|\,\varrho\,d\zeta\le C\,h\,\varepsilon$.
\end{enumerate}
(C) follows from (A), (B) and \eqref{eq:tail} by interpolating the histories, $A_\lambda=A_0+\lambda(A_1-A_0)$: at points where
$\lambda\mapsto\widetilde X^{A_\lambda}_t(\zeta)$ is differentiable,
$$|\partial_\lambda\widetilde X^{A_\lambda}_t|\le C|w|\,h\,\varepsilon\,e^{C|w|h}\prod_k\dfrac{C_L}{\cos\varphi_k},$$ and this map is absolutely
continuous for almost every $\zeta$. Since the number of contacts is $O(|w|h)$, the loss $e^{C|w|h}$ is absorbed by the tail
\eqref{eq:tail} when $h$ is small compared with $\alpha$.
\begin{oss}[The tail assumption]\label{rem:tail}
Condition \eqref{eq:tail} is the weighted $L^\infty$ bound $\|e^{\alpha|v|}f_0\|_{L^\infty}<\infty$ and it is satisfied by Maxwellian data\footnote{That is a density with the form $f_0(x,v)=\rho_0(x)e^{-\beta |v|^2}$} and
by any superexponential velocity tail. It is also compatible with the
weak-field stationary speed distribution predicted for the thermostatted
Lorentz gas in \cite{BoCh2}, whose large-$N$ limit behaves as
$e^{-c|v|^3}$. It is much weaker than a compact velocity support, but it is stronger than the integral condition
$\int e^{\alpha|v|}d\mu_0<\infty$ becuase the products of the grazing factors are bounded with respect to the Lebesgue measure on speed shells, and we need a pointwise bound on the density in each shell. More precisely, the proof below uses the following: with
\[
b_j:=\operatorname*{ess\,sup}_{x\in{\bf D},\ j\le|v|<j+1}f_0(x,v),\qquad j\ge0,
\]
it suffices that $$\sum_{j\ge0}(j+1)^2e^{\beta j}b_j<\infty$$ for some $\beta>0$. Indeed, along every trajectory
$\bigl||W_t|-|v|\bigr|\le c_T$, so the density at time $t_0$ on the shell $\{j\le|w|<j+1\}$ is at most
$\max_{|i-j|\le c_T+1}b_i$, and this is the only requirement on $f_0$ used in Proposition~\ref{prop:fastcomp} (besides $\sup_jb_j<\infty$,
used for the slow particles). For a polynomial tail, or a stretched exponential tail $e^{-\alpha|v|^\gamma}$ with $0<\gamma<1$, the
method fails for every window length $h>0$, since a particle of speed $V$ has $O(Vh)$ contacts in the window and the averaged
product of the grazing factors grows like $e^{CVh}$. Nevertheless we expect uniqueness to be true for such data as well.
\end{oss}

\smallskip
\emph{Notational conventions.} 
\begin{itemize}
\item[-] We write $
\Dd:=({\bf D},|E_0|,T),$ where ${\bf D}$ denotes the entire periodic billiard table, including the
number, shapes and relative positions of the obstacles in a fundamental
cell. Unless stated otherwise, constants depend only on $\Dd$. In
particular, they are uniform with respect to the history $A\in\cA_T$, the
initial time $t_0$, the window length $h$, the velocity, and the initial
state. Any dependence on the lower velocity bound $c_0$ or on the tail
parameters $C_0,\alpha$ will be indicated explicitly.
\item[-] Every $A\in\cA_T$ is extended to $[0,T+2]$ by $A(t)=A(T)$ for $t\ge T$ and the extension is
$|E_0|$-Lipschitz with $\lambda_-\le e^{\pm A}\le\lambda_+$. All constructions and lemmas of Section~\ref{sec:flow} hold
on $[0,T+2]$ and for trajectories issued at any time $t_0$ instead of $0$, with the same proofs.
 For $\zeta=(x,w)\in\mathring{\bf D}\times\bR^2$ we denote by $(\widetilde X^A_t(\zeta),W^A_t(\zeta))$ the lifted state at
time $t\ge t_0$ of the trajectory issued from $\zeta$ at time $t_0$, which will be clear from the context. 
\item[-] For a unit vector $u$ we denote by $u^\perp$ the vector obtained by rotating $u$ by $\pi/2$. At a point $y\in\partial Q$ we set
$\mathbf t(y):=n(y)^\perp$ and parametrise $\partial Q$ by arclength $r$ in the direction of $\mathbf t$, so that the
angle of $n$ increases with $r$ at rate $\kappa$. 
\item[-] Finally, we write $E_*:=\lambda_+^2|E_0|$.
\end{itemize}

\subsection{Geometry of fast particles}\label{sec:prelim}

Assumption \rm{(FH)} is formulated for straight rays, whereas
the auxiliary free flights of high-velocity particles are only
approximately straight. We thus show that every sufficiently long
straight segment contains a point lying at positive distance inside an
obstacle. Consequently, every curve sufficiently close to that segment
must intersect the interior of an obstacle, and the associated physical
trajectory must have a boundary contact before that happens. This is a kind of roboustness
property of assumption {\rm(FH)}.

\begin{lem}\label{lem:FH}
Under {\rm(FH)} there are $\ell_H,\delta_H>0$ such that every segment $\{x+se:\ 0\le s\le\ell_H\}$, $|e|=1$, contains
a point $q$ with $B(q,\delta_H)$ contained in an obstacle of the cover.
\end{lem}

\begin{proof}
Let $g(q)$ be the distance from $q$ to the complement of the union of the obstacles; $g$ is continuous and
$\bZ^2$-periodic. For every $(x,e)$ there is $s(x,e)\ge0$ with $g(x+s(x,e)e)>0$, by (FH). By continuity, $g(x'+s(x,e)e')>
\tfrac12g(x+s(x,e)e)$ for $(x',e')$ in a neighbourhood of $(x,e)$. Finitely many such neighbourhoods cover the compact
set $\bT^2\times S^1$; take $\ell_H$ the maximum of the corresponding $s(x,e)$ and $\delta_H$ the minimum of
$\tfrac12g(x+s(x,e)e)$.
\end{proof}

The flow constructed in Section~\ref{sec:flow} is defined only for regular
initial states. In the analysis below we must also follow families of
trajectories through isolated tangencies. We therefore introduce an
extended trajectory by continuing a tangential flight without changing
its velocity. Proposition~\ref{prop:fast} will show that, for sufficiently
large velocities, this recursive construction is well defined and has
only finitely many contacts.

\begin{defn}[Trajectories of fast states]\label{def:fast-traj}
Let $A\in\cA_T$, $t_0\in[0,T]$ and $\zeta=(x,w)$ with $x\in\mathring{\bf D}$. On $[t_0,t_0+1]$ the trajectory issued
from $\zeta$ at time $t_0$ follows the free flights \eqref{eq:free-flight}; at a contact time $\tau$, with incoming
state $(y,w^-)$, we set $w^+=\cR_yw^-$ if the contact is transversal and $w^+=w^-$ if it is tangential, and restart
the free flight. By Proposition~\ref{prop:fast} below, if $|w|\ge W_0\ge W_1$ this defines the trajectory on the whole of
$[t_0,t_0+1]$, with finitely many contacts. On $\Reg_A$ (shifted to time $t_0$) it coincides with the trajectory of Section~\ref{sec:flow}.
\end{defn}
For high-velocity particles, time is not the most convenient scale because
the number and geometry of the contacts are controlled by the distance
travelled. For a trajectory defined on $[a,b]$, its \emph{path length} is
\[
\mathcal L_{[a,b]}
:=
\int_a^b|\dot{\widetilde X}_t|\,dt
=
\int_a^b e^{-A(t)}|W_t|\,dt .
\]
The next proposition is the main result of this subsection, and it contains the geometric info used throughout
the high-velocity argument. It shows that a sufficiently fast trajectory
is nearly straight so that (FH) forces it to meet an obstacle
within a uniformly bounded path length. On the other hand, convexity prevents it from
returning immediately to the same obstacle. Therefore, consecutive
contacts are spatially separated and their number on a time window is
$O(|w|h)$.
\begin{prop}[Geometry of fast particles]\label{prop:fast}
There is $W_{\rm geo}\ge1$ such that the following holds for $W_0\ge W_{\rm geo}$, every $A\in\cA_T$, $t_0\in[0,T]$, and every
$\zeta=(x,w)$ with $x\in\mathring{\bf D}$ and $|w|=:\rho_0\ge W_0$. Along the trajectory on $[t_0,t_0+1]$, and along the
free flights issued from its states at times in $[t_0,t_0+1]$, continued through the obstacles up to time $t_0+2$:
\begin{enumerate}
\item[(a)] $\tfrac34\rho_0\le|W_t|\le\tfrac54\rho_0$;
\item[(b)] between contacts, $t\mapsto\widetilde X_t$ is a $\cC^2$ curve with curvature at most
$\epsilon(\rho_0):=2E_*/\rho_0^2\le\epsilon_0:=2E_*/W_0^2$.
\end{enumerate}
Moreover, for the contacts at times in $[t_0,t_0+1]$:
\begin{enumerate}
\item[(c)] the first contact occurs within path length $\ell_H$, and so does the next contact after any contact;
\item[(d)] after a contact on an obstacle $Q$, transversal or tangential, the distance to $Q$ along the free flight (continued
through the other obstacles) is strictly increasing for path length $L_*:=4(\ell_H+1)$; the next contact is on a different
obstacle of the cover, and consecutive contacts are separated by path length at least $d_{\min}$;
\item[(e)] at a tangential contact the trajectory touches $\partial Q$ and stays outside $Q$ near the contact point.
\end{enumerate}
In particular the trajectory is defined on $[t_0,t_0+1]$, and the number of contacts in an interval
$[t_0,t_0+h]$, $h\le1$, is at most $1+\tfrac54\lambda_+\rho_0h/d_{\min}$.
\end{prop}

\begin{proof}
(a) Between contacts $|\frac d{dt}|W||\le|\dot W|=e^{A}|E_0|\le\lambda_+|E_0|$, and $|W|$ is continuous across contacts, so
$\bigl||W_t|-\rho_0\bigr|\le2\lambda_+|E_0|\le\rho_0/4$ on $[t_0,t_0+2]$ if $W_{\rm geo}\ge8\lambda_+|E_0|$. We also take
$W_{\rm geo}\ge2L_*/\lambda_-$: since $|\dot{\widetilde X}|=e^{-A}|W|\ge\frac34\lambda_-\rho_0$, a path length $L_*$ issued from a time in
$[t_0,t_0+1]$ is covered before time $t_0+2$, so that (a) and (b) apply to all the flights considered in (c)--(e).

(b) In the time $s=\int_{t_0}^te^{-A}$ the flight satisfies $x'=w$, $w'=e^{2A}E_0$ ($'=d/ds$), with continuous
right-hand side; hence it is $\cC^2$, and its curvature is 
$$\dfrac{|\langle w',w^\perp\rangle|}{|w|^3}\le \dfrac{E_*}{|w|^2}\le
\dfrac{\tfrac{16}9E_*}{\rho_0^2}.$$

(d), (e) Let $\gamma$ be the arclength parametrisation of a free flight issued from a contact point on $Q$ (a $\cC^2$ curve
by (b)), and $\varrho(\ell):=d^s_Q(\gamma(\ell))$ the signed distance to $Q$ (positive outside), of class $\cC^2$ on
$\bR^2\setminus\mathring Q$ and near $\partial Q$. At a point outside $Q$ at distance $d$, whose nearest point on
$\partial Q$ has curvature $\kappa$, the Hessian of $d^s_Q$ has eigenvalues $0$ (along $\nabla d^s_Q$) and
$\kappa/(1+\kappa d)\ge\kappa_{\min}/(1+\kappa_{\min}d)$. Since $|\gamma'|=1$ and $|\gamma''|\le\epsilon_0$,
\[
\varrho''\ \ge\ \frac{\kappa_{\min}}{1+\kappa_{\min}\varrho}\bigl(1-\varrho'^2\bigr)-\epsilon_0 .
\]
Assume $\epsilon_0\le\kappa_{\min}/(2+2\kappa_{\min}L_*)$. At a tangential contact $\varrho=\varrho'=0$ and
$\varrho''\ge\kappa_{\min}-\epsilon_0>0$, which gives (e). After a contact on $Q$ (with $\varrho'\ge0$), suppose that
$\varrho'$ vanishes at a first $\ell_1\in(0,L_*]$ (at $\ell_1>0$ in the tangential case, where $\varrho'>0$ just
after the contact). On $[0,\ell_1]$ we have $0\le\varrho\le L_*$, so that near $\ell_1$, where $\varrho'$ is small,
$\varrho'^2\le\tfrac14$ and hence 
$$\varrho''\ge\tfrac{\kappa_{\min}}{4(1+\kappa_{\min}L_*)}>0.$$ 
Then
$\varrho'$ is increasing near $\ell_1$, which contradicts $\varrho'(\ell_1)=0<\varrho'(\ell_1-\delta)$. Hence
$\varrho$ is increasing on $[0,L_*]$, and $\gamma$ does not meet $Q$ again within path length $L_*$.

(c) Let $\gamma$ be the arclength parametrisation of a flight, issued from a point of $\mathring{\bf D}$ or from a contact
point. By (b), 
$$|\gamma(\ell)-\gamma(0)-\ell\gamma'(0)|\le\epsilon_0\ell_H^2/2<\delta_H$$
 for $\ell\le\ell_H$, if
$\epsilon_0<2\delta_H/\ell_H^2$. By Lemma~\ref{lem:FH} the segment $\{\gamma(0)+\ell\gamma'(0)\}_{\ell\le\ell_H}$
contains a point $q$ such that $B(q,\delta_H)$ is contained in an
obstacle.  Since
$
|\gamma(\ell)-\gamma(0)-\ell\gamma'(0)|<\delta_H,
$
the auxiliary continuation $\gamma$ intersects the interior of an
obstacle before path length $\ell_H$. Consequently, the corresponding
physical trajectory has a contact no later than that path length. By the
previous paragraph, after a contact on $Q$ this next contact cannot be
with $Q$ itself. This proves (c) and the rest of (d).
The bound on the number of contacts follows from (d) and from
$$\int_{t_0}^{t_0+h}e^{-A}|W_t|\,dt\le\tfrac54\lambda_+\rho_0h.\qedhere$$
\end{proof}
Finally, we will need the following Lemma on the exponential tail at time $t_0$. 
Indeed, in the uniqueness argument the two histories are compared starting from
the last time $t_0$ at which they agreed. Hence we need the
velocity tail of the common distribution at time $t_0$ and the next lemma shows that the
exponential tail is preserved, up to a constant, by volume preservation and the uniform bound on the change of speed.
\begin{lem}\label{lem:tail}
Let $\mu_0$ be as in Theorem~\ref{thm:C-tail}, $A\in\cA_T$ and $t_0\in[0,T]$. The law of the state at time $t_0$, in the
variables $(x,w)$, has a density $\varrho$ with
\[
\varrho(x,w)\le C_0\,e^{\alpha c_T}\,e^{-\alpha|w|},\qquad\|\varrho\|_\infty\le C_0 .
\]
\end{lem}

\begin{proof}
By Lemma~\ref{lem:F1}, the map $\Psi^{t_0}$ preserves $dx\,dw$, and $w=v$ at time $0$; hence
$\varrho=f_0\circ(\Psi^{t_0})^{-1}$ on the image and $0$ elsewhere. Since $|\dot w|=e^{A}|E_0|\le\lambda_+|E_0|$, the state
$(x,w)=\Psi^{t_0}(x_0,v)$ satisfies $$\bigl||w|-|v|\bigr|\le\lambda_+|E_0|T=c_T,$$ so that 
$$f_0(x_0,v)\le C_0e^{-\alpha|v|}\le
C_0e^{\alpha c_T}e^{-\alpha|w|}.$$
\end{proof}

\subsection{Linearised dynamics}\label{sec:lin} To compare the trajectories generated by two histories, we interpolate
between them and differentiate the trajectory with respect to the
interpolation parameter. Between contacts this derivative satisfies an
inhomogeneous linear equation whose source is proportional to
$\|A_1-A_0\|_\infty$. At a contact, however, both the reflected velocity
and the contact time vary, producing a jump in the derivative. The purpose
of this subsection is to compute this jump and to isolate its singular
dependence on the angle of incidence. Ultimately, each contact will
contribute at most one factor of order $(\cos\varphi)^{-1}$.\\
Throughout this subsection, let $A_0,A_1\in\cA_T$ coincide on
$[0,t_0]$, and set
\[
\delta A:=A_1-A_0,
\qquad
A_\lambda:=A_0+\lambda\delta A,
\qquad \lambda\in[0,1].
\]
By convexity, $A_\lambda\in\cA_T$. Let $h\in(0,1]$ be such that
\[
I:=[t_0,t_0+h]\subset[0,T],
\qquad
\varepsilon:=\sup_I|\delta A|.
\]
Differentiability with respect to $\lambda$ can fail when a tangency is
created or when a contact occurs exactly at the observation time $t$.
We therefore first work on the open set on which the collision pattern is
transversal and stable.
\begin{defn}[good pairs]\label{def:good}
Fix $t\in I$. A pair $(\lambda,\zeta)\in[0,1]\times(\mathring{\bf D}\times\bR^2)$ is \emph{good} if the trajectory of
$\zeta$ under $A_\lambda$ has finitely many contacts in $(t_0,t]$, all transversal, and none at time $t$.
\end{defn}
At a transversal contact, the variation does not jump simply by the
differential of the reflection map, one must also consider the
variation of the contact time. The resulting operator is the
following linear map.
For $\xi=(\xi_x,\xi_w)\in\bR^2\times\bR^2$ and an incoming state $(y,w^-)$, $y\in\partial Q$, at time $\tau$, set
$R:=\cR_y$, $n:=n(y)$, $\mathbf t:=n^\perp$, $\kappa:=\kappa(y)$ and
\begin{equation}\label{eq:L}
\begin{gathered}
\sL\xi:=\Bigl(R\xi_x,\ R\xi_w+G\xi_x-\frac{\langle n,\xi_x\rangle}{\langle n,w^-\rangle}\bigl(Gw^-+e^{2A(\tau)}(RE_0-E_0)\bigr)\Bigr),\\
G\xi_x:=-2\kappa\langle\xi_x,\mathbf t\rangle\bigl(\langle w^-,\mathbf t\rangle n+\langle w^-,n\rangle\mathbf t\bigr).
\end{gathered}
\end{equation}

\begin{lem}[Variational equation]\label{lem:var}
The set of good pairs is relatively open in $[0,1]\times(\mathring{\bf D}\times\bR^2)$, and on it the map
$(\lambda,\zeta)\mapsto(\widetilde X^{A_\lambda}_t(\zeta),W^{A_\lambda}_t(\zeta))$ is of class $\cC^1$. Along the
trajectory of a good pair, $\eta(s):=\partial_\lambda(\widetilde X^{A_\lambda}_s,W^{A_\lambda}_s)$, $s\in[t_0,t]$, satisfies
$\eta(t_0)=0$,
\begin{equation}\label{eq:var-source}
\dot\eta_x=e^{-A_\lambda}\eta_w-\delta A\,e^{-A_\lambda}W,\qquad\dot\eta_w=\delta A\,e^{A_\lambda}E_0
\end{equation}
between contacts, and $\eta(\tau^+)=\sL\eta(\tau^-)$ at each contact time $\tau$, with $A=A_\lambda$ in \eqref{eq:L}. The
derivative in the direction $\delta\zeta$ of the initial state satisfies the same equations without the terms containing
$\delta A$, with initial value $\delta\zeta$.
\end{lem}

\begin{proof}
The coefficients $\alpha_{A_\lambda},\beta_{A_\lambda},\gamma_{A_\lambda}$ of the free flights depend smoothly on
$\lambda$, because $e^{\pm A_\lambda}=e^{\pm A_0}e^{\pm\lambda\delta A}$; hence $F^{s',s}_{A_\lambda}$ is $\cC^1$ in
$(\lambda,s,s',x,w)$, and $\partial_\lambda F^{s',s}_{A_\lambda}(x,w)$, as a function of $s'$, solves
\eqref{eq:var-source} with initial value $0$ at $s'=s$. Openness and the $\cC^1$ property follow as in Lemmas~\ref{lem:stability}
and~\ref{lem:F1}, with $\lambda$ as an additional parameter: near a
good pair the map is a composition of free flights, extended reflections $\widetilde R(x,w)=(x,\cR_{\pi(x)}w)$ and
contact times, which are $\cC^1$ in $(\lambda,\zeta)$ by the implicit function theorem.

Consider one contact, at time $\tau=\tau(\lambda,\zeta)$, and write as in the proof of Lemma~\ref{lem:F1},
 $\Psi(\lambda,\zeta)=G_\lambda(\tau(\lambda,\zeta),\zeta)$ with
$G_\lambda(\sigma,\zeta):=F^{t,\sigma}_{A_\lambda}\circ\widetilde R\circ F^{\sigma,t_0}_{A_\lambda}(\zeta)$. Let
$p=(y,w^-)$ be the incoming state and $X_\tau(x,w)=(e^{-A_\lambda(\tau)}w,e^{A_\lambda(\tau)}E_0)$. Then
\begin{gather*}
\partial_\lambda\Psi=DF^{t,\tau}\Bigl[D\widetilde R(p)\,\partial_\lambda F^{\tau,t_0}(\zeta)+
\bigl(D\widetilde R(p)X_\tau(p)-X_\tau(\widetilde Rp)\bigr)\partial_\lambda\tau\Bigr]+(\partial_\lambda F^{t,\tau})(\widetilde Rp),\\
\partial_\lambda\tau=-\frac{\langle n,\pi_x\partial_\lambda F^{\tau,t_0}(\zeta)\rangle}{e^{-A_\lambda(\tau)}\langle n,w^-\rangle}.
\end{gather*}
Hence the bracket is $\sL'\xi$ with $\xi:=\partial_\lambda F^{\tau,t_0}(\zeta)$ and
\begin{equation}\label{eq:sL}
\sL':=D\widetilde R(p)-\bigl(D\widetilde R(p)X_\tau(p)-X_\tau(\widetilde Rp)\bigr)\otimes
\frac{\langle n,\pi_x\,\cdot\,\rangle}{e^{-A_\lambda(\tau)}\langle n,w^-\rangle} .
\end{equation}
We have $D\widetilde R(p)(\xi_x,\xi_w)=(\xi_x,R\xi_w+G\xi_x)$, where
$$G\xi_x=D_y(\cR_{\pi(y)}w^-)[\xi_x]=-2\langle w^-,D\nu\,\xi_x\rangle n-2\langle w^-,n\rangle D\nu\,\xi_x$$ and
$\nu=n\circ\pi$ has derivative $D\nu\,\xi_x=\kappa\langle\xi_x,\mathbf t\rangle\mathbf t$ at $y$; this is the $G$ of
\eqref{eq:L}. Moreover,
$$D\widetilde R(p)X_\tau(p)-X_\tau(\widetilde Rp)=\bigl(e^{-A_\lambda}(w^--w^+),\
e^{-A_\lambda}Gw^-+e^{A_\lambda}(RE_0-E_0)\bigr)$$ and $w^--w^+=2\langle w^-,n\rangle n$. Substituting gives
$\sL'=\sL$: the $x$-component is $\xi_x-2\langle n,\xi_x\rangle n=R\xi_x$, and the $w$-component is the one in
\eqref{eq:L}. Iterating over the contacts gives the claim. For derivatives in $\zeta$ the same computation applies with
$\partial_\lambda F$ replaced by $DF$, which has no source term.
\end{proof}

\subsubsection{The jump in the frame of the velocity}
The coordinate-free formula~\eqref{eq:L} is convenient for deriving the
linearised collision map, but it does not display its geometric structure.
We now express the variation in the frame determined by the velocity.
Below, $J$ and $Z$ are respectively the transverse and longitudinal
components of the positional variation, while $K$ and $Y$ are the
corresponding velocity variations, normalised by the speed. The ratio
$B=K/J$ is the curvature of the associated infinitesimal wave front.\\
Given a unit vector $u$ and $\rho>0$, the {components} of $\xi=(\xi_x,\xi_w)$ are
\begin{equation}\label{eq:components}
J:=\langle\xi_x,u^\perp\rangle,\qquad Z:=\langle\xi_x,u\rangle,\qquad K:=\frac{\langle\xi_w,u^\perp\rangle}{\rho},
\qquad Y:=\frac{\langle\xi_w,u\rangle}{\rho}.
\end{equation}
At a contact with incoming state $(y,w^-)$ and outgoing velocity $w^+=\cR_yw^-$ let $\rho:=|w^\pm|$,
$u^\pm:=w^\pm/\rho$, and let $\varphi\in(-\frac\pi2,\frac\pi2)$ be the angle of incidence, defined by
\begin{equation}\label{eq:phi}
u^-=-\cos\varphi\,n+\sin\varphi\,\mathbf t,\qquad u^+=\cos\varphi\,n+\sin\varphi\,\mathbf t ,
\end{equation}
so that $\cos\varphi=|\langle w^-,n\rangle|/\rho$.

\begin{lem}\label{lem:jump}
Let $(J,Z,K,Y)^-$ be the components of $\xi$ with respect to $(u^-,\rho)$ and $(J,Z,K,Y)^+$ those of $\sL\xi$ with
respect to $(u^+,\rho)$. Then, with $e_n:=e^{2A(\tau)}\langle E_0,n\rangle/\rho^2$,
\begin{align*}
J^+&=-J^-,\qquad Z^+=Z^-,\\
K^+&=-K^--\frac{2\kappa}{\cos\varphi}J^--2e_n\bigl(J^-\sin\varphi+Z^-\cos\varphi\bigr)\tan\varphi,\\
Y^+&=Y^-+2e_n\bigl(J^-\sin\varphi+Z^-\cos\varphi\bigr).
\end{align*}
Consequently, if $J^-\ne0$, the ratio $B:=K/J$ satisfies
\begin{equation}\label{eq:B-jump}
B^+=B^-+\frac{2\kappa+2e_n\sin^2\varphi}{\cos\varphi}+2e_n\sin\varphi\,\frac{Z^-}{J^-}.
\end{equation}
\end{lem}

\begin{proof}
Rotation by $\pi/2$ anticommutes with reflections, so $R(u^\perp)=-(Ru)^\perp$; with $u^+=Ru^-$ this gives
$$J^+=\langle R\xi_x,(u^+)^\perp\rangle=-\langle\xi_x,(u^-)^\perp\rangle=-J^- \qquad \text{and}\qquad Z^+=Z^-.$$ The term $R\xi_w$
contributes $-K^-$ to $K^+$ and $Y^-$ to $Y^+$. From \eqref{eq:phi}, using $n^\perp=\mathbf t$ and $\mathbf t^\perp=-n$,
\begin{gather*}
(u^-)^\perp=-\cos\varphi\,\mathbf t-\sin\varphi\,n,\qquad(u^+)^\perp=\cos\varphi\,\mathbf t-\sin\varphi\,n,\\
\langle w^-,n\rangle=-\rho\cos\varphi,\qquad\langle w^-,\mathbf t\rangle=\rho\sin\varphi .
\end{gather*}
Hence $G\xi_x=-2\kappa\rho\langle\xi_x,\mathbf t\rangle(\sin\varphi\,n-\cos\varphi\,\mathbf t)=2\kappa\rho\langle\xi_x,\mathbf t\rangle(u^+)^\perp$, and
$Gw^-=2\kappa\rho^2\sin\varphi\,(u^+)^\perp$. Therefore
\[
G\xi_x-\frac{\langle n,\xi_x\rangle}{\langle n,w^-\rangle}Gw^-
=\frac{2\kappa\rho}{\cos\varphi}\bigl\langle\xi_x,\cos\varphi\,\mathbf t+\sin\varphi\,n\bigr\rangle(u^+)^\perp
=-\frac{2\kappa\rho}{\cos\varphi}J^-(u^+)^\perp .
\]
Next, $RE_0-E_0=-2\langle E_0,n\rangle n$ and $\langle n,\xi_x\rangle=-J^-\sin\varphi-Z^-\cos\varphi$, so the field term
of \eqref{eq:L} equals 
$$\frac{2e^{2A}\langle E_0,n\rangle}{\rho\cos\varphi}(J^-\sin\varphi+Z^-\cos\varphi)\,n.$$
Projecting on $(u^+)^\perp$ and $u^+$, with $\langle n,(u^+)^\perp\rangle=-\sin\varphi$ and $\langle n,u^+\rangle=\cos\varphi$, and
dividing by $\rho$ gives the formulas for $K^+$ and $Y^+$. Finally $B^+=K^+/J^+=-K^+/J^-$.
\end{proof}
By the previous Lemma the only singularity created by a transversal contact is the grazing
factor $(\cos\varphi)^{-1}$. For high velocities the additional field
terms are small, since $e_n=O(\rho^{-2})$. The exact formula for $B$ will
be used in the cone argument of Section~\ref{sec:product}, while the 
 operator estimate below will be used to compare two histories.\\
For $E_0=0$, \eqref{eq:B-jump} is the classical relation $B^+=B^-+2\kappa/\cos\varphi$ between the curvatures of the
incoming and outgoing wave fronts.

\subsubsection{Norm bounds} For the comparison argument we do not need the full collision formula
but only an upper bound on the growth of the variation. We normalise the
velocity component by the initial speed $\rho_0$ so that free
flights produce an exponential factor in $\rho_0h$, while each contact
produces only a fixed constant and the grazing factor
$(\cos\varphi)^{-1}$.\\
For a fast state $\zeta$, $|w|=\rho_0\ge W_0$, we use the norm $\|\xi\|:=|\xi_x|+|\xi_w|/\rho_0$, and set
\[
C_L:=1+3\kappa_{\max}+3E_* .
\]

\begin{lem}\label{lem:norms}
Let $W_0\ge W_{\rm geo}$, $\zeta$ fast and $\lambda\in[0,1]$. Along the trajectory, on $[t_0,t_0+1]$:
\begin{enumerate}
\item[(a)] if $\xi$ solves the homogeneous equation $\dot\xi_x=e^{-A_\lambda}\xi_w$, $\dot\xi_w=0$ on $[s,s']$, then
$\|\xi(s')\|\le e^{\lambda_+\rho_0(s'-s)}\|\xi(s)\|$;
\item[(b)] at a transversal contact with angle $\varphi$, $\|\sL\xi\|\le C_L\|\xi\|/\cos\varphi$;
\item[(c)] the source term $S:=(-\delta A\,e^{-A_\lambda}W,\ \delta A\,e^{A_\lambda}E_0)$ in \eqref{eq:var-source} satisfies
$\|S(s)\|\le3\lambda_+\rho_0\,\varepsilon$ for $s\in I$, if $W_0\ge|E_0|$.
\end{enumerate}
\end{lem}

\begin{proof}
(a) The claim follows directly from 
$$\frac{d}{ds}|\xi_x|\le\lambda_+|\xi_w|=\lambda_+\dfrac{\rho_0 |\xi_w|}{\rho_0}$$
and the fact that $|\xi_w|$ is constant.\\
(b) By Lemma~\ref{lem:jump} and its proof, $\sL\xi=(R\xi_x,R\xi_w+q)$ with
$$|q|\le\frac{2\kappa\rho}{\cos\varphi}|J^-|+\frac{2E_*}{\rho\cos\varphi}|\langle n,\xi_x\rangle|\le
\bigl(2\kappa_{\max}\rho+\frac{2E_*}{\rho}\bigr)\frac{|\xi_x|}{\cos\varphi}.$$
 By Proposition~\ref{prop:fast}(a),
$\tfrac34\rho_0\le\rho\le\tfrac54\rho_0$, so, recalling that $\rho_0\ge1$ $$|q|/\rho_0\le(3\kappa_{\max}+3E_*)|\xi_x|/\cos\varphi$$
and 
$$\|\sL\xi\|\le\|\xi\|+|q|/\rho_0\le C_L\|\xi\|/\cos\varphi.$$
(c) 
On the one hand, 
$$|\delta A\,e^{-A_\lambda}W|\le\varepsilon\lambda_+\tfrac54\rho_0.$$
 On the other hand, since $\rho_0\ge W_0\ge\max(1,|E_0|)$,
\[
\frac{|\delta A\,e^{A_\lambda}E_0|}{\rho_0}
\le\varepsilon\lambda_+.
\]
Consequently,
\[
\|S(s)\|
\le
\frac54\lambda_+\rho_0\varepsilon+\lambda_+\varepsilon
\le
\frac94\lambda_+\rho_0\varepsilon
\le
3\lambda_+\rho_0\varepsilon .
\]
\end{proof}

\subsection{The product lemma}\label{sec:product}
The derivative estimate of the preceding subsection contains one factor
$(\cos\varphi_k)^{-1}$ for each contact that cannot be controlled
pointwise, since a trajectory may come arbitrarily close to grazing.
The purpose of this subsection is to show, in the next Proposition, that their product is nevertheless
integrable with respect to the initial state, with a loss which is only
exponential in the number of contacts.

The proof is based on one-dimensional families of trajectories, which we call
\emph{fronts}, transverse to the initial velocity. Along each smooth piece
of a front, the product of the grazing factors is controlled by the
transverse expansion of the front. The last expansion factor is then
integrated by changing variables from the initial parameter to the contact
point. Tangencies divide a front into several smooth pieces, and the main
remaining task is to show that the number of pieces grows at most
exponentially with the number of contacts.

\begin{prop}[Product lemma]\label{prop:product}
Assume {\rm(FH)}. There are $W_{\rm prod}\ge W_{\rm geo}$ and $C_P,\Lambda_P\ge1$, depending only on $\Dd$, such that for every
$A\in\cA_T$, $t_0\in[0,T]$, $V\ge W_{\rm prod}$ and every integer $M$ with $1\le M\le M_V:=\lfloor\lambda_-V/(2\ell_H)\rfloor$,
\[
\int_{\{(x,w):\ x\in{\bf D},\ V\le|w|\le V+1\}}\ \prod_{k=1}^M\frac1{\cos\varphi_k(\zeta)}\,d\zeta\ \le\ C_P\,(1+V)\,\Lambda_P^M ,
\]
where $\varphi_k(\zeta)$ is the angle of the $k$-th contact of the trajectory issued from $\zeta$ at time $t_0$. For almost
every $\zeta$ in the domain, the first $M$ contacts are transversal and occur before time $t_0+\tfrac23$.
\end{prop}
By Proposition~\ref{prop:fast}, the path length up to the $M$-th contact
is at most $M\ell_H$. Since
\[
M\ell_H\le\frac{\lambda_-V}{2},
\qquad
|\dot{\widetilde X}|\ge\frac34\lambda_-V,
\]
the $M$-th contact occurs no later than
\[
t_0+\frac{M\ell_H}{(3/4)\lambda_-V}
\le t_0+\frac23.
\]
Moreover, Lemma~\ref{lem:pieces}, followed by Fubini's theorem in the
variables $(\rho_0,\theta,a)$, shows that the first $M$ contacts are
transversal for almost every initial state.\\
 For $E_0=0$ the proof is based on the estimate for the billiard map, with $\Lambda_P$ given in terms of
$d_{\min}$, $\kappa_{\min}$ and the number of obstacles visible from a given one.

\subsubsection{Set-up}
Fix $A$, $t_0$, $V\ge W_{\rm prod}$ and $M\le M_V$. We use the time $s:=\int_{t_0}^te^{-A}$, in which, between contacts,
\[
x'=w,\qquad w'=E(s):=e^{2A(t(s))}E_0,\qquad|E|\le E_* ;
\]
$E$ is continuous, and the contacts, their angles and the reflection law are the same as in the time $t$. We write $\rho:=|w|$,
$u:=w/\rho$, $\ell$ for the path length ($d\ell=\rho\,ds$), and
\[
\epsilon:=2E_*/V^2,\qquad e_\perp:=\langle E,u^\perp\rangle/\rho^2,\qquad e_\parallel:=\langle E,u\rangle/\rho^2 .
\]
By Proposition~\ref{prop:fast}(a), $\tfrac34V\le\rho\le\tfrac54(V+1)\le2V$ before time $t_0+2$, so $|e_\perp|,|e_\parallel|\le\epsilon$, the
coefficient $e_n$ of Lemma~\ref{lem:jump} satisfies $|e_n|\le\epsilon$, $E_*/\rho^2\le\epsilon$, and the curvature of the flights is at
most $\epsilon$. Every smallness condition on $\epsilon$ below is a condition $V\ge W_{\rm prod}$.\\
We disintegrate the shell $\{V\le |w|\le V+1\}$ into one-dimensional
families of initial states obtained by varying the initial position in the
direction transverse to the velocity.

\emph{Initial fronts.}
Fix
\[
\rho_0\in[V,V+1],\qquad \theta\in[0,2\pi),\qquad a\in\bR,
\]
and set
\[
x_p:=a\,e_\theta+p\,e_\theta^\perp,\qquad
w_0:=\rho_0e_\theta,\qquad
\zeta_p:=(x_p,w_0).
\]
We represent the position in a fixed unit-square fundamental domain of
$\bT^2$. The set of parameters $p$ for which
$x_p\in\mathring{\bf D}$ is a union of finitely many open intervals. Each
such interval is called an \emph{initial front} and is denoted by $I_0$.
There is a geometric constant $N_\Sigma$ such that their number is at most
$N_\Sigma$: one may take $N_\Sigma-1$ to be the maximal number of
obstacles of the cover intersected by a segment of length $\sqrt2$.

These coordinates satisfy
\[
dx\,dw=\rho_0\,d\rho_0\,d\theta\,da\,dp.
\]

For $p\in I_0$, we denote the data of the $k$-th contact by
\[
s_k(p),\quad y_k(p),\quad Q_k(p),\quad
\varphi_k(p),\quad \rho_k(p),\quad u_k^\pm(p),
\]
corresponding respectively to its time, contact point, obstacle, angle,
speed, and incoming and outgoing directions. Whenever the obstacle
$Q_k$ is fixed, we denote by $r_k(p)$ a continuous lift to $\bR$ of the
arclength coordinate of $y_k(p)\in\partial Q_k$.

\emph{Smooth pieces.}
Tangencies are precisely the parameters at which the collision pattern
may change. For $1\le k\le M$, a \emph{piece of level $k$} is therefore
a connected component of the open set
\[
\{p\in I_0:\text{the first $k$ contacts of $\zeta_p$ are transversal}\}.
\]
We regard $I_0$ itself as the unique piece of level $0$. On every piece
of level $k$, the obstacles $Q_1,\dots,Q_k$ are fixed and all the contact
data above depend in a $\cC^1$ way on $p$, by Lemma~\ref{lem:var}.
Moreover, every piece of level $k+1$ is contained in a piece of level
$k$. Thus every parameter whose first $M$ contacts are transversal
determines a nested sequence
\[
I_0=P_0(p)\supset P_1(p)\supset\cdots\supset P_M(p),
\]
where $P_k(p)$ is the piece of level $k$ containing $p$.

\emph{Variation along a front.}
Varying $p$ moves the initial position at unit speed in the direction
$e_\theta^\perp$, while leaving the initial velocity fixed. Accordingly,
let $\xi(s;p)=(\xi_x(s;p),\xi_w(s;p))$ be the solution of the homogeneous
variational equation along the trajectory of $\zeta_p$, with
$
\xi(0;p)=(e_\theta^\perp,0).
$
Thus
\[
\xi_x'=\xi_w,\qquad \xi_w'=0
\]
between contacts, while $\xi$ is transformed by $\sL$ (recall \eqref{eq:sL}) at each contact.
We denote its components in the moving frame \eqref{eq:components} by
$
(J,Z,K,Y)(s;p),
$
and write $J_k,Z_k,K_k,Y_k$ for their incoming values at the $k$-th
contact. Whenever $J_k\ne0$, we set
\[
B_k:=\frac{K_k}{J_k}.
\]

More precisely, on a piece of level $k$, let
$(\widetilde x(s,p),\widetilde w(s,p))$ be the free flight issued from
the $(k-1)$-st outgoing state, or from $\zeta_p$ when $k=1$, and continued
past the $k$-th contact while ignoring that reflection. Then, by Lemma~\ref{lem:var},
\[
\partial_p(\widetilde x,\widetilde w)(s,p)=\xi(s;p)
\]
between the $(k-1)$-st and the $k$-th contacts, including the incoming
one-sided value at the $k$-th contact.

\subsubsection{The cone}

\begin{lem}[Flight equations]\label{lem:flight}
Between contacts, as functions of the path length,
\[
\frac{dJ}{d\ell}=K-e_\perp Z,\qquad\frac{dZ}{d\ell}=Y+e_\perp J,\qquad\frac{dK}{d\ell}=-e_\perp Y-e_\parallel K,\qquad
\frac{dY}{d\ell}=e_\perp K-e_\parallel Y .
\]
\end{lem}

\begin{proof}
In the time $s$: $\xi_x'=\xi_w$, $\xi_w'=0$, $u'=(\langle E,u^\perp\rangle/\rho)\,u^\perp$, $(u^\perp)'=-(\langle E,u^\perp\rangle/\rho)\,u$ and
$\rho'=\langle E,u\rangle$. Differentiating \eqref{eq:components} gives $J'=\rho K-\rho e_\perp Z$, $Z'=\rho Y+\rho e_\perp J$,
$K'=-\rho e_\perp Y-\rho e_\parallel K$, $Y'=\rho e_\perp K-\rho e_\parallel Y$ and dividing by $\rho=d\ell/ds$ gives the claim.
\end{proof}

Thus the flight equations are a perturbation of size $O(\epsilon)=O(V^{-2})$ of the Jacobi equations of the billiard,
$J'=K$, $K'=0$, and the jump of Lemma~\ref{lem:jump} is a perturbation of size $O(\epsilon)$ of the classical jump. Let
\[
c_*:=\tfrac12d_{\min}\kappa_{\min} .
\]

\begin{lem}[Cone]\label{lem:cone}
There are $\epsilon_*>0$ and $C_*\ge1$, depending only on $\Dd$, such that the following holds if $\epsilon\le\epsilon_*$. Let $p\in I_0$
have its first $k$ contacts transversal, $k\le M$. Then $J(\cdot;p)$ does not vanish on $[0,s_k(p)]$, and for $1\le j\le k$:
\begin{enumerate}
\item[(i)] $B_j\ge-C_*\epsilon$, $|Z_j/J_j|+|Y_j/J_j|\le C_*\epsilon$, and $|Z_j/J_j|,|Y_j/J_j|\le1$;
\item[(ii)] just after the $j$-th contact, $B^+_j:=K/J\ge\kappa_{\min}/\cos\varphi_j$;
\item[(iii)] $|J_1|\ge\frac12$, and $|J_{j+1}|\ge c_*|J_j|/\cos\varphi_j$ if $j<k$.
\end{enumerate}
\end{lem}

\begin{proof}
Let $\epsilon_3:=4(1+\ell_H)^3\epsilon$, $\Lambda_0:=1+d_{\min}\kappa_{\min}$, $a:=2(\Lambda_0+4)/(\Lambda_0-1)$ and
$b:=2(\ell_H(a+4)+\Lambda_0)/(\Lambda_0-1)$. We assume
$\epsilon_3\le\min\{(\Lambda_0-1)/(2\Lambda_0),\ \kappa_{\min}/8,\ 1/(a+4),\ 1/b\}$, which defines $\epsilon_*$, and prove (i) in the
form $B_j\ge-2\epsilon_3$, $|Z_j/J_j|\le b\epsilon_3\le1$, $|Y_j/J_j|\le a\epsilon_3\le1$; this gives (i) with
$C_*:=4(1+\ell_H)^3(a+b+2)$.

\emph{A flight.} Consider a flight of path length $L\le\ell_H$. Let its initial components be
$v_0=(J_0,Z_0,K_0,Y_0)$, with $J_0>0$ (the case $J_0<0$ is symmetric), $K_0=B_0J_0$ with $B_0\ge0$, and $|Z_0|,|Y_0|\le J_0$. Write $v=(J,Z,K,Y)$ and
$\frac{dv}{d\ell}=(\mathsf A_0+\mathsf A_1)v$, where $\mathsf A_0$ is the free part ($J\mapsto K$, $Z\mapsto Y$) and
$\|\mathsf A_1\|\le2\epsilon$ in the max norm (Lemma~\ref{lem:flight}). The free propagator satisfies
$\|\Phi_0(\ell)\|\le1+\ell$. From $v(\ell)=\Phi_0(\ell)v_0+\int_0^\ell\Phi_0(\ell-l)\mathsf A_1v(l)\,dl$ and Gr\"onwall's lemma,
$\sup_{[0,\ell]}\|v\|\le(1+\ell)e^{2\epsilon\ell(1+\ell)}\|v_0\|$ and
\[
\|v(\ell)-\Phi_0(\ell)v_0\|\le2\epsilon\ell(1+\ell)^2e^{2\epsilon\ell(1+\ell)}\|v_0\|\le4\epsilon(1+\ell_H)^2\,\ell\,(1+B_0)J_0
\le\epsilon_3(1+\ell B_0)J_0 ,
\]
using $e^{2\epsilon\ell_H(1+\ell_H)}\le2$ and $\ell(1+B_0)\le(1+\ell_H)(1+\ell B_0)$. Since $\Phi_0(\ell)v_0=(J_0+\ell K_0,\ Z_0+\ell Y_0,\
K_0,\ Y_0)$ and $J_0+\ell K_0=(1+\ell B_0)J_0$:
\begin{align}
&(1-\epsilon_3)(1+\ell B_0)J_0\le J(\ell)\le(1+\epsilon_3)(1+\ell B_0)J_0, \label{eq:F1}\\
&B(\ell):=K(\ell)/J(\ell)\ge-2\epsilon_3, \label{eq:F2}\\
&|Z(\ell)|\le|Z_0|+\ell|Y_0|+\epsilon_3(1+\ell B_0)J_0,\qquad|Y(\ell)|\le|Y_0|+\epsilon_3(1+\ell B_0)J_0 . \label{eq:F3}
\end{align}
For \eqref{eq:F2}: $K(\ell)\ge B_0J_0-\epsilon_3(1+\ell B_0)J_0\ge-\epsilon_3(1+\ell B_0)J_0$, and if $K(\ell)<0$ then by \eqref{eq:F1}
$B(\ell)\ge-\epsilon_3/(1-\epsilon_3)\ge-2\epsilon_3$.

\emph{First flight.} Up to the first contact $\xi_w\equiv0$ and $\xi_x\equiv e_\theta^\perp$, so $K=Y=0$, $J=\langle e_\theta,u\rangle$ and
$Z=\langle e_\theta^\perp,u\rangle$. The direction $u$ turns at rate at most $\epsilon$ per unit path length and the first contact occurs
within path length $\ell_H$; hence $J_1\ge\cos(\epsilon\ell_H)\ge\frac12$, $|Z_1/J_1|\le2\epsilon\ell_H\le\epsilon_3\le b\epsilon_3$, $Y_1=0$
and $B_1=0$.

\emph{Induction.} Assume the bounds at the $j$-th contact. By Lemma~\ref{lem:jump}, $J^+=-J_j$, $|Z^+/J^+|=|Z_j/J_j|\le b\epsilon_3\le1$,
\[
|Y^+/J^+|\le|Y_j/J_j|+2\epsilon(1+|Z_j/J_j|)\le(a+4)\epsilon_3\le1 ,
\]
and, by \eqref{eq:B-jump}, $B^+_j\ge-2\epsilon_3+(2\kappa_{\min}-2\epsilon)/\cos\varphi_j-2\epsilon\ge(2\kappa_{\min}-6\epsilon_3)/\cos\varphi_j
\ge\kappa_{\min}/\cos\varphi_j$, which is (ii). The next flight has path length $\ell_j\in[d_{\min},\ell_H]$
(Proposition~\ref{prop:fast}(c),(d)), and $1+\ell_jB^+_j\ge\Lambda_0$. The flight estimates apply with $J_0=|J_j|$, $B_0=B_j^+$. By
\eqref{eq:F1}, $J$ does not vanish on the flight and
$|J_{j+1}|\ge(1-\epsilon_3)(1+\ell_jB_j^+)|J_j|\ge\frac12d_{\min}\kappa_{\min}|J_j|/\cos\varphi_j$, which is (iii). By \eqref{eq:F2},
$B_{j+1}\ge-2\epsilon_3$. By \eqref{eq:F3} and \eqref{eq:F1},
\[
\Bigl|\frac{Y_{j+1}}{J_{j+1}}\Bigr|\le\frac{(a+4)\epsilon_3/\Lambda_0+\epsilon_3}{1-\epsilon_3},\qquad
\Bigl|\frac{Z_{j+1}}{J_{j+1}}\Bigr|\le\frac{(b+\ell_H(a+4))\epsilon_3/\Lambda_0+\epsilon_3}{1-\epsilon_3}.
\]
Since $1-\epsilon_3\ge(\Lambda_0+1)/(2\Lambda_0)$, the choice of $a$ and $b$ gives
$(a+4+\Lambda_0)/\Lambda_0\le a(\Lambda_0+1)/(2\Lambda_0)$ and $(b+\ell_H(a+4)+\Lambda_0)/\Lambda_0\le b(\Lambda_0+1)/(2\Lambda_0)$, so the
bounds hold at the $(j+1)$-th contact.
\end{proof}

\subsubsection{Contact arcs}

\begin{lem}[Contact data]\label{lem:contact-data}
On a piece of level $k$, with the components just before the $k$-th contact and $\varphi=\varphi_k$, $\rho=\rho_k$,
$e_\perp=\langle E(s_k),(u_k^-)^\perp\rangle/\rho^2$, $\kappa=\kappa(y_k)$:
\begin{gather*}
\frac{dr_k}{dp}=-\frac{J_k}{\cos\varphi},\qquad\frac{ds_k}{dp}=-\frac{J_k\sin\varphi+Z_k\cos\varphi}{\rho\cos\varphi},\qquad
\frac{d\rho_k}{dp}=\rho Y_k+\langle E(s_k),u_k^-\rangle\frac{ds_k}{dp},\\
\frac{d\varphi_k}{dp}=\Bigl[\kappa+B_k\cos\varphi-e_\perp\Bigl(\sin\varphi+\frac{Z_k}{J_k}\cos\varphi\Bigr)\Bigr]\frac{dr_k}{dp}.
\end{gather*}
\end{lem}

\begin{proof}
With $\widetilde x,\widetilde w$ as above, $s_k$ is defined by $d^s_{Q_k}(\widetilde x(s_k(p),p))=0$ and
$y_k=\widetilde x(s_k(p),p)$. Differentiating, with $\xi=\xi(s_k^-;p)$ and $w^-=\rho u^-$,
$\langle n,\xi_x+w^-\tfrac{ds_k}{dp}\rangle=0$, so $\frac{ds_k}{dp}=\langle n,\xi_x\rangle/(\rho\cos\varphi)$ and, with the
expressions of $n,(u^-)^\perp$ in the proof of Lemma~\ref{lem:jump}, $\langle n,\xi_x\rangle=-J_k\sin\varphi-Z_k\cos\varphi$. The
vector $\frac{dy_k}{dp}=\xi_x+w^-\frac{ds_k}{dp}=\frac{dr_k}{dp}\mathbf t$ is tangent to $\partial Q_k$; its scalar product with
$(u^-)^\perp$ is $J_k$, and $\langle\mathbf t,(u^-)^\perp\rangle=-\cos\varphi$; this gives $dr_k/dp$. The speed is
$\rho_k=|\widetilde w(s_k(p),p)|$, with $\partial_p|\widetilde w|=\langle u,\xi_w\rangle=\rho Y$ and $\partial_s|\widetilde w|=\langle u,E\rangle$.
Finally, by \eqref{eq:phi} the angle of $u^-$ is the angle of $n$ plus $\pi-\varphi$. The angle of $n(y_k)$ has derivative
$\kappa\,dr_k/dp$; the angle of $u^-=\widetilde w/|\widetilde w|$ at $(s_k(p),p)$ has derivative
$\langle\partial_pu,u^\perp\rangle+\langle\partial_su,u^\perp\rangle\frac{ds_k}{dp}=K_k+\frac{\langle E,u^\perp\rangle}\rho\frac{ds_k}{dp}
=K_k-e_\perp\frac{J_k\sin\varphi+Z_k\cos\varphi}{\cos\varphi}$. Since $K_k=B_kJ_k=-B_k\cos\varphi\,dr_k/dp$ and
$J_k/\cos\varphi=-dr_k/dp$, the formula for $d\varphi_k/dp$ follows.
\end{proof}

\begin{lem}[Bounded arcs]\label{lem:arcs}
Assume $\epsilon\le\epsilon_*$ and $(C_*+4)\epsilon\le\kappa_{\min}/2$. On a piece $P$ of level $k\ge1$, $r_k$ is strictly monotone,
$\varphi_k$ is a strictly increasing function of $r_k$ with derivative at least $\kappa_{\min}/2$, and
\[
\int_P\frac{|J_k|}{\cos\varphi_k}\,dp=|r_k(P)|\le\frac{2\pi}{\kappa_{\min}} .
\]
\end{lem}

\begin{proof}
By Lemma~\ref{lem:cone}, $J_k$ does not vanish on $P$, so $dr_k/dp$ has a constant sign. By Lemma~\ref{lem:contact-data} and
Lemma~\ref{lem:cone}(i), $d\varphi_k/dr_k\ge\kappa_{\min}-C_*\epsilon-\epsilon(1+|Z_k/J_k|)\ge\kappa_{\min}-(C_*+2)\epsilon\ge\kappa_{\min}/2$.
Since $\varphi_k\in(-\frac\pi2,\frac\pi2)$, the length of $r_k(P)$ is at most $\pi/(\kappa_{\min}/2)$.
\end{proof}

\subsubsection{Tangencies}

A free flight is \emph{tangent to $Q$ on the side $\sigma\in\{+1,-1\}$} at a point $q$ if it touches $\partial Q$ at $q$ without
entering $\mathring Q$ near $q$, and $n(q)=-\sigma u^\perp$, $u$ being its direction at $q$ (so that $Q$ lies on the side of
$\sigma u^\perp$). By Proposition~\ref{prop:fast}(d),(e), a tangential contact is such a point.

\begin{lem}[Initial fronts]\label{lem:tan0}
For every initial front $I_0$ and every obstacle $Q$ of the cover, there are at most two $p\in I_0$ such that the trajectory of
$\zeta_p$ has its first contact tangential and on $Q$.
\end{lem}

\begin{proof} The free flights issued from the states $\zeta_p$ are translates:
\[
x_p(s)
 =a e_\theta+p e_\theta^\perp+c(s),
\qquad
c(s)
 :=s\rho_0e_\theta+\int_0^s(s-\sigma)E(\sigma)\,d\sigma .
\]
Set
$
s_c:=\frac{2\ell_H}{V}.
$
For $0\le s\le s_c$,
\[
|c'(s)-\rho_0e_\theta|
 \le E_*s
 \le \frac{2E_*\ell_H}{V}.
\]
Consequently,
\[
\rho_0s_c-\frac{E_*s_c^2}{2}
 \le \int_0^{s_c}|c'(s)|\,ds
 \le \rho_0s_c+\frac{E_*s_c^2}{2}.
\]
Since $s_c=2\ell_H/V$ and $\rho_0\in[V,V+1]$, after increasing $W_{\rm prod}$
the path length belongs to $[\ell_H,4\ell_H]$. The same estimate gives
$\langle c'(s),e_\theta\rangle>0$.

Every point $q$ of 
$$\Sigma:=\{q:\langle q-a e_\theta,e_\theta\rangle\in\langle c([0,s_c]),e_\theta\rangle\}$$ is of the form $a\,e_\theta+p(q)e_\theta^\perp+c(s)$ for
a unique $(p(q),s)$, and $q\mapsto p(q)$ is continuous. The set $S:=\{p:\ x_p([0,s_c])\cap Q\ne\emptyset\}=p(\Sigma\cap Q)$ is an
interval, as a continuous image of the convex set $\Sigma\cap Q$. If the first contact of $\zeta_{p_0}$ is tangential and on $Q$,
then $p_0\in S$, and $p_0$ is an endpoint of $S$: otherwise there are $q_1,q_2\in\Sigma\cap Q$ with $p(q_1)<p_0<p(q_2)$, and by
continuity of $p$ on the segment $[q_1,q_2]\subset\Sigma\cap Q$ there is $q$ in the open segment $(q_1,q_2)$, which lies in
$\mathring Q$ by strict convexity, with $p(q)=p_0$; then $x_{p_0}$ meets $\mathring Q$ within path length $4\ell_H$, which is impossible:
before its first contact the flight does not meet $Q$, and after it the distance to $Q$ increases (Proposition~\ref{prop:fast}(d)).$\qedhere$
\end{proof}

For the flights issued after a contact we use the fan of all flights from one point. Let
\[
d_c:=2(\ell_H+1+2D_{\max}).
\]
For $y\in\partial Q_k$, $\rho\in[\frac34V,2V]$ and a time $s$, set
$\Gamma_s(\tau):=\int_s^{s+\tau}(s+\tau-\sigma)E(\sigma)\,d\sigma$ and
\[
\Psi(z):=y+z+\Gamma_s(|z|/\rho),\qquad z\in B(0,d_c).
\]
The free flight issued at time $s$ from $(y,\rho e_\theta)$ is $s'\mapsto\Psi\bigl(\rho(s'-s)e_\theta\bigr)$. Along it the path length
$\ell$ and the parameter $d=\rho(s'-s)$ satisfy $|d\ell/dd-1|\le\epsilon d_c$. Since $E$ is continuous,
\begin{equation}\label{eq:Gamma-der}
\partial_\tau\Gamma_s(\tau)=\int_s^{s+\tau}E,\qquad\partial_\tau^2\Gamma_s(\tau)=E(s+\tau),\qquad
\partial_s\Gamma_s(\tau)=\int_s^{s+\tau}\bigl(E(\sigma)-E(s)\bigr)d\sigma,
\end{equation}
and $\partial_s\partial_\tau\Gamma_s(\tau)=E(s+\tau)-E(s)$ are continuous in $(s,\tau)$.

\begin{lem}[Fans]\label{lem:fan}
There is $\epsilon_4>0$, depending only on $\Dd$, such that if $\epsilon\le\epsilon_4$ the following holds for $y$, $\rho$, $s$ as above
and every obstacle $Q\ne Q_k$ at distance at most $\ell_H+1$ from $Q_k$.
\begin{enumerate}
\item[(a)] $\Psi$ is a $\cC^2$ diffeomorphism of $B(0,d_c)$ onto its image, which contains $Q$; $H:=\Psi^{-1}(Q)$ is a strictly
convex compact set with $0\notin H$.
\item[(b)] For each $\sigma$ there is exactly one direction $\Theta_\sigma=\Theta_\sigma(y,\rho,s)\in\bR/2\pi\bZ$ such that the flight issued
at time $s$ from $(y,\rho e_{\Theta_\sigma})$ is tangent to $Q$ on the side $\sigma$ at a parameter $d<d_c$; the tangency point is at
parameter $d_\sigma\ge d_{\min}/2$. Every tangency of such a flight with $Q$ on the side $\sigma$ within path length $\ell_H+1$ is of
this kind.
\item[(c)] $\Theta_\sigma$ is of class $\cC^1$ in $(y,\rho,s)$. With $r$ the arclength on $\partial Q_k$,
\[
\Bigl|\partial_r\Theta_\sigma+\frac{\langle e_{\Theta_\sigma}^\perp,\mathbf t(y)\rangle}{d_\sigma}\Bigr|\le\frac{16\epsilon d_c}{d_{\min}},\qquad
|\partial_r\Theta_\sigma|\le\frac4{d_{\min}},\qquad|\partial_\rho\Theta_\sigma|\le\frac{4\epsilon d_c^2}{\rho\,d_{\min}},\qquad
|\partial_s\Theta_\sigma|\le\frac{8E_*d_c}{\rho\,d_{\min}} .
\]
\end{enumerate}
\end{lem}

\begin{proof}
(a) Write $\Psi(z)=y+z+\gamma(|z|)$ with $\gamma(d):=\Gamma_s(d/\rho)$. By \eqref{eq:Gamma-der}, 
$$|\gamma'(d)|\le E_*d/\rho^2\le\epsilon d \qquad \text{and} \qquad  \gamma''(d)=E(s+d/\rho)/\rho^2,$$
and $|\gamma''|\le\epsilon$.  Since $\gamma'(d)/d\to\gamma''(0)$ as $d\to0$, $\Psi\in\cC^2$, with
$|D\Psi-I|\le\epsilon d_c$ and $|D^2\Psi|\le2\epsilon$ on $B(0,d_c)$.  Moreover, for $\epsilon d_c\le\frac14$, $\Psi$ is injective because
$$|\Psi(z)-\Psi(z')-(z-z')|\le\epsilon d_c|z-z'|$$
 and since 
 $$|\Psi(z)-y-z|\le\epsilon d_c^2/2\le d_c/8,$$ its image contains
$B(y,\frac34d_c)\supset Q$, because every point of $Q$ is within distance $\ell_H+1+2D_{\max}=d_c/2$ of $y$. The inverse
$\Phi:=\Psi^{-1}$ satisfies $|D\Phi-I|\le2\epsilon d_c$ and $|D^2\Phi|\le16\epsilon$. If $\beta=\Phi\circ c$ with $c$ a unit speed
parametrisation of $\partial Q$ (curvature $\ge\kappa_{\min}$), then 
$$\beta'=D\Phi c', \quad , \quad \beta''=D\Phi c''+D^2\Phi[c',c'],$$ and for $\epsilon$ small the
curvature of $\beta$ is at least 
$$\dfrac{(1-2\epsilon d_c)^2\kappa_{\min}-(1+2\epsilon d_c)16\epsilon}{(1+2\epsilon d_c)^3}\ge	 \dfrac{\kappa_{\min}}{2}.$$ 
Hence, denoting by $H$ the compact set bounded by the simple closed curve $\Phi(\partial Q)$, we have that $H$ is strictly convex and $0\notin H$
because $\Psi(0)=y\notin Q$, which concludes the proof of (a).\\

(b) The flight in the direction $\theta$ is the image under $\Psi$ of the ray $\{de_\theta:0\le d<d_c\}$, and $\Psi$, being a
homeomorphism close to the identity in $\cC^1$, maps tangencies to tangencies and preserves the side ($n(q)$ is within angle
$2\epsilon d_c$ of the normal to $\partial H$ at $\Phi(q)$, and the direction of the flight at $q$ is within angle $\epsilon d_c$ of
$e_\theta$). A tangency within path length $\ell_H+1$ occurs at $d\le(\ell_H+1)/(1-\epsilon d_c)<d_c$. \\
Since $H$ is compact and strictly convex, has nonempty interior, and $0\notin H$, there are unique angles $\theta_-<\theta_+$,
chosen so that $\theta_+-\theta_-<\pi$, for which
\[
\bigl\{\theta\in[\theta_-,\theta_-+2\pi):
\{de_\theta:d\ge0\}\cap H\ne\varnothing\bigr\}
=[\theta_-,\theta_+].
\]
For $\theta\in(\theta_-,\theta_+)$ the ray $\{de_\theta:d\ge0\}$
meets $\mathring H$, whereas the two endpoint rays are supporting rays
of $H$ and touch $\partial H$ at exactly one point. Moreover,
\[
H\subset
\{z:\langle z,e_{\theta_-}^{\perp}\rangle\ge0\},
\qquad
H\subset
\{z:\langle z,e_{\theta_+}^{\perp}\rangle\le0\}.
\]
Finally, since $H\subset\Phi(B(y,d_c/2))\subset B(0,d_c)$, this gives existence and uniqueness of $\Theta_\sigma$. Since $|q-y|\ge d_{\min}$
and $|\Psi(z)-y|\le(1+\epsilon d_c)|z|$, we conclude that $d_\sigma\ge d_{\min}/2$.\\

(c) Let 
$$G(y',\theta,d;\rho',s'):=d^s_Q\bigl(y'+de_\theta+\Gamma_{s'}(d/\rho')\bigr).$$
 By \eqref{eq:Gamma-der}, $G$ is $\cC^1$ in all variables,
$\cC^2$ in $d$, and $\partial_dG$ is $\cC^1$ in all variables. At $(y,\Theta_\sigma,d_\sigma;\rho,s)$ we have $G=\partial_dG=0$ and
$$\partial_d^2G\ge\kappa_{\min}(1-\epsilon d_c)^2-\epsilon\ge\kappa_{\min}/2.$$
By the implicit function theorem, near this point
$\min_dG$ is attained at a unique $d_m$, of class $\cC^1$ in $(y',\theta,\rho',s')$, and $g:=G(\cdot,d_m)$ is $\cC^1$, with the derivatives of
$G$ at $d=d_m$ given by
$$\partial_\theta g=\langle n(q),d_\sigma e_\theta^\perp\rangle  , \quad \partial_{y'}g=n(q) ,\quad  \partial_{\rho'}g=\langle n(q),\partial_\rho[\Gamma_s(d_\sigma/\rho)]\rangle ,\quad \partial_{s'}g=\langle n(q),\partial_s\Gamma_s(d_\sigma/\rho)\rangle$$ where
$q$ is the tangency point. \\
Near $(y,\Theta_\sigma;\rho,s)$, a tangency on the side $\sigma$ is equivalent to $g=0$. Since $n(q)=-\sigma u^\perp(q)$
and $u(q)$ is within angle $2\epsilon d_c$ of $e_\theta$, we have 
$$\langle n(q),e_\theta^\perp\rangle=-\sigma+\delta_1$$ and
$$\langle n(q),\mathbf t\rangle=-\sigma\langle e_\theta^\perp,\mathbf t\rangle+\delta_2$$
 with $|\delta_1|,|\delta_2|\le2\epsilon d_c$; in particular
$$|\partial_\theta g|\ge d_\sigma/2\ge d_{\min}/4.$$
 The implicit function theorem gives $\Theta_\sigma\in\cC^1$ with
$$\partial_r\Theta_\sigma=-\frac{\langle n(q),\mathbf t\rangle}{\langle n(q),d_\sigma e_\theta^\perp\rangle}.$$
Therefore, $|\partial_r\Theta_\sigma|\le4/d_{\min}$ and
$$|\partial_r\Theta_\sigma+\dfrac{\langle e_\theta^\perp,\mathbf t\rangle}{d_\sigma}|\le \dfrac{|\delta_1|+|\delta_2|}{\frac{d_\sigma}{2}}\le\frac{16\epsilon d_c}{d_{\min}}.$$
 Finally,
$$|\partial_\rho[\Gamma_s(d/\rho)]|=|\partial_\tau\Gamma_s(d/\rho)|\,d/\rho^2\le E_*d^2/\rho^3\le\epsilon d_c^2/\rho$$
and, by \eqref{eq:Gamma-der},
$$|\partial_s\Gamma_s(d/\rho)|\le2E_*d/\rho\le2E_*d_c/\rho.$$
 Dividing by $|\partial_\theta g|\ge d_{\min}/4$ gives the last two bounds.$\qedhere$
\end{proof}

\begin{lem}[Tangencies after a contact]\label{lem:tank}
Let $$C_\theta:=2+\kappa_{\max}/\kappa_{\min}+5/(\kappa_{\min}d_{\min}).$$
 There is $W_{\rm prod}\ge W_{\rm geo}$, depending only on $\Dd$, such that for
$V\ge W_{\rm prod}$ the following holds. Let $P$ be a piece of level $k\ge1$, $Q$ an obstacle and $\sigma\in\{\pm1\}$. Then there are at most
$C_\theta$ values of $p\in P$ such that the $(k+1)$-th contact of $\zeta_p$ is tangential, on $Q$, on the side $\sigma$.
\end{lem}

\begin{proof}
Choose $W_{\rm prod}$ so that, for $V\ge W_{\rm prod}$, the assumptions of Lemmas~\ref{lem:cone}, \ref{lem:arcs} and~\ref{lem:fan} hold and the error term $R$
below satisfies $|R|\le\kappa_{\min}/2$, and $C_*\epsilon\le\kappa_{\min}/4$. If the $(k+1)$-th contact is tangential on $Q$ then $Q$ is
within distance $\ell_H$ of $Q_k$ (Proposition~\ref{prop:fast}(c)), so Lemma~\ref{lem:fan} applies with $y=y_k(p)$, $\rho=\rho_k(p)$,
$s=s_k(p)$. Let $\theta(p)$ be a continuous lift of the angle of $u_k^+(p)$, and
$$h(p):=\theta(p)-\Theta_\sigma(y_k(p),\rho_k(p),s_k(p)),$$ a $\cC^1$ function on $P$ (after choosing a continuous lift of
$\Theta_\sigma$). If the $(k+1)$-th contact of $\zeta_p$ is tangential on $Q$ on the side $\sigma$, then $h(p)\in2\pi\bZ$, by
Lemma~\ref{lem:fan}(b).

By \eqref{eq:phi} the angle of $u_k^+$ is the angle of $n(y_k)$ plus $\varphi_k$, so by Lemma~\ref{lem:contact-data}
$$\frac{d\theta}{dp}=\bigl(2\kappa+B_k\cos\varphi_k-e_\perp(\sin\varphi_k+\frac{Z_k}{J_k}\cos\varphi_k)\bigr)\frac{dr_k}{dp}.$$
 By
Lemma~\ref{lem:cone}(i), 
$$|ds_k/dp|\le2|dr_k/dp|/\rho_k \quad \text{and}\quad |d\rho_k/dp|\le\rho_k(C_*+2)\epsilon\,|dr_k/dp|.$$ Let $p_0$ be such that
$h(p_0)\in2\pi\bZ$. At $p_0$ the direction $e_{\Theta_\sigma}$ is $u_k^+$, and $\langle(u_k^+)^\perp,\mathbf t\rangle=\cos\varphi_k$. By the chain
rule and Lemma~\ref{lem:fan}(c),
\[
\frac{dh}{dp}(p_0)=\frac{dr_k}{dp}(p_0)\Bigl[2\kappa+B_k\cos\varphi_k+\frac{\cos\varphi_k}{d_\sigma}+R\Bigr],
\]
\[
|R|\le\epsilon(1+|Z_k/J_k|)+\frac{16\epsilon d_c}{d_{\min}}+\frac{4\epsilon d_c^2}{d_{\min}}(C_*+2)\epsilon+\frac{16E_*d_c}{\rho_k^2d_{\min}}
\le2\epsilon+\frac{32\epsilon d_c}{d_{\min}}+\frac{4(C_*+2)\epsilon^2d_c^2}{d_{\min}} .
\]
For $V\ge W_{\rm prod}$, $|R|\le\kappa_{\min}/2$ and $B_k\ge-C_*\epsilon\ge-\kappa_{\min}/4$, so the bracket is at least $\kappa_{\min}$. Thus $dh/dp$
has the sign of $dr_k/dp$ at every point where $h\in2\pi\bZ$, and this sign is the same on all of $P$ (Lemma~\ref{lem:arcs}).
Consequently each value $2\pi m$ is taken by $h$ at most once on $P$, and the number of such values is at most
$1+\operatorname{osc}_P(h)/(2\pi)$.\footnote{Here, as usual,
$
\operatorname{osc}_P f
:=\sup_{p\in P}f(p)-\inf_{p\in P}f(p).$}
 Now, 
 $$ \operatorname{osc}_P\theta\le\kappa_{\max}|r_k(P)|+\pi$$ and, by Lemma~\ref{lem:fan}(c) and the bounds above,
$$|\frac{d}{dp}\Theta_\sigma(y_k(p),\rho_k(p),s_k(p))|\le\frac5{d_{\min}}|\frac{dr_k}{dp}|$$
Lemma~\ref{lem:arcs} gives $|r_k(P)|\le2\pi/\kappa_{\min}$ so that
 $$1+\operatorname{osc}_P(h)/(2\pi)\le\frac32+\frac{\kappa_{\max}}{\kappa_{\min}}+\frac5{\kappa_{\min}d_{\min}}\le C_\theta$$
 which conlcudes the proof.$\qedhere$
 
\end{proof}

\subsubsection{Number of pieces and proof of the product lemma}

Let $K_0$ be the maximal number of obstacles of the cover within distance $\ell_H$ of a segment of length $\sqrt2$, and $K$ the
maximal number of obstacles within distance $\ell_H$ of a given obstacle.

\begin{lem}[Pieces]\label{lem:pieces}
Let $V\ge W_{\rm prod}$ and let $N_k$ be the number of pieces of level $k$ in an initial front. Then $N_1\le1+2K_0$ and
$N_{k+1}\le(1+2C_\theta K)N_k$ for $1\le k<M$.
\end{lem}

\begin{proof}
Let $P$ be a piece of level $k$ and let $G_P$ be the set of $p\in P$ whose $(k+1)$-th contact is tangential. By stability
(Lemma~\ref{lem:stability}), $P\setminus G_P$ is open, and its connected components are the pieces of level $k+1$
contained in $P$. For $p\in G_P$ the tangential contact is on an obstacle within distance $\ell_H$ of $Q_k$ (of the initial segment if
$k=0$), on one of the two sides. By Lemmas~\ref{lem:tan0} and~\ref{lem:tank}, $G_P$ has at most $2K_0$ elements if $k=0$ and at most
$2C_\theta K$ if $k\ge1$. Removing $m$ points from an interval leaves at most $m+1$ intervals. $\qedhere$
\end{proof}

\begin{proof}[Proof of Proposition~\ref{prop:product}]
Let $I_0$ be an initial front. By Lemma~\ref{lem:pieces}, all $p\in I_0$ except finitely many belong to a piece of level $M$. Let $P$ be
such a piece. By Lemma~\ref{lem:cone}(iii), for $p\in P$,
\[
\prod_{k=1}^M\frac1{\cos\varphi_k}\le c_*^{1-M}\frac{|J_M|}{|J_1|}\cdot\frac1{\cos\varphi_M}\le2c_*^{1-M}\frac{|J_M|}{\cos\varphi_M},
\]
and by Lemma~\ref{lem:arcs} the integral over $P$ is at most $2c_*^{1-M}\cdot2\pi/\kappa_{\min}$. By Lemma~\ref{lem:pieces},
\[
\int_{I_0}\prod_{k=1}^M\frac{dp}{\cos\varphi_k}\le\frac{4\pi}{\kappa_{\min}}(1+2K_0)\,(1+2C_\theta K)^{M-1}c_*^{1-M} .
\]
Integrating over $a$ (an interval of length at most $2\sqrt2$), $\theta\in[0,2\pi)$ and $\rho_0\in[V,V+1]$ with weight $\rho_0$, and
summing over the at most $N_\Sigma$ initial fronts, we obtain the claim with
$$\Lambda_P:=(1+2C_\theta K)/\min(1,c_*) \qquad \text{and} \qquad C_P:=16\sqrt2\,\pi^2N_\Sigma(1+2K_0)/\kappa_{\min}. \qedhere$$
\end{proof}

\subsection{Comparison for fast particles}\label{sec:fastcomp}

We keep the setting of Section~\ref{sec:lin}: $A_0,A_1\in\cA_T$ coincide on $[0,t_0]$, $A_\lambda=A_0+\lambda\,\delta A$,
$I=[t_0,t_0+h]\subset[0,T]$, $h\le1$, and $\varepsilon=\sup_I|\delta A|$. 
Let $W_{\rm cvx}$ be the threshold provided by
Lemma~\ref{lem:convex} below. We henceforth fix
\begin{equation}\label{eq:choice-W0}
W_0:=
\max\left\{
W_{\rm prod},
W_{\rm cvx},
|E_0|,
\frac{8\ell_H}{\lambda_-}
\right\}.
\end{equation}
In particular, $W_0\ge W_{\rm geo}$. The threshold $W_0$ will not be
enlarged further.\\
The proof of the following proposition occupies the rest of this section.
\begin{prop}[Comparison for fast particles]\label{prop:fastcomp}
Assume {\rm(FH)}. There are $h_f>0$ and $C_f$, depending only on $\Dd$, $\alpha$ and $C_0$, such that, if $h\le h_f$ and $\varrho$ is a
density on ${\bf D}\times\bR^2$ with $\varrho(x,w)\le C_0e^{\alpha c_T}e^{-\alpha|w|}$, then for every $t\in I$
\[
\int_{\{|w|\ge W_0\}}\bigl|\widetilde X^{A_1}_t(\zeta)-\widetilde X^{A_0}_t(\zeta)\bigr|\,\varrho(\zeta)\,d\zeta\le C_f\,h\,\varepsilon .
\]
\end{prop}
We prove the proposition by interpolating between the two histories through
\[
A_\lambda=A_0+\lambda(A_1-A_0),\qquad \lambda\in[0,1].
\]
At values of $\lambda$ for which the collision pattern is transversal,
the trajectory is differentiable with respect to $\lambda$, and the
variational equation gives a bound containing one factor
$(\cos\varphi_k)^{-1}$ for each contact
(Lemma~\ref{lem:deriv}).

The collision pattern may, however, change at exceptional values of
$\lambda$, corresponding to a tangency or to a contact at the observation
time. Lemmas~\ref{lem:convex} and~\ref{lem:cont} show that the position
remains continuous across such parameters, while Lemma~\ref{lem:exc}
shows that, for almost every initial state, the exceptional set is
countable. Finally, we recover the difference
between the two trajectories by integrating the derivative with respect
to $\lambda$, and this is formulated in Proposition~\ref{prop:ac}.


\subsubsection{The derivative bound}

\begin{lem}[Derivative bound]\label{lem:deriv}
Let $t\in I$, and let $(\lambda,\zeta)$ be a good pair with $\zeta$ fast, $|w|=\rho_0$. Let
$\varphi_1,\dots,\varphi_m$ be the angles of the contacts of the $A_\lambda$-trajectory of $\zeta$ in $(t_0,t)$. Then
\[
\bigl|\partial_\lambda\widetilde X^{A_\lambda}_t(\zeta)\bigr|\le3\lambda_+\rho_0\,h\,\varepsilon\,e^{\lambda_+\rho_0h}
\prod_{k=1}^m\frac{C_L}{\cos\varphi_k}.
\]
\end{lem}

\begin{proof}
For $t_0\le s\le s'\le t$ let $\mathcal G(s',s)$ be the propagator of the homogeneous equation along the trajectory: free
evolution between contacts and $\sL$ at the contacts in $(s,s']$. By Lemma~\ref{lem:var}, $\eta$ solves the linear equation
with source $S$, jumps $\sL$ and $\eta(t_0)=0$, so that $\eta(t)=\int_{t_0}^t\mathcal G(t,s)S(s)\,ds$ (both sides solve the
same equation with the same initial value). By Lemma~\ref{lem:norms}(a),(b), and since every factor
$C_L/\cos\varphi_k$ is at least $1$, 
$$\|\mathcal G(t,s)\|\le e^{\lambda_+\rho_0h}\prod_{k=1}^mC_L/\cos\varphi_k.$$
Consequently,
\[
\begin{aligned}
|\eta_x(t)|
&\le \int_{t_0}^t\|\mathcal G(t,s)\|\,\|S(s)\|\,ds\\
&\le
3\lambda_+\rho_0h\varepsilon\,
e^{\lambda_+\rho_0h}
\prod_{k=1}^m\frac{C_L}{\cos\varphi_k}.
\end{aligned}
\]
 where we used Lemma~\ref{lem:norms}(c) and $|\eta_x|\le\|\eta\|$, proving the claim.
\end{proof}

\subsection{Absolute continuity in the history}\label{subsec:ac}

Fix $t\in I$ and let $\cU:=\{(x,w):\ x\in\mathring{\bf D},\ |w|>W_0\}$. The trajectories of states in $\cU$ are those of
Definition~\ref{def:fast-traj}. For $\zeta\in\cU$ set
\[
F_\zeta(\lambda):=\widetilde X^{A_\lambda}_t(\zeta),\qquad\lambda\in[0,1].
\]

\begin{lem}[Convexity near a tangency]\label{lem:convex}
There is $W_{\rm cvx}\ge W_{\rm geo}$ such that, if $W_0\ge W_{\rm cvx}$, the following holds. Let $(\lambda_*,\zeta_*)\in[0,1]\times\cU$ and
let $\tau_*$ be the time of a tangential contact on $Q$ of the $A_{\lambda_*}$-trajectory of $\zeta_*$. There are
$\delta,c_2>0$ such that every free flight $x(\cdot)$ of a history $A\in\cA_T$ which is $\cC^1$-close enough to the
flight of $\zeta_*$ on $[\tau_*-\delta,\tau_*+\delta]$ has $\varrho:=d^s_Q\circ x$ with $\dot\varrho$ Lipschitz and
$\ddot\varrho\ge c_2$ a.e.\ on $[\tau_*-\delta,\tau_*+\delta]$.
\end{lem}

\begin{proof} We have
$$\ddot\varrho=D^2d^s_Q(x)[\dot x,\dot x]+\langle\nabla d^s_Q(x),\ddot x\rangle$$ a.e., with $\ddot x=-\dot A\,\dot x+E_0$.
At $\tau_*$, $\dot x$ is tangent to $\partial Q$ at $x(\tau_*)$ and $|\dot x|\ge\tfrac34\lambda_-W_0$ by Proposition~\ref{prop:fast}(a), and
$$D^2d^s_Q[\dot x,\dot x]=\kappa|\dot x|^2.$$ By continuity, for $\delta$ small and for
flights close to the given one, 
$$D^2d^s_Q(x)[\dot x,\dot x]\ge\tfrac{\kappa_{\min}}2|\dot x|^2 \qquad \text{and} \qquad |\langle\nabla d^s_Q(x),\dot x\rangle|\le|\dot x|,$$ so that, if $W_0$ is large enough,
$$\ddot\varrho\ge\tfrac{\kappa_{\min}}2|\dot x|^2-|E_0|\,|\dot x|-|E_0|\ge c_2>0. \qedhere$$ 
\end{proof}

\begin{lem}[Continuity]\label{lem:cont}
The map $(\lambda,\zeta)\mapsto\widetilde X^{A_\lambda}_t(\zeta)$ is continuous on $[0,1]\times\cU$. In
particular $F_\zeta$ is continuous on $[0,1]$ for every $\zeta\in\cU$.
\end{lem}

\begin{proof}
Let $(\lambda_n,\zeta_n)\to(\lambda,\zeta)$. We argue by induction on the contacts of the limit trajectory, which are
finitely many (Proposition~\ref{prop:fast}), as in Lemma~\ref{lem:stability}: if the data at the
beginning of a flight converge, the free flights converge in $\cC^1$ on compact time intervals, together with the
time derivatives at the initial times. Away from the contact times of the limit, the limit flight stays at positive
distance from the obstacles, and so do the approximating flights. Near a transversal contact the argument of Lemma~\ref{lem:stability}
applies. Let $\tau_*$ be a tangential contact time of the limit, on $Q$, and let $\varrho_n,\varrho$ be the signed
distances to $Q$ along the approximating flights and along the limit flight, so that $\varrho_n\to\varrho$ in
$\cC^1$ on $[\tau_*-\delta,\tau_*+\delta]$. By Lemma~\ref{lem:convex}, 
$$\varrho(\tau_*)=\dot\varrho(\tau_*)=0, \quad \ddot\varrho_n\ge c_2 \quad \text{a.e.}\quad \text{and}\quad \varrho_n(\tau_*\pm\delta)\to\varrho(\tau_*\pm\delta)>0.$$ If $\varrho_n>0$ on
$[\tau_*-\delta,\tau_*+\delta]$ there is no contact there, and the flight converges to the limit flight. Otherwise
$\varrho_n$ has a first zero $\tau_n\in(\tau_*-\delta,\tau_*+\delta)$; since $\min\varrho_n\to0$ and the minimum
points converge to $\tau_*$, we have $\tau_n\to\tau_*$ and $\dot\varrho_n(\tau_n)\to0$, that is
$\langle w_n^-,n\rangle\to0$. The outgoing velocity $w_n^+\in\{w_n^-,\cR w_n^-\}$ satisfies
$|w_n^+-w_n^-|\le2|\langle w_n^-,n\rangle|\to0$, so the new flight issued at $\tau_n$ converges to the limit flight
on $[\tau_n,\tau_*+\delta]$. Along it, the signed distance to $Q$ vanishes at $\tau_n$, has nonnegative derivative there,
and is strictly convex; hence there is no further contact with $Q$ on $(\tau_n,\tau_*+\delta]$, and no other obstacle is
met for $n$ large. Thus the states at time $\tau_*+\delta$ converge, the positions converge uniformly on
$[\tau_*-\delta,\tau_*+\delta]$, and the induction proceeds. Positions converge also when a contact of the limit
occurs at time $t$.
\end{proof}

Recalling definition \ref{def:good}, for $\zeta\in\cU$ let
\[
\cE(\zeta):=\{0,1\}\cup\{\lambda\in[0,1]:\ (\lambda,\zeta)\text{ is not good}\}.
\]
By Lemma~\ref{lem:var}, $\cE(\zeta)$ is closed and $F_\zeta$ is of class $\cC^1$ on $[0,1]\setminus\cE(\zeta)$, with
$F_\zeta'(\lambda)=\partial_\lambda\widetilde X^{A_\lambda}_t(\zeta)$.

\begin{lem}[Exceptional parameters]\label{lem:exc}
For almost every $\zeta\in\cU$ the set $\cE(\zeta)$ is countable.
\end{lem}

\begin{proof}
The idea of the proof is quite simple: every bad pair \((\lambda,\zeta)\) satisfies one scalar equation \(g(\lambda,\zeta)=0\).
Locally, each bad set is therefore contained in a four-dimensional hypersurface of the five-dimensional space \((\lambda,\zeta)\).
For almost every fixed \(\zeta\), the vertical line \([0,1]\times\{\zeta\}\) meets each such hypersurface only finitely many times locally.
Since there are countably many hypersurfaces, the exceptional set of \(\lambda\)'s is countable.\\

A pair $(\lambda,\zeta)\in(0,1)\times\cU$ is not good if and only if,
after finitely many transversal contacts, either the next contact is
tangential at some time in $(t_0,t]$, or it is transversal and occurs
exactly at time $t$. We show that, locally, each of these events is
contained in a $\cC^1$ hypersurface of
$(0,1)\times\cU$. Fix a non-good pair $(\lambda_*,\zeta_*)$, and let $\tau_*$ be the time
of its first tangential contact, or $\tau_*=t$ if all the contacts up to
time $t$ are transversal and the last one occurs at time $t$. Denote by
$Q$ the obstacle involved in this contact. Choose
$t'\in(t_0,\tau_*)$ after the preceding contact, or any
$t'\in(t_0,\tau_*)$ if there is no preceding contact.

By Lemma~\ref{lem:var}, after restricting to a neighbourhood $O$ of
$(\lambda_*,\zeta_*)$, the state
$
S(\lambda,\zeta)
$
at time $t'$ is of class $\cC^1$ in $(\lambda,\zeta)$, and
$D_\zeta S$ is invertible. Let
\[
x(r;\lambda,\zeta)
:=
\pi_xF^{r,t'}_{A_\lambda}\bigl(S(\lambda,\zeta)\bigr)
\]
be the corresponding free flight, continued through the obstacles.
Then $x$ is of class $\cC^1$ in its variables, and
\[
L(r;\lambda,\zeta)
:=
\pi_xDF^{r,t'}_{A_\lambda}\,D_\zeta S
\]
has rank $2$. Suppose first that the relevant event is a transversal contact with
$Q$ at time $t$. On a sufficiently small neighbourhood of
$(\lambda_*,\zeta_*)$, define
$
g(\lambda,\zeta)
:=
d^s_Q\bigl(x(t;\lambda,\zeta)\bigr).
$
Every nearby pair having the corresponding contact at time $t$ belongs
to $\{g=0\}$. Moreover, at such a pair,
$
D_\zeta g
=
n\bigl(x(t;\lambda,\zeta)\bigr)\,L(t;\lambda,\zeta)
\ne0,
$
because $L(t;\lambda,\zeta)$ has rank $2$.

Suppose instead that the relevant event is a tangential contact with
$Q$ at time $\tau_*$. By Lemma~\ref{lem:convex}, there are
$\delta>0$ and, after shrinking $O$, a unique point
$
r_m(\lambda,\zeta)\in(\tau_*-\delta,\tau_*+\delta)
$
at which
$
r\longmapsto
d^s_Q\bigl(x(r;\lambda,\zeta)\bigr)
$
attains its minimum. The map $r_m$ is continuous. With a slight abuse of notation, set
$
g(\lambda,\zeta)
:=
d^s_Q\bigl(x(r_m(\lambda,\zeta);\lambda,\zeta)\bigr).
$
The minimizing property shows that $g$ is of class $\cC^1$. Indeed, if
$q=(\lambda,\zeta)$ and $q'=(\lambda',\zeta')$, then
\[
\begin{aligned}
&d^s_Q\bigl(x(r_m(q');q')\bigr)
 -d^s_Q\bigl(x(r_m(q');q)\bigr)\le g(q')-g(q)\le
d^s_Q\bigl(x(r_m(q);q')\bigr)
 -d^s_Q\bigl(x(r_m(q);q)\bigr).
\end{aligned}
\]
Since $r_m$ is continuous, both bounds have the same first-order
expansion as $q'\to q$. Consequently,
\[
D_\zeta g
=
n\bigl(x(r_m(\lambda,\zeta);\lambda,\zeta)\bigr)
L(r_m(\lambda,\zeta);\lambda,\zeta)
\ne0.
\]
Every nearby pair having a tangential contact with $Q$ belongs to
$\{g=0\}$.

We have therefore shown that every non-good pair has a neighbourhood
in which the set of non-good pairs is contained in a $\cC^1$
hypersurface. There are only countably many possible contact indices
and obstacles, and the space $(0,1)\times\cU$ is second countable.
Hence the set of non-good pairs is contained in a countable union
$
\bigcup_{j\ge1}M_j
$
of $\cC^1$ hypersurfaces of $(0,1)\times\cU$.

It remains to examine their sections at fixed $\zeta$. Let
$
\pi:(\lambda,\zeta)\longmapsto\zeta.
$
For every compact subset $K$ of a coordinate neighbourhood of $M_j$,
the area formula gives
\[
\int_{\bR^4}
\#\bigl(M_j\cap K\cap\pi^{-1}(\zeta)\bigr)\,d\zeta
\le
\mathcal H^4(M_j\cap K)<\infty.
\]
Thus, for almost every $\zeta$, the fiber
$M_j\cap K\cap\pi^{-1}(\zeta)$ is finite. Exhausting the coordinate
neighbourhoods by countably many compact sets and then taking the
countable union over $j$, we conclude that, for almost every $\zeta$,
the set of $\lambda\in(0,1)$ for which $(\lambda,\zeta)$ is not good is
countable. Adding the two endpoints $0$ and $1$ proves the claim.
\end{proof}

\begin{prop}\label{prop:ac}
For almost every $\zeta\in\cU$ such that
\[
\int_0^1
\bigl|\partial_\lambda\widetilde X^{A_\lambda}_t(\zeta)\bigr|
\,d\lambda<\infty,
\]
one has
\[
\bigl|
\widetilde X^{A_1}_t(\zeta)
-\widetilde X^{A_0}_t(\zeta)
\bigr|
\le
\int_0^1
\bigl|\partial_\lambda\widetilde X^{A_\lambda}_t(\zeta)\bigr|
\,d\lambda .
\]
\end{prop}

\begin{proof}
Fix $\zeta$ for which the conclusion of Lemma~\ref{lem:exc} holds and
the integral above is finite. Set
\[
F(\lambda):=\widetilde X^{A_\lambda}_t(\zeta),
\qquad
\Gamma(\lambda):=
F(\lambda)-F(0)-\int_0^\lambda F'(s)\,ds,
\]
where $F'$ is defined arbitrarily on the null set $\cE(\zeta)$.
By Lemma~\ref{lem:cont}, $\Gamma$ is continuous. On every connected
component $(a,b)$ of $[0,1]\setminus\cE(\zeta)$, the function $F$ is of
class $\cC^1$ and $\Gamma'=0$. Hence $\Gamma$ is constant on $[a,b]$.

Since the endpoints of every such interval belong to $\cE(\zeta)$, it
follows that
\[
\Gamma([0,1])\subset\Gamma(\cE(\zeta)).
\]
The set on the right is countable, whereas $\Gamma([0,1])$ is connected
because $\Gamma$ is continuous. Therefore $\Gamma([0,1])$ consists of
one point. Since $\Gamma(0)=0$, we obtain
\[
F(1)-F(0)=\int_0^1F'(\lambda)\,d\lambda.
\]
Taking absolute values gives the claim.
\end{proof}

The function $(\lambda,\zeta)\mapsto\partial_\lambda\widetilde X^{A_\lambda}_t(\zeta)$ is continuous on the open set of good pairs,
whose complement in $(0,1)\times\cU$ is null (for each $\lambda$ its section is null by Proposition~\ref{lem:ae-regular},
 applied to the trajectories issued at time $t_0$, and by Fubini); hence it is measurable, and Tonelli's theorem
applies to it.

\subsection{Proof of Proposition~\ref{prop:fastcomp}}
For $j\ge0$ let $V_j:=W_0+j$, $S_j:=\{(x,w):\ x\in{\bf D},\ V_j\le|w|<V_j+1\}$ and
$m_j:=1+\lfloor\frac54\lambda_+(V_j+1)h/d_{\min}\rfloor$. By Proposition~\ref{prop:fast}, for $\zeta\in S_j$ and $\lambda\in[0,1]$ the
$A_\lambda$-trajectory of $\zeta$ has at most $m_j$ contacts in $I$. Let
\[
h_f:=\min\Bigl\{\frac{\lambda_-d_{\min}}{10\lambda_+\ell_H},\ \frac{\alpha\,d_{\min}}{2\lambda_+\bigl(d_{\min}+\frac54\ln(C_L\Lambda_P)\bigr)}\Bigr\}.
\]
Since $V_j\ge W_0\ge8\ell_H/\lambda_-$, the first condition gives $m_j\le1+\frac52\lambda_+V_jh/d_{\min}\le\frac{\lambda_-V_j}{4\ell_H}+1\le M_{V_j}$.
At a good pair $(\lambda,\zeta)$ with $\zeta\in S_j$, Lemma~\ref{lem:deriv}, extended to the first $m_j$ contacts (all factors are
at least $1$), gives
\[
\bigl|\partial_\lambda\widetilde X^{A_\lambda}_t(\zeta)\bigr|\le3\lambda_+(V_j+1)\,h\,\varepsilon\,e^{\lambda_+(V_j+1)h}\,C_L^{m_j}
\prod_{k=1}^{m_j}\frac1{\cos\varphi^\lambda_k(\zeta)} .
\]
The pairs which are not good form a null set. By Tonelli's theorem and Proposition~\ref{prop:product} applied to each
$A_\lambda\in\cA_T$,
\begin{multline*}
\int_{S_j}\varrho(\zeta)\int_0^1\bigl|\partial_\lambda\widetilde X^{A_\lambda}_t(\zeta)\bigr|\,d\lambda\,d\zeta\\
\le C_0e^{\alpha c_T}e^{-\alpha V_j}\cdot3\lambda_+(V_j+1)\,h\,\varepsilon\,e^{\lambda_+(V_j+1)h}\,C_PC_L^{m_j}\Lambda_P^{m_j}(V_j+1).
\end{multline*}
By the second condition on $h$,
\[
e^{\lambda_+(V_j+1)h}(C_L\Lambda_P)^{m_j}\le C_L\Lambda_P\,e^{\frac\alpha2(V_j+1)},
\]
so the right-hand side is at most $C\,h\,\varepsilon\,(V_j+1)^2e^{-\alpha V_j/2}$, whose sum over $j$ is at most $C_fh\varepsilon$. In particular
$\int_0^1|\partial_\lambda\widetilde X^{A_\lambda}_t(\zeta)|\,d\lambda<\infty$ for $\varrho$-a.e.\ $\zeta$ with $|w|\ge W_0$, and
Proposition~\ref{prop:ac} gives
$|\widetilde X^{A_1}_t-\widetilde X^{A_0}_t|\le\int_0^1|\partial_\lambda\widetilde X^{A_\lambda}_t|\,d\lambda$ for $\varrho$-a.e.\ such $\zeta$
(the sphere $|w|=W_0$ is null). \qed

\subsection{Proof of Theorem~\ref{thm:C-tail}}\label{subsec:proof-tail}

Let $W_0$ be as at the beginning of this section. Let $A_0,A_1$ be the histories of two Lagrangian solutions with initial datum $\mu_0$,
$t_0:=\sup\{t\in[0,T]:A_0=A_1\text{ on }[0,t]\}$, and assume by contradiction that $t_0<T$. For $h\in(0,1]$ let
$I:=[t_0,\min(t_0+h,T)]$ and $\varepsilon:=\sup_I|A_1-A_0|$. Since all the conditions on $h$ below are upper bounds, we may
replace $h$ by $T-t_0$ if $T-t_0<h$, and assume $I=[t_0,t_0+h]\subset[0,T]$. As in the proof of Theorem~\ref{thm:C-fast}, by
Lemma~\ref{lem:equivalent-fp}(ii), and since the two trajectories of a point coincide up to time $t_0$, \eqref{eq:basic} holds
for $t\in I$, with $\varrho$ the density of the common state at time $t_0$ (Lemma~\ref{lem:tail}). For almost every $\zeta$ with
$|w|>W_0$ the trajectories in \eqref{eq:basic} are those of Definition~\ref{def:fast-traj}.

\begin{lem}[Slow particles]\label{lem:slow}
There are $h_s>0$ and $C_s$, depending only on $\Dd$, $c_0$ and $C_0$, such that for $h\le h_s$ and $t\in I$
\[
\int_{\{|w|<W_0\}}\bigl|\widetilde X^1_t-\widetilde X^0_t\bigr|\,\varrho\,d\zeta\le C_s\,h\,\varepsilon .
\]
\end{lem}

\begin{proof}
Let $\zeta$ be in the support of $\varrho$, with $|w|<W_0$; the corresponding initial velocity satisfies $|v|\ge c_0$. Under both
histories, $|V_t|\ge c_1$ for $t\in[0,T]$ (the lower bound in Lemma~\ref{lem:speed-window}, whose proof uses only $|v|\ge c_0$),
and $|V_t|\le\lambda_+(W_0+\lambda_+|E_0|)=:L_s$ for $t\in I$. The proofs of Lemmas~\ref{lem:one-contact} and~\ref{lem:comparison},
and of the estimate concluding the proof of Theorem~\ref{thm:C-fast}, use the upper bound $|v|\le L$ only through the bound
$|V_t|\le L_1$ for $t$ in the window and the bound $|w|\le W$ at time $t_0$, and they use the density only through
$\|\varrho\|_\infty\le\|f_0\|_\infty$. With $L_1$ and $W$ replaced by $L_s$ and $W_0$ respectively, and with $\|\varrho\|_\infty\le C_0$
(Lemma~\ref{lem:tail}), they give the claim, the integral being restricted to $\{|w|<W_0\}$ and $h_s$ beign the constant $h_0$ of
Lemma~\ref{lem:one-contact}, computed with $L_s$ in place of $L_1$.
\end{proof}
We have at last all the ingredients to prove Theorem~\ref{thm:C-tail}.
\begin{proof}[Proof of Theorem~\ref{thm:C-tail}]
By Lemma~\ref{lem:tail}, $\varrho\le C_0e^{\alpha c_T}e^{-\alpha|w|}$, so Proposition~\ref{prop:fastcomp} applies and the sphere
$\{|w|=W_0\}$ is null. By Lemma~\ref{lem:slow}, Proposition~\ref{prop:fastcomp} and \eqref{eq:basic}, for $h\le\min(h_s,h_f)$ and every
$t\in I$, 
$$|A_1(t)-A_0(t)|\le|E_0|(C_s+C_f)\,h\,\varepsilon.$$ Taking the supremum over $I$, taking $\varepsilon\le|E_0|(C_s+C_f)h\varepsilon$, and
for $h<1/(|E_0|(C_s+C_f))$ we get $\varepsilon=0$, that is $A_0=A_1$ on $[0,t_0+h]$, which contradicts the definition of $t_0$. Since the
constants do not depend on $t_0$, the proof is concluded. $\qedhere$
\end{proof}

\begin{oss}\label{rem:nobound-comments} We conclude with two remarks on the assumptions used in this section.\\
(a) \emph{The role of the assumptions.} The fast regime \eqref{eq:fast-regime} and the lower bound $c_0$ are used only for the slow
particles (Lemma~\ref{lem:slow}); the fast particles are in the dispersing regime because their trajectories have curvature
$O(|E_0|/|w|^2)$. The tail \eqref{eq:tail} is used only to sum over the speed shells in Proposition~\ref{prop:fastcomp}; the precise
condition is given in Remark~\ref{rem:tail}. Finite horizon, which is not needed in Section~\ref{sec:fast}, enters through $\ell_H$:
the flights have bounded length, which is used in Proposition~\ref{prop:fast}(c), in Lemma~\ref{lem:cone}), in the count
of the obstacles that can be hit next (Lemma~\ref{lem:pieces}), and in the bound $M\le M_V$. In infinite horizon a front can be cut
into infinitely many pieces at one collision, and a different argument would be needed.

(b) \emph{The product lemma.} For $E_0=0$, Proposition~\ref{prop:product} says that
$\int\prod_{k=1}^M\cos^{-1}\varphi_k\,d\mu\le C\Lambda^M$ for the invariant measure $\mu=\cos\varphi\,dr\,d\varphi$ of the billiard map.
The proof is somwhow elementary: the product is dominated by the growth $|J_M|/|J_1|$ of the perpendicular Jacobi field along a front, and
$\int|J_M|/\cos\varphi_M\,dp$ is the length of the arc of $\partial Q_M$ hit by a piece of the front, which is bounded because along a
dispersing front $\varphi$ increases with the arclength at rate at least $\kappa_{\min}$. No distortion bound, homogeneity strips or
mixing is used and the exponential growth comes only from the number of pieces.

\end{oss}

\section{Open problems}\label{sec:further} We conclude the paper with a list of open problem, left for future investigations:\\

\emph{Non-constant fields.}
We chose a constant field to keep the length of the paper reasonable. However, the key idea (the
reduction to the integrated thermostat history) is not specific to this case. For a bounded field
$E=E(t,x)$, $\bZ^2$-periodic and of class $\cC^1$ in $x$, the change of variables $W=e^{A(t)}V$
transforms the equations between contacts into
\[
\dot X=e^{-A(t)}W,\qquad \dot W=e^{A(t)}E(t,X),
\]
which, as in the constant field case, depend on $A$ but not on $\dot A$, and whose vector field is still divergence free. Under the
normalisation of the energy, the self-consistency condition becomes
\[
A(t)=\int_\Omega\int_0^t\bigl\langle E(s,X_s^A(z)),V_s^A(z)\bigr\rangle\,ds\,d\mu_0(z),
\]
so that $A(t)$ is again the average work done by the external field. The inner integral is a
continuous function of the history and of the initial state at regular points, so the existence
and compactness arguments of the present paper should extend with minor modifications.
If the field is spatially homogeneous and absolutely continuous in time, then
\[
\int_0^t\langle E(s),V_s^A(z)\rangle\,ds
=\langle E(t),\Delta_t^A(z)\rangle-\int_0^t\langle\dot E(s),\Delta_s^A(z)\rangle\,ds,
\]
and if $E=-\nabla U$ with $\nabla U$ periodic 
the work equals $U(x)-U(\widetilde X_t^A(z))$. In both cases the fixed point involves positions only,
as in the present paper. For a homogeneous field $E(t)$ the free flights remain affine, and we expect
all the proofs, including those of Section~\ref{sec:nobound}, to carry over essentially verbatim. If the field is no longer affine, the linearised equations contain terms
involving the spatial derivatives of $E$, and the flights of a front are no longer translates of each
other, so extending the uniqueness estimates, and in particular Lemmas~\ref{lem:tan0} and~\ref{lem:fan},
would require more work to control these terms.\\

\emph{Larger classes of weak solutions.} Proposition~\ref{prop:weak-lagr} identifies weak and Lagrangian
solutions for absolutely continuous measures. Whether uniqueness holds
among measure-valued weak solutions is a different question, which we do not address; note that the
empirical measures of the particle system are atomic weak solutions.\\

\emph{Rates and propagation of chaos.} Theorem~\ref{thm:B} is obtained by compactness and gives no rate. In
the fast regime, the comparison argument of Section~\ref{sec:fast} is quantitative, and combined with a
law of large numbers uniform over the compact set of histories it should give a rate of convergence for
independent initial data.\\

\emph{Long-time behaviour.} The stationary solutions of the kinetic equation and their stability are open;
they should be related to the nonequilibrium steady states of the particle system.\\

\emph{Self-consistent transfer operators.} The relation with globally coupled systems and self-consistent
transfer operators, and a possible route to the long-time behaviour, are interesting enough to deserve their own appendix (Appendix~\ref{app:STO}).

\appendix

\section{Proof of Lemma~\ref{lem:invariant-test}}
\label{app:transported-tests}

This appendix proves Lemma~\ref{lem:invariant-test}. The main issue is the
construction of local inverse branches of the billiard flow near a
transversal contact. The incoming and outgoing contact states are considered
separately and uniform transversality ensures that the corresponding
neighbourhoods can be chosen disjoint.

\begin{proof}
Let $K:=\operatorname{supp}\Theta$, a compact subset of the open set $\Reg_A$. The function is well defined
because $\Phi^t_A$ is injective on $\Reg_A$, by the same argument of the proof in Lemma~\ref{lem:F1}.

\emph{Step 0: uniformity.} By Lemma~\ref{lem:stability} the number of contacts is locally constant on
$\Reg_A$, and the contact times and the incoming states depend continuously on the initial point. Hence there
is an open set $U$ with $K\subset U$ and $\overline U\subset\Reg_A$ compact, and constants $\nu_K,\sigma_K>0$,
such that for every $z'\in U$: all contacts have normal velocity at least $\nu_K$; the first contact time is at
least $\sigma_K$, the last one at most $T-\sigma_K$, and two consecutive contact times differ by at least
$\sigma_K$. We regard $\Theta$ as a $\mathcal C^1$ function on $U$, vanishing outside
$K$.

\emph{Step 1: local inverse branches.} We use the variables $(x,w)$, $w=e^{A(t)}v$, in which the free-flight
maps $F^{t,s}_A$ are of class $\cC^1$ both in time and state. Let $z_0\in U$ and $t_*\in[0,T]$.
We construct a relatively open neighbourhood $\mathcal N$ of the state $(t_*,\Phi^{t_*}_A(z_0))$ in
$[0,T]\times\Omega$ and a map $\widehat Z(t,x,w)$ with values in $U$, of class $\cC^1$ in $(t,x,w)$ on an open
subset of $\bR\times\bR^4$, such that
\begin{equation}\label{eq:branch}
Z(t,(x,v)):=\widehat Z\bigl(t,x,e^{A(t)}v\bigr)\quad\text{satisfies}\quad
z=\Phi^t_A\bigl(Z(t,z)\bigr)\quad\text{for all }(t,z)\in\mathcal N .
\end{equation}
In particular $Z$ is of class $\cC^1$ in $z$ up to the boundary of $\Omega$, $D_zZ$ is continuous in $(t,z)$, and
$Z$ is Lipschitz in $t$ and differentiable in $t$ at every time at which $A$ is differentiable. We distinguish three cases:

\emph{(I) The time $t_*$ is not a contact time of $z_0$.}
Near $(t_*,z_0)$, the map
$
(t,z')\longmapsto \bigl(t,\Psi^t(z')\bigr)
$
is of class $\cC^1$. Indeed, after restricting to a sufficiently small
neighbourhood, the collision pattern is fixed and $\Psi^t$ is a composition
of free flights, reflections, and $\cC^1$ contact times. Moreover,
$D_{z'}\Psi^{t_*}(z_0)$ is invertible by Lemma~\ref{lem:F1}. The inverse
function theorem therefore gives a map
$
\widehat Z_{\mathrm{nc}}(t,x,w)
$
defined near $(t_*,\Psi^{t_*}(z_0))$ such that
$
\Psi^t\bigl(\widehat Z_{\mathrm{nc}}(t,x,w)\bigr)=(x,w).
$
This gives the desired local inverse branch at every noncontact state.

\emph{(II) The time $t_*$ is a contact time; incoming state $(y,w^-)$.}
Set
$
q:=t_*-\frac{\sigma_K}{4}.
$
By the choice of $\sigma_K$, the trajectories issued from points sufficiently
close to $z_0$ have no contact near time $q$. By Lemma~\ref{lem:F1}, after
shrinking the neighbourhood of $z_0$ if necessary, there are open sets
$U_q\ni z_0$ and $V_q\ni\Psi^q(z_0)$ such that
$
\Psi^q:U_q\longrightarrow V_q
$
is a $\cC^1$ diffeomorphism. Denote its inverse by
$
B_q:=\bigl(\Psi^q|_{U_q}\bigr)^{-1}.
$
For $(t,x,w)$ close to $(t_*,y,w^-)$, define
\[
\widehat Z_{\mathrm{in}}(t,x,w)
   :=B_q\bigl(F_A^{q,t}(x,w)\bigr).
\]
This is a $\cC^1$ map on a full neighbourhood of
$(t_*,y,w^-)$ in $\bR\times\bR^4$.\\
If $(x,e^{-A(t)}w)\in\Omega$ and $(t,x,w)$ is sufficiently close to
$(t_*,y,w^-)$, the free flight from time $q$ to time $t$ does not meet the
boundary before time $t$. Indeed, it is close to the corresponding free
flight of $z_0$, which remains in $\mathring{\bf D}$ on $[q,t_*)$, and the
signed distance to the obstacle is strictly decreasing near $t_*$ by
transversality. Since $x\in{\bf D}$, the first zero of this signed distance
cannot occur before $t$. Consequently,
$
\Psi^t\bigl(\widehat Z_{\mathrm{in}}(t,x,w)\bigr)=(x,w)$. This gives \eqref{eq:branch}.

\emph{(III) The time $t_*$ is a contact time; outgoing state $(y,w^+)$.}
For $(t,x,w)$ close to $(t_*,y,w^+)$, consider the free flight $ s\longmapsto F_A^{s,t}(x,w)$ run backwards from time $t$. By transversality and the implicit function
theorem, this flight meets the boundary of the relevant obstacle at a unique
time $\tau=\tau(t,x,w)$ close to $t_*$. The function $\tau$ is of class $\cC^1$. Write $ (y',w'):=F_A^{\tau(t,x,w),t}(x,w).$
Thus $(y',w')$ is the outgoing state at the backward crossing, and
$(y',\cR_{y'}w')$ is the corresponding incoming state. Define
\[
\widehat Z_{\mathrm{out}}(t,x,w)
   :=
   \widehat Z_{\mathrm{in}}
   \bigl(
      \tau(t,x,w),\,y',\,\cR_{y'}w'
   \bigr).
\]
After shrinking the neighbourhood if necessary, the argument of
$\widehat Z_{\mathrm{in}}$ belongs to its domain. Hence
$\widehat Z_{\mathrm{out}}$ is of class $\cC^1$ and satisfies
$
\Psi^t\bigl(\widehat Z_{\mathrm{out}}(t,x,w)\bigr)=(x,w)
$
for physical states sufficiently close to the outgoing state.

Since reflection is an involution, for $y'\in\partial{\bf D}$ and $w$ close
to $w^+$ one has
\[
\widehat Z_{\mathrm{out}}(t,y',w)
   =
\widehat Z_{\mathrm{in}}
   \bigl(t,y',\cR_{y'}w\bigr).
\]
Moreover,
\[
|w^+-w^-|
   =2|\langle w^-,n(y)\rangle|
   \ge 2\nu_K,
\]
so the incoming and outgoing neighbourhoods may be chosen disjoint.
Consequently, no compatibility between their derivatives is required. 

\emph{Step 2: regularity of $\varphi_\Theta$.} Let $\mathcal K$ be the set of $(t,z)$ such that
$z=\Phi^t_A(z_0)$ for some $z_0\in K$. By continuity of the branches and compactness of $[0,T]\times K$,
$\mathcal K$ is compact, and it is covered by finitely many neighbourhoods $\mathcal N_i$ as in Step~1, with
maps $Z_i$. We claim that
\[
\varphi_\Theta=\Theta\circ Z_i\quad\text{on }\mathcal N_i .
\]
Indeed, let $(t,z)\in\mathcal N_i$. If $Z_i(t,z)\in K$ this is the definition of $\varphi_\Theta$, by
\eqref{eq:branch}. If $Z_i(t,z)\notin K$, then $\Theta(Z_i(t,z))=0$; and $z=\Phi^t_A(z_0)$ with $z_0\in K$ is
impossible, since by \eqref{eq:branch} and the injectivity of $\Phi^t_A$ on $\Reg_A$ it would give
$z_0=Z_i(t,z)$. Thus $\varphi_\Theta=0$ as well. Since $\Theta\in\cC^1(U)$, the claim shows that $\varphi_\Theta$
has the stated regularity on $\bigcup_i\mathcal N_i$, an open neighbourhood of $\mathcal K$; and $\varphi_\Theta$
vanishes on the open complement of $\mathcal K$. Hence $\varphi_\Theta$ is continuous on $[0,T]\times\Omega$,
supported in $\mathcal K$, of class $\cC^1$ in $z$ with $\nabla_z\varphi_\Theta$ continuous, and Lipschitz in
$t$. The symmetry follows from the relation between the branches (II) and (III).

\emph{Step 3: the transport identity.} A point of $\mathcal K$ with $z\in\mathring\Omega$ is of type (I).
There $Z(t,\Phi^t_A(z'))=z'$, and between contacts $t\mapsto\Phi^t_A(z')$ solves $\dot x=v$,
$\dot v=E_0-\dot A\,v$ at every $t$ at which $A$ is differentiable. Since $Z$ is differentiable in $t$ at such
times and $\cC^1$ in $z$, differentiating in $t$ gives
$\partial_tZ+D_zZ\,(v,E_0-\dot Av)=0$ at such times, hence \eqref{eq:transported} for $\varphi_\Theta=\Theta\circ Z$.
Outside $\mathcal K$ the identity is trivial.
\end{proof}

\section{Globally coupled flows and self-consistent transfer operators}\label{app:STO}

The particle system studied in this paper is a globally coupled system, i.e. the particles evolve independently
except for a coupling through a global observable. For globally coupled {maps} the mean-field limit is
described by a self-consistent (nonlinear) transfer operator (see the review \cite{TanziReview} and, for
recent results on the existence, stability and convergence to equilibrium of its fixed points,
\cite{CasGaTa,CasGaTa1,BaLiSa}.) For {flows}, and in particular for flows with singularities such as
billiards, we are not aware of any result of this type; the smooth case is treated in \cite{Chau}, where
nonlinear Vlasov equations driven by smooth chaotic flows are studied with anisotropic spaces and microlocal
methods, which are not suitable in the presence of discontinuities. In this appendix we rephrase the results of the paper in this language and describe how they could
enter a study of the long-time behaviour. Except where a result of the paper is quoted, the discussion is
formal and nothing in this appendix is used elsewhere.

\paragraph{\bf The self-consistent semigroup}
For a history $A\in\cA_T$ let
\[
\mathcal L^t_A\mu:=(\Phi^t_A)_\#\mu
\]
be the transfer operator of the non-autonomous billiard flow of Section~\ref{sec:flow}; on densities,
$$\mathcal L^t_Af=e^{2A(t)}f\circ(\Phi^t_A)^{-1}$$ by Proposition~\ref{lem:ae-regular}. A Lagrangian solution is a curve
$\mu_t=\mathcal L^t_A\mu_0$ whose history is the one it generates, $\dot A=j_{\mu_t}/u_{\mu_t}$. In other
words it is a fixed point of the map $\{\mu_t\}\mapsto\{\mathcal L^t_{A[\mu]}\mu_0\}$ acting on curves of
measures. The coupling is of rank one in the sense that $\mathcal L_A$ depends on the curve only through the scalar function
$A[\mu]$, and this is why the fixed point problem can be written in $\cC^0([0,T])$, as we did in
Section~\ref{sec:lagr}. Where the solution is unique we may define
\[
\mathfrak L_t\mu_0:=\mu_t .
\]
Since the equation is autonomous, $\mathfrak L_{t+s}=\mathfrak L_t\circ\mathfrak L_s$ whenever all terms are
defined. This is the self-consistent transfer operator of the problem. Its infinitesimal form is
$$\partial_t\mu=\mathcal G_{a(\mu)}\mu, \quad \text{where}\quad \mathcal G_a\mu:=-\operatorname{div}\bigl((v,E_0-av)\mu\bigr),$$
with specular reflection, is the generator of the billiard flow with constant friction $a$, and
$a(\mu):=j_\mu/u_\mu$. A genuine self-consistent operator in the ``classical" sense (see e.g. \cite{TanziReview}), with an explicit
map, is obtained by freezing the friction over a time step $\tau$, and reads as
$$\mathfrak F_\tau\mu:=(\Psi^\tau_{a(\mu)})_\#\mu.$$
Here, $\Psi^\tau_a$ is the time-$\tau$ map of the flow with
constant friction $a$. One thus expects $\mathfrak F_\tau^{\,n}\to\mathfrak L_t$ as $\tau\to0$, $n\tau\to t$.

In these terms the results of the paper read as follows.
\begin{itemize}
\item By Theorem~\ref{thm:C-fast} and Corollary~\ref{cor:weak-unique}, $\mathfrak L_t$ is well defined on bounded
densities supported in $\{c_0\le|v|\le L\}$, for $t\le T$ in the fast regime, and, if the table has finite horizon, on
bounded densities supported in $\{|v|\ge c_0\}$ with an exponential tail (Theorem~\ref{thm:C-tail}). It is \emph{not} known to be a
semigroup for all times, because the fast regime is not preserved by the dynamics.
\item By Lemma~\ref{lem:empirical}, the particle system is the same evolution on empirical measures:
$\delta_{N,t}$ is a weak solution with datum $\delta_{N,0}$. This is the exact analogue of the relation between
$N$ globally coupled maps and their self-consistent transfer operator.
\item Theorem~\ref{thm:B} gives sequential stability along admissible
empirical approximations on the energy shell: every subsequential limit is
a Lagrangian solution. At initial data for which the Lagrangian solution is
unique, and in particular in the fast regime, this becomes
\[
\delta_{N,t}\rightharpoonup\mathfrak L_t\mu_0
\qquad\text{for every }t\in[0,T].
\]
\item The comparison estimate of Section~\ref{sec:fast} is a weak--strong Lipschitz bound for the dependence
of the linear operator on the coupling: for histories $A_0,A_1$ which agree up to the beginning of a short
window $I$, and $\mu_0$ with bounded density,
\[
\Bigl|\int\langle E_0,x\rangle\,d\bigl(\mathcal L^t_{A_1}-\mathcal L^t_{A_0}\bigr)\mu_0\Bigr|
\le C\,\|f_0\|_\infty\,|I|\,\sup_I|A_1-A_0| ,\qquad t\in I,
\]
with $x$ understood in the cover. Estimates of the form
$\|(\mathcal L_\mu-\mathcal L_\nu)f\|_{\mathrm w}\le C\|\mu-\nu\|_{\mathrm w}\|f\|_{\mathrm s}$, for a weak and a
strong norm, are the basic hypothesis in the theory of self-consistent transfer operators. Here the strong
space is $L^\infty$ with compactly supported velocities, the weak ``norm'' is a single observable, and the
estimate is obtained without any hyperbolicity, from the fact that the singularity of one collision is
integrable against the flux measure.
\end{itemize}

\paragraph{\bf Stationary states}
The fixed points of $\mathfrak L_t$ are the stationary solutions of the kinetic equation. They have a simple
description. If $\mu$ is stationary then $a:=a(\mu)$ is constant, and $\mu$ is invariant for the billiard flow
with constant friction $a$,
\begin{equation}\label{eq:const-friction}
\dot x=v,\qquad\dot v=E_0-a\,v,\qquad\text{with specular reflection.}
\end{equation}
Conversely, let $\mu$ be an invariant measure of \eqref{eq:const-friction}, for some $a\in\bR$, with finite
positive energy and such that the flow is defined $\mu$-almost everywhere. The energy identity
$\tfrac{d}{dt}|v|^2=2\langle E_0,v\rangle-2a|v|^2$, which is not affected by the reflections, gives
$0=2j_\mu-2a\,u_\mu$ after integration. Hence $a=a(\mu)$ automatically: {every invariant measure of
\eqref{eq:const-friction} is a stationary solution of the kinetic equation with its own energy}, and the
normalisation $u_\mu=1$ is the only condition left. If, for each $a$ in a suitable range, $\mu_a$ denotes a
physical invariant measure of \eqref{eq:const-friction}, the stationary states with unit
energy correspond to the solutions of the scalar equation
\begin{equation}\label{eq:scalar-stationary}
u(\mu_a)=1 .
\end{equation}
One expects $a>0$: since $\tfrac{d}{dt}|v|\le|E_0|-a|v|$, for $a>0$ the region $\{|v|\le|E_0|/a\}$ is absorbing,
so that $u(\mu_a)\le|E_0|^2/a^2$ and necessarily $a\le|E_0|$. The flow \eqref{eq:const-friction} contracts the
phase volume at the constant rate $2a$, and it is conformally symplectic (compare \cite{LiMi}); its physical
invariant measures are therefore singular. Thus, as for self-consistent transfer operators of maps, the
problem of the stationary states splits into a \emph{linear} problem (existence, uniqueness and regularity in
$a$ of the physical measure of a dissipative billiard flow) and a \emph{scalar} self-consistency equation.
These stationary states are the natural candidates for the limit, as $N\to\infty$, of the nonequilibrium
steady states studied in \cite{BoDaLeRi,BoCh1,BoCh2,BoCh3} and, for a stochastic collision mechanism, in
\cite{BoCAEsLe}.

\paragraph{\bf Linearisation: a rank-one perturbation}
Let $\mu_*$ be a stationary state with friction $a_*$ and unit energy, and write
$\mathcal L^t_*:=e^{t\mathcal G_{a_*}}$ for the transfer operator of \eqref{eq:const-friction} with $a=a_*$.
Formally, the differential of $\mu\mapsto a(\mu)$ at $\mu_*$ is the linear functional
\[
\ell(\nu):=\int\bigl(\langle E_0,v\rangle-a_*|v|^2\bigr)\,d\nu ,
\]
and $\partial_a\mathcal G_a\mu=\operatorname{div}_v(v\mu)=:\mathcal P\mu$. A perturbation $\nu_t=\mu_t-\mu_*$
therefore satisfies, to first order,
\[
\partial_t\nu_t=\mathcal G_{a_*}\nu_t+\ell(\nu_t)\,\mathcal P\mu_* ,
\]
so that the differential of the self-consistent semigroup at $\mu_*$ is the linear transfer operator plus a
perturbation of rank one, exactly as for the differential of a self-consistent transfer operator of maps
\cite{CasGaTa}. By Duhamel's formula the scalar $\alpha(t):=\ell(\nu_t)$, which is the perturbation of the
friction, solves the Volterra equation
\begin{equation}\label{eq:volterra}
\alpha(t)=\ell\bigl(\mathcal L^t_*\nu_0\bigr)+\int_0^tK(t-s)\,\alpha(s)\,ds,\qquad
K(t):=\ell\bigl(\mathcal L^t_*\,\mathcal P\mu_*\bigr).
\end{equation}
The kernel $K$ satisfies formally
\[
\left.
\frac{d}{da}\ell\bigl(\mathcal L_a^t\mu_*\bigr)
\right|_{a=a_*}
=\int_0^tK(s)\,ds .
\]
If $\mathcal L^t_*$ has a spectral gap on a
Banach space containing $\nu_0$ and $\mathcal P\mu_*$, then $K$ and the forcing term in \eqref{eq:volterra}
decay exponentially, and the linear stability of $\mu_*$ is decided by the scalar condition
\[
1-\widehat K(\lambda)\ne0\qquad\text{for }\operatorname{Re}\lambda\ge0 ,
\]
where $\widehat K$ is the Laplace transform of $K$. This is the analogue, for the present rank-one coupling,
of the stability conditions for stationary solutions of Vlasov equations, and of the conditions on the
differential in \cite{CasGaTa,CasGaTa1}. In the smooth setting a programme of this type is carried out in
\cite{Chau}.

The link with the present paper is the following. The kernel $K$ is the derivative with respect to the
friction of $t\mapsto\ell(\mathcal L^t_a\mu_*)$, and the comparison Lemma~\ref{lem:comparison} is a
finite-time statement of exactly this kind: the response of an observable to a change of the history is
Lipschitz, uniformly over the histories, when the measure has a bounded density. The singular factor
$1/\nu$ of a grazing collision, which makes the derivative of the flow unbounded, is integrable in the mean.

\paragraph{\bf What is missing}
To turn this into a theory of the long-time behaviour one needs, in our opinion, the following.

\emph{(1) A functional setting for the linear flow.} The density of $\mathcal L^t_*f$ grows like $e^{2a_*t}$,
and $\mu_*$ is singular: bounded densities, which are the strong space of this paper, are not preserved
uniformly in time. The natural replacement is an anisotropic Banach space of distributions which are regular
along unstable curves, as constructed for the finite horizon Sinai billiard flow in \cite{BaDeLi}. Here the
flow \eqref{eq:const-friction} differs from a billiard flow because it is dissipative, the speed is
not constant (so that the phase space is four-dimensional and non-compact with two neutral directions) and
the field bends the trajectories. Nevertheless, n the variables of Section~\ref{sec:flow} the field does not appear in the
linearised free flight, and enters the linearised collision map through a term whose size, relative to the
curvature term, is of order $|E_0|/(\kappa|v|^2)$. For fast particles the tangent dynamics is therefore a
small perturbation of that of the Sinai billiard, and one may hope to adapt the construction.

\emph{(2) Slow particles.} For $|v|^2\lesssim|E_0|/\kappa_{\min}$ the dynamics is not dispersing, and particles
may hop along the obstacles (Remark~\ref{rem:global-uniqueness}). In a stationary state all speeds up to $|E_0|/a_*$ are
allowed, so that slow particles cannot be excluded by an assumption on the support, as we did for finite
times. This is the main obstruction, both to global uniqueness and to a spectral theory.

\emph{(3) Regularity in the friction.} Equation \eqref{eq:scalar-stationary} and the kernel $K$ require the
differentiability of $a\mapsto\mu_a$ and of $a\mapsto\mathcal L^t_a$, i.e.\ linear response for the dissipative
billiard flow. The perturbation $\mathcal P\mu_*=\operatorname{div}_v(v\mu_*)$ is a derivative of a singular
measure in a direction which is not an unstable one, so that $\mathcal P\mu_*$ need not belong to the space
on which $\mathcal L^t_*$ has a spectral gap; as usual in linear response theory, one expects that only the
pairing with regular observables makes sense. The weak--strong estimate of this paper should be compared with
this: its analogue in an anisotropic setting would have the regularity along unstable curves in place of the
bounded density, and the integrability of $1/\nu$ along unstable curves, which is the content of the growth
lemma for billiards, in place of the flux identity.

\emph{(4) The nonlinear step.} Given (1)--(3), the convergence to the stationary state for initial data
close to $\mu_*$ would follow from a fixed point or cone contraction argument for the self-consistent
operator, in the spirit of \cite{CasGaTa,CasGaTa1,BaLiSa}. The scalar structure should help since, by
\eqref{eq:volterra}, the nonlinear problem can be closed on the friction alone, and the measure is then
recovered from the non-autonomous linear flow, which is asymptotically autonomous if the friction converges.
Combining this with a mean-field limit uniform in time would give propagation of chaos towards the
nonequilibrium steady state.\\


\noindent \textbf{Acknowledgements}
The author is immensely grateful to Carlangelo Liverani for introducing him to the problem and for the many extremely helpful discussions on the subject. He also thanks Yann Chaubet, Théophile Dolmaire, Stefano Galatolo and Chiara Saffirio for useful discussions on the topic.
The author acknowledges the UMI Group “DinAmicI” ({www.dinamici.org}) and the INdAM group GNFM. Finally, much of this work was carried out thanks to the hospitality of Scuola Normale Superiore di Pisa and Università di Pisa, for which the author is grateful.

\end{document}